\documentclass[12pt, twoside]{article}
\usepackage{amsmath,amsthm,amssymb}
\usepackage{times}
\usepackage{enumerate}

\def\titlerunning#1{\gdef\titrun{#1}}
\makeatletter
\def\author#1{\gdef\autrun{\def\and{\unskip, }#1}\gdef\@author{#1}}
\def\address#1{{\def\and{\\\hspace*{18pt}}\renewcommand{\thefootnote}{}%
\footnote {#1}}%
\markboth{\autrun}{\titrun}}
\makeatother
\def\email#1{e-mail: #1}
\def\subjclass#1{{\renewcommand{\thefootnote}{}%
\footnote{\emph{Mathematics Subject Classification (2010):} #1}}}
\def\keywords#1{\par\medskip
\noindent\textbf{Keywords.} #1}

\theoremstyle{definition}

\numberwithin{equation}{section}

\usepackage{pb-diagram}
\usepackage{color}
\usepackage[normalem]{ulem}

\newcommand{\bb}[1]{\ensuremath{\mathbb{#1}}}

\newcommand{\Bad}{\operatorname{Bad}}

\newcommand{\es}{\ensuremath{\emptyset}}

 \DeclareMathOperator{\pol}{pol}
 \DeclareMathOperator{\fr}{fr}
 \DeclareMathOperator{\intr}{int}
  \DeclareMathOperator{\bdry}{bd}

 \DeclareMathOperator{\acl}{acl}
\DeclareMathOperator{\dcl}{dcl}

\DeclareMathOperator{\Stab}{Stab}
\DeclareMathOperator{\stab}{Stab}

 \DeclareMathOperator{\dom}{dom}

\DeclareMathOperator{\tp}{tp}

\DeclareMathOperator{\mr}{RM}
\DeclareMathOperator{\md}{MD}

\DeclareMathOperator{\cl}{cl}

\newtheorem{theorem}{Theorem}[section]

\newtheorem{claim}[theorem]{Claim}

\newtheorem{corollary}[theorem]{Corollary}

\newtheorem{fact}[theorem]{Fact}
\newtheorem{lemma}[theorem]{Lemma}

\newtheorem{proposition}[theorem]{Proposition}

\theoremstyle{definition}
\newtheorem{defn}[theorem]{Definition}
\newtheorem{definition}[theorem]{Definition}
\newtheorem{example}[theorem]{Example}
\newtheorem{remark}[theorem]{Remark}
\newcommand{\CL}{{\mathcal L}}

\newcommand{\CN}{{\mathcal N}}
\newcommand{\CH}{{\mathcal H}}
\newcommand{\CR}{{\mathcal R}}
\newcommand{\CM}{{\mathcal M}}

\newcommand{\CC}{{\mathcal C}}

\newcommand{\CF}{{\mathcal F}}
\newcommand{\CD}{{\mathcal D}}

\newcommand{\CS}{\mathcal S}

\newcommand{\0}{\emptyset}

\newcommand{\UB}[1]{{#1}_{\mbox{\tiny {$\pol$}} }}

\newcommand{\sub}{\subseteq}
\usepackage[matrix, arrow,all,cmtip,color]{xy}
\newcommand{\groupconf}[6]{\xy <1.5cm,0cm>: (1,0.2)*{#4}; (2,0.2)*{#5}; (.9,-.9)*{#6};
    (-0.2,0)*{#1}; (-0.2,-1)*{#2}; (-0.2,-2)*{#3}; (1,0)*{}; (0,-2)*{}
    **\dir{-}; (0,-2)*{}; (0,0)*{} **\dir{-}; (0,0)*{}; (2,0)*{}
    **\dir{-}; (2,0)*{}; (0,-1)*{} **\dir{-};
    \endxy }
\newcommand{\groupconflarge}[6]{\xy <1.5cm,0cm>: (1,0.2)*{#4}; (2,0.2)*{#5}; (.9,-.9)*{#6};
    (-.6,0)*{#1}; (-0.6,-1)*{#2}; (-0.6,-2)*{#3}; (1,0)*{}; (0,-2)*{}
    **\dir{-}; (0,-2)*{}; (0,0)*{} **\dir{-}; (0,0)*{}; (2,0)*{}
    **\dir{-}; (2,0)*{}; (0,-1)*{} **\dir{-};
    \endxy }

\renewcommand{\phi}{\varphi}

\long\def\symbolfootnote[#1]#2{\begingroup%
\def\thefootnote{\fnsymbol{footnote}}\footnote[#1]{#2}\endgroup}

\makeatletter

\def\Ind#1#2{#1\setbox0=\hbox{$#1x$}\kern\wd0\hbox to 0pt{\hss$#1\mid$\hss}
\lower.9\ht0\hbox to 0pt{\hss$#1\smile$\hss}\kern\wd0}

\def\Notind#1#2{#1\setbox0=\hbox{$#1x$}\kern\wd0\hbox to 0pt{\mathchardef
\nn=12854\hss$#1\nn$\kern1.4\wd0\hss}\hbox to
0pt{\hss$#1\mid$\hss}\lower.9\ht0 \hbox to
0pt{\hss$#1\smile$\hss}\kern\wd0}

\def\la{\langle}
\def\ra{\rangle}

\def\op{\oplus}
\def\om{\ominus}
\def\mfR{\mathfrak R}
\def\sub{\subseteq}
\def\sm{\setminus}

\begin{document}

%%%%% To ease editing, add:

\baselineskip=17pt

%%%%%%%%%%%%%%%%

%% In the running head, give an abbreviation of the title.
\titlerunning{Strongly minimal groups in o-minimal structures}

\title{Strongly minimal groups in o-minimal structures}

\author{Pantelis E. Eleftheriou
\and
Assaf Hasson
\and
Ya'acov Peterzil}

\date{}

\maketitle

\address{P. E. Eleftheriou: Department of Mathematics and Statistics, University of Konstanz, Box 216, 78457 Konstanz, Germany; \email{panteleimon.eleftheriou@uni-konstanz.de}
\and
A. Hasson: Department of Mathematics, Ben Gurion University of the Negev, Be'er-Sheva 84105, Israel; \email{hassonas@math.bgu.ac.il}
\and
Y. Peterzil: Department of Mathematics, University of Haifa, Haifa, Israel; \email{kobi@math.haifa.ac.il}}


%%%%%%%%

\begin{abstract}
We prove Zilber's Trichotomy Conjecture for strongly minimal expansions of
$2$-dimensional groups,  definable in o-minimal structures:
\smallskip

\noindent{\bf Theorem.} Let $\CM$ be an o-minimal expansion of a real closed field,
$\la G;+\ra $ a 2-dimensional group definable in $\CM$, and $\CD=\la G;+, \ldots\ra$ a
strongly minimal structure, all of whose atomic relations are definable in $\CM$. If
$\CD$ is not locally modular, then an algebraically closed field $K$ is
interpretable in $\CD$, and the group $G$, with all its induced $\CD$-structure, is
definably isomorphic in $\CD$ to an algebraic $K$-group with all its induced
$K$-structure.

%% Keywords are optional
\end{abstract}

\section{Introduction}

\subsection{Zilber's Conjecture (ZC)}
In \cite{ZilTri}, Boris Zilber formulated the following conjecture.\\

%{\bf do we write Zil'ber or Zilber?}

\noindent {\bf Zilber's Trichotomy Conjecture.} {\em The geometry of every strongly
minimal structure $\CD$ is either (i) trivial, (ii) non-trivial and locally modular,
or (iii) isomorphic to the geometry of an algebraically closed field $K$ definable
in $\CD$. Moreover, in (iii)  the structure induced on $K$ from $\CD$ is already
definable in $K$ (that is, the field $K$ is ``pure'' in $\CD$)}.\smallskip

The conjecture reduces by \cite{Hr4} to: {\em if a strongly minimal structure $\CD$
is not locally modular, then it interprets a field $K$, and the field $K$ is pure in
$\CD$}.\smallskip

In the early 1990s, Hrushovski refuted both parts of the conjecture. Using his
amalgamation method he showed the existence of a strongly minimal structure which is
not locally modular and yet does not interpret any group (so certainly not a field),
see \cite{Hr1}. In addition he showed the existence of {\em a proper} strongly
minimal expansion of a field, see \cite{Hr2}, thus disproving also the purity of the
field. Nevertheless, Zilber's Conjecture stayed alive since it turned out to be true
in various restricted settings, and moreover its verification in those settings gave
rise to important applications (such as Hrushovski's proof of the function field
Mordell-Lang conjecture in all characteristics \cite{HrMordellLang}).

A common feature to many cases where the conjecture is true is the presence of an
underlying geometry putting strong restrictions on the definable sets in the
strongly minimal structure $\CD$. This is for example the case when $\CD$ is
definable in an algebraically closed field (\cite{Hasson-Sustretov},  \cite{MaPi}
and \cite{Ra}), in a differentially closed field (\cite{MarZariski}),
separably closed field (\cite{HrMordellLang}), or in algebraically closed valued field
(\cite{KowRand}). This is also the case when $\CD$ is endowed with a Zariski
geometry (\cite{HrZil}).

Thus, it is interesting to examine the conjecture in various geometric settings. In this paper, we consider Zilber's Conjecture in the o-minimal geometric setting, introduced in the
1980s (\cite{vdd-tarski, kps, ps}). O-minimality imposes strong conditions on definable complex analytic objects forcing them in many cases to be algebraic (see \cite{PetStaICM} for
a survey, and \cite{BeKliTsi} for a recent application). The results of this paper can be seen among others as another manifestation of the same phenomenon.

%This paper is devoted to the Zilber Trichotomy Conjecture for strongly minimal structures definable in an o-minimal one.

%In this paper, we settle a main case of the Zilber Trichotomy Conjecture for strongly minimal structures definable in an o-minimal one.

\subsection{The connection to o-minimality} The  complex field is an example of a strongly minimal structure definable in the o-minimal $\la \mathbb R;+,\cdot,<\ra$, and indeed, the
underlying Euclidean geometry is an important component in understanding complex
algebraic varieties. This leads to examining in greater generality those strongly
minimal structures definable in o-minimal ones, and  to the following restricted
variant of Zilber's Conjecture, formulated by the third author in a model theory
conference at East Anglia in 2005.
\\

\noindent{\bf The o-minimal ZC. }{\em Let $\CM$ be an o-minimal structure and $\CD$
a strongly minimal structure whose underlying set and atomic relations are definable
in $\CM$. If $\CD$ is not locally modular, then an algebraically closed field $K$ is
interpretable in $\CD$, and moreover, $K$ is a pure field in $\CD$.}
\\

\begin{remark}
\begin{enumerate}
\item Because every algebraically closed field of characteristic zero (ACF$_0$) is
definable in an o-minimal real closed field, Zilber's Conjecture for reducts of
algebraically closed fields of characteristic zero is a special case of the
o-minimal ZC. This variant of the conjecture is still open for reducts whose
universe is not an algebraic curve.

 \item The purity of the field in the o-minimal setting was already
proven in \cite{PeStExpansions}, thus the o-minimal ZC reduces to proving the
interpretability of a field in $\CD$.

\item Since every definable algebraically closed field in an o-minimal structure has
dimension $2$ (see \cite{PeSte}), it is not hard to see that the above conjecture
implies that the underlying universe of $\CD$ must be $2$-dimensional in $\CM$.
Therefore, it is natural to consider the o-minimal ZC under the $2$-dimensional
assumption on $\CD$, which is the case of our Theorem \ref{main} below.

\item By \cite{HaOnPe}, if $\CD$ is strongly minimal, interpretable in an o-minimal
structure and in addition $\dim_{\CM}D=1$, then $\CD$ must be locally modular, thus
trivially implying the o-minimal ZC in the case when $\dim_{\CM} D=1$.

\item The theory of  compact complex manifolds, denoted by CCM, (see
\cite{ZilZariski}) is the multi-sorted theory of the structure whose sorts are all
compact complex manifolds, endowed with all analytic subsets and analytic maps. It
is known (\cite[Theorems 3.4.3 and 3.2.8]{ZilZariski}) that each sort in this
structure has finite Morley rank, and also that the structure is interpretable in
the o-minimal $\mathbb R_{an}$. Hence, every sufficiently saturated structure
elementarily equivalent to a CCM is interpretable in an o-minimal structure.

By \cite{MooRE}, every set of Morley rank one in any model of CCM is definably
isomorphic to an algebraic curve. Thus, Zilber's conjecture for reducts of CCM whose
universe is analytically $1$-dimensional reduces to the work in
\cite{Hasson-Sustretov}. The higher dimensional cases may also reduce to the
conjecture for ACF$_0$ but this is still open.
\end{enumerate}
\end{remark}

In \cite{HaKo} the following case of the o-minimal ZC was proven.
\begin{theorem} Let $\CR:=\la R, +, \cdot, <, \dots \ra$ be an o-minimal expansion of a real closed field, $K:=R[i]$ its algebraic closure. Let $f: K\to K$ be an $\CR$-definable function  If $\CD=\la  K;+,f\ra$ is strongly
minimal and non-locally modular (equivalently, $f$ is not an affine map),  then up to
conjugation by an invertible $2\times 2$ $R$-matrix and finitely many corrections,
$f$ is a $K$-rational function. In addition, a function $\odot: K^2\to
K$ is definable in $\CD$, making $\la K;+,\odot\ra$ an algebraically
closed field.
\end{theorem}

In our current result below we replace the additive group of $K$ above by an
arbitrary $\CR$-definable $2$-dimensional group $G$.  Moreover,
we let $\CD$ be an arbitrary expansion of $G$ and not only by a map $f:G\to G$.
Since strongly minimal groups are abelian (\cite[Corollary 3.1]{PoiGroups}), we
write the group below additively. Here is the main theorem of our article.

\begin{theorem}\label{main} Let $\CM$ be an o-minimal expansion of a real closed field $R$, and
let $\la G;\oplus \ra $ be a $2$-dimensional group definable in $\CM$. Let $\CD=\la
G;\oplus,\ldots\ra$ be a strongly minimal structure expanding $G$, all of whose
atomic relations are definable in $\CM$.

Then there are in $\CD$ an interpretable algebraically closed field $K$, a
$K$-algebraic group $H$ with $\dim_K H=1$, and a definable isomorphism $\phi:G\to
H$, such that the definable sets in $\CD$ are precisely those of the form
$\phi^{-1}(X)$ for $X$ a $K$-constructible subset of $H^n$.

In fact, the structure $\CD$ and the field $K$ are bi-interpretable.
\end{theorem}

Note that the  theorem implies in particular that $G$ is definably isomorphic in
$\CD$ to either $\la K;+\ra$, $\la K^\times;\cdot\ra$ or to an elliptic curve over
$K$.

\subsection{The general strategy: from real geometry and strong minimality to complex algebraic geometry}\label{sec-overview} Let $\CM$, $G$ and $\CD$ be as in Theorem \ref{main}.
Since $G$ is a group definable  in an o-minimal expansion of a real closed field
$R$, it admits a differentiable structure which makes it into a Lie group with
respect to $R$ (see \cite{Pi5}). We let $\mathfrak F$ be the collection of all
differentiable (with respect to that Lie structure) partial functions $f:G\to G$,
with $f(0_G)=0_G$, such that for some $\CD$-definable strongly minimal $S_f\sub
G^2$, we have $graph(f)\sub S_f$. We let $J_0 f$ denote the Jacobian matrix of $f$
at $0$. The following is easy to verify, using the chain rule for differentiable
functions:
$$J_0(f\oplus g)=J_0f+J_0 g\,\,\,;\,\,\, J_0(f\circ g)=J_0 f\cdot J_0 g,$$
where on the left hand side of each equation we use the group operation and functional composition, and on the
right hand side the usual matrix operations in $M_2(R)$. Let also
$$\mathfrak R=\{J_0f\in
M_{2}(R): f\in \mathfrak F\}.$$

The key observation, going back to Zilber, is that via the above equations we can
recover a ring structure on $\mathfrak R$ by performing addition and composition of
curves in $\CD$. Most importantly, for the ring structure to be $\CD$-definable,
%to do that $\CD$-definably, namely for $\mathfrak$ to become a $\CD$-definable ring,
 one needs to recognize tangency of curves at a point $\CD$-definably. The geometric idea
for that goes back to Rabinovich's work \cite{Ra}, and requires us to develop a
sufficient amount of intersection theory for $\CD$-definable sets, so as to
recognize ``combinatorially'' when two curves are tangent.

%As we shall see, our goal is to show that the combinatorial restrictions of strong
%minimality, together with the geometric restrictions of o-minimality make the
%geometry and Intersection Theory of definable sets in $\CD$ resemble those of
%complex analytic sets.

%\subsection{The structure of the proof}

This paper establishes in several distinct steps the necessary ingredients for the
proof. In each of these steps we prove an additional property of $\CD$-definable
sets which shows their resemblance to complex algebraic sets. We briefly describe
these steps.

We call  $S\sub G^2$ {\em a plane curve} if it is $\CD$-definable and $\mr(S)=1$ (we
recall the definition of Morley rank in Section \ref{sec-morley}). In Section
\ref{sec-frontier} we investigate the frontier of plane curves, where the frontier
of a set $S$ is $\cl(S)\setminus S$. We prove that every plane curve has finite
frontier in the group topology on $G$.

In Section \ref{sec-poles} we consider the poles of plane curves, where a pole of
$S\sub G^2$ is a point $a\in G$, such that for every neighborhood $U\ni a$, the set
$(U\times G)\cap S$ is ``unbounded''. We prove that every plane curve has at most
finitely many poles.

As a corollary of the above two results we establish in Section \ref{sec-topology}
another geometric property which is typically true for complex analytic curves.
Namely, we show that every plane curve $S$ whose projection on both coordinates is
finite-to-one, is locally, outside finitely many points, the graph of a
homeomorphism.

Next, we discuss the differential properties of plane curves, and consider in
Section \ref{sec-jacobians} the collection, $\mathfrak R$, of all Jacobian matrices at $0$ of local
smooth maps from $G$ to $G$ whose graph is contained in a plane curve. Using our
previous results we prove that this collection forms an algebraically closed
subfield $K$ of $M_2(R)$, and thus up to conjugation by a fixed invertible matrix,
every such Jacobian matrix at $0$ satisfies the Cauchy-Riemann equations.

In Section \ref{sec-intersection} we establish elements of complex intersection
theory, showing that if two plane curves $E$ and $X$ are tangent at some point, then
by varying $E$ within a sufficiently well-behaved family, we gain additional
intersection points with $X$. This allows us to identify tangency of curves in $\CD$
by counting intersection points.

Finally, in Section \ref{sec-last} we use the above results in order to interpret an
algebraically closed field in $\CD$ and prove our main theorem.\\

\section{Preliminaries} We review  briefly the basic model theoretic notions
appearing in the text. We refer to any standard textbook in model theory (such as
\cite[\S6, \S7]{MaBook}) for more details. Standard facts on o-minimality can be
found in \cite{vdDries} whose Sections 1.1 and 1.2 provide most of the basic
background needed on structures and definability.

\subsection{Strong minimality and related notions}

\label{sec-morley} Throughout the text, given a structure $\CN$, by
\emph{$\CN$-definable} we mean definable
 in $\CN$ with parameters, unless stated otherwise. We drop the index `$\CN$-' if it is clear from the context.
  In the next subsection, we will adopt a global convention about this index to be enforced in Sections
  \ref{sec-frontier} - \ref{sec-last}.

Let $\CN=\la N, \dots\ra$ be an $\omega$-saturated structure. A definable set $S$
is \emph{strongly minimal} if every definable subset of $S$ is  finite or co-finite.
We call $\CN$ \emph{strongly minimal} if $N$ is a strongly minimal set.

Let $\CN=\la N, <, \dots\ra$ be an expansion of a dense linear order without
endpoints. We call $\CN$ \emph{o-minimal} if every definable subset of $N$ is a
finite union of points from $N$ and open intervals whose endpoints lie in
$N\cup\{\pm \infty\}$. The standard topology in $\CN$ is the order topology on $N$
and the product topology on $N^n$.

Now let $\CN$ be a strongly minimal or an o-minimal structure. The algebraic closure
operator $\acl$ in both cases is known to give rise to a pregeometry. We refer to
\cite[\S 6.2]{MaBook} and \cite[\S 1]{Pi5} for all details, and recall here only
some.  Given $A\sub N$ and $a\in N^n$, we let $\dim(a/A)$ be the size of a maximal
$\acl$-independent subtuple of $a$ over $A$. Given a  set $C\sub N^n$, definable over $A$ we
let
$$\dim(C)=\max\{\dim(a/A):a\in C\},$$ and we call an element $a\in C$ {\em generic in $C$ over $A$ in $\CN$} if
$\dim(a/A)=\dim(C)$. We also note that $\CN$ eliminates the $\exists^{\infty}$ quantifier. Namely, if $\varphi(x,y)$ is a formula, then the set of all $x$ for which there are infinitely many $y$ such that  $\varphi(x,y)$ holds is a definable set. We say that $\exists^{\infty} y\varphi(x,y)$ defines that set.

If $\CN$ is a strongly minimal structure, then $\dim(C)$ coincides with the Morley
rank of $C$, and we denote $\dim(a/A)$ and $\dim(C)$ by $\mr(a/A)$ and $\mr(C)$,
respectively. We denote the Morley degree of $C$ by $\md(C)$.
 In the o-minimal case, $\dim C$ coincides with topological dimension
of $C$, and we keep the notation $\dim(a/A)$ and $\dim(C)$.

%It is defined recursively as follows: For  $C\subseteq D$, $\mr(C)=0$ if and only if
%$C$ is finite and otherwise $\mr(C)=1$.

%For $C\sub D^n$, $\mr(C)<n$ if and only if there exists a finite partition
%$C=\bigcup_{i=1}^k C_i$ into definable sets, and for each $i=1,\ldots, k$ there
%exist $d_i<n$ and coordinates $j_1,\ldots, j_{d_i}\in \{1,\ldots, n\}$, such that
%the projection of $C_i$ into these coordinates is finite-to-one. In this case
%$\mr(C)$ is the minimum of $\max\{d_i:i=1,\ldots, k\}$ when we consider all possible
%partitions of $C$ . If there is no such partition then $\mr(C)=n$. The Morley degree
%of $C$ is the maximal $d$, such that $C$ can be partitioned into definable sets, all
%of the same Morley rank as $C$.

%{\bf I don't like the above. But i could not find anything better. Is it clear that
%we need it?}

% if $k$ is maximal such that for some projection
%on $k$ of the coordinates, $\pi: M^n\to M^k$, the image $\pi(S)$ contains an open
%set.

Let $\CN$ be any structure. Given a definable set $X$, {\em a canonical parameter
for $X$} is an element in $\CN^{eq}$ which is inter-definable with the set $X$,
namely $\bar a$ is a canonical parameter
 for $X$ if $\phi(\bar x, \bar a)$ defines $X$ and $\phi(\bar x, \bar a')\neq X$ for all $\bar a'\neq \bar a$.
Any two canonical parameters are inter-definable over $\emptyset$, and so we use
$[X]$ to denote any such parameter. Note that if $X=X_{t_0}$ for some definable
family
 of sets over $\0$, $\{X_t:t\in T\}$, then $[X]\in \dcl(t_0)$, but $t_0$ need not be a canonical parameter for $X$.

% If $C$ is a definable set in a strongly
%minimal structure $\CD$ an element $c\in C$ is \emph{$\CD$-generic\footnote{This is
%a different notion than that of generic points in algebraic geometry.}} over a
%parameter set $A$ if there is no $A$-definable set $C'$ of Morley rank smaller than
%the Morley rank of $C$ such that $c\in C'$.

A structure $\CN=\la N, \dots\ra$ is \emph{interpretable} in  $\CM$ if there is an
isomorphism of structures $\alpha:\CN\to \CN'$, where the universe of $\CN'$ and all
$\CN'$-atomic relations are interpretable $\CM$.

If $\CN$ is interpretable in $\CM$ via $\alpha$ and $\CM$ is interpretable in $\CN$
via $\beta$, and if in addition $\beta\circ \alpha$ is definable in $\CN$ and
$\alpha\circ \beta$ is definable in $\CM$, then we say that {\em $\CM$ and $\CN$ are
bi-interpretable}.

Note that  if $\CM$ is an o-minimal expansion of an ordered group, then by Definable
Choice, every interpretable structure in $\CM$ is also definable in $\CM$.\\

\subsection{The setting}\label{sec-setting} Throughout Sections \ref{sec-frontier} - \ref{sec-last}, we
 fix a sufficiently saturated o-minimal expansion  $\CM=\la R;+,\cdot,<,\ldots \ra$ of a real closed
field.
As described in \cite[Chapters 6-7]{vdDries}, definable sets in $\CM$ admit various topological properties with respect to the underlying order topology on $R$ and the product topology on $R^n$. In addition, a theory of differentiability
with respect to $R$ is developed there, allowing notions which are analogous to classical ones, such as  manifolds, differentials of definable maps, jacobian matrices, etc. We are going to exploit this theory heavily, similarly to the way $\mathbb R$-differentiability  is often used when developing complex algebraic geometry.

Throughout the same sections, we also fix a $2$-dimensional $\CM$-definable group $G$.
By \cite{Pi5}, the group $G$ admits a definable $C^1$-manifold structure with respect to the
field $R$, such that  the group operation and inverse function are $C^1$ maps with
respect to it. The topology and differentiable structure which we refer to below is
always that of this smooth group structure on $G$. Note that the group $G$ is definably isomorphic, as a topological group, to a
definable group whose domain is a closed subset of some $R^r$, endowed with the
$R^r$-topology (see, for example, \cite[Claim 3.1]{PetStaStab}). Thus, we assume that $G$ is a closed subset  of $R^r$ and its
topology is the subspace topology.

Finally, throughout Sections \ref{sec-frontier} - \ref{sec-last}, we  fix a strongly
minimal non-locally modular structure
 $\CD=\la G;  \dots\ra$ definable in $\CM$. We treat $\CM$ as the default
structure and thus use ``definable'' to mean ``definable in $\CM$'', and use
``$\CD$-definable'' to mean ``definable in $\CD$''. Similarly, we use $\acl$, $\dim$
and `generic' to denote the
 corresponding
 notions in $\CM$, and let $\acl_\CD$, $\mr$, `$\CD$-generic' and `$\CD$-canonical parameter' denote the corresponding
  notions in $\CD$.

 Since the underlying universe of the strongly minimal $\CD$
is the $2$-dimensional set $G$, it follows that for every $\CD$-definable set $X\sub
G^n$, we have
$$\dim X =2\mr(X).$$ Also, for $a\in G^n$ and $A\sub G$, we have $$\dim(a/A)\leq
2\mr(a/A),$$ and in particular, if $X\sub G^n$ is definable in $\CD$ and $a\in X$ is
generic in $X$ over $A$, then it is also $\CD$-generic in $X$ over $A$. The converse
fails: indeed, let $\CM$ be the real field and $\CD$ the complex field,
interpretable in the real field $\CM$.  The element $\pi\in \mathbb C$ is
$\CD$-generic in
 $\mathbb C$ over $\emptyset$ but it is not generic in $\mathbb C$
 over $\emptyset$ because it is contained in the definable, $1$-dimensional
 set $\mathbb R$.

\subsection{The field configuration}
Recall the following definition.

\begin{definition}\label{D:gpconf}     Let $\CN$ be a strongly minimal structure.
A set $\{a,b,c,x,y,z\}$ of tuples is called a \emph{field
        configuration in $\CN$} if
    \[
    \groupconf{a}{b}{c}{x}{y}{z}
    \]
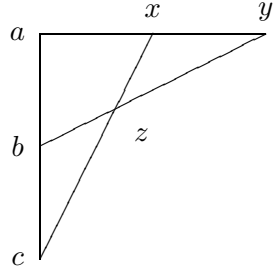
    \begin{enumerate}
        \item all elements of the diagram are pairwise independent and
        $\mr(a,b,c,x,y,z)=5$;
        \item $\mr(a) = \mr(b)= \mr(c)=2$, $\mr(x) = \mr(y) = \mr(z) = 1$;
        \item all triples of tuples lying on the same line are
        dependent,  and moreover, $\mr(a,b,c)=4$,
        $\mr(a,x,y)=\mr(b,z,y)=\mr(c,x,z)=3$;
        \item $\mr(\mathrm{Cb}(x,y)/a)=\mr(a)$, $\mr(\mathrm{Cb}(y,z)/c)=\mr(b)$ and $\mr(\mathrm{Cb}(x,z)/c)=\mr(c)$.
    \end{enumerate}
(For the notion  $\mathrm{Cb}$ of a canonical base, see \cite[page 19]{Pi}.)
\end{definition}

\begin{remark}\label{R:minimality} Consider the following minimality condition on a set $\{a,b,c, x, y,z\}$  of tuples in $\CN$:
\begin{enumerate}
  \item[(4)$'$] there are no $a'\in \acl(a)$, $b'\in \acl(b)$ and $c'\in \acl(c)$ with $\mr(a')=\mr(b')=\mr(c')=1$ such that the above (1) - (3) hold with $a',b',c'$ replacing $a,b,c$.
\end{enumerate}
    Standard Morley rank calculations show that the above conditions (1) - (4)  are equivalent to (1) - (3) and (4)$'$.
\end{remark}

For a proof of the  following theorem, see \cite[Main Theorem, Proposition
2]{BousGpCf}
and the discussion following Proposition 2 there.

\begin{fact}\label{Hrushovski}(Hrushovski)
    If a strongly minimal structure $\CN$ admits a field configuration,
    then $\CN$
    interprets an algebraically closed field.
    %(A field configuration for $k>3$ does not exist).
\end{fact}

Let $\mathbb G_m$ and $\mathbb G_a$ denote the multiplicative and additive groups of
an algebraically closed field $K$. The action of $\mathbb G_m \ltimes \mathbb G_a$
on $ \mathbb G_a$ (defined  by $(a,c)\cdot b=ab+c$) gives rise, naturally, to a field
configuration on the structure $(K,+,\cdot)$ as follows: take $g,h\in \mathbb G_m
\ltimes \mathbb G_a$ independent generics (in $K$), and $b\in \mathbb G_a$ generic
over $g,h$. Then
\[
\mathcal F:=\{h,g,gh,b,h\cdot b, gh \cdot b\}
\]

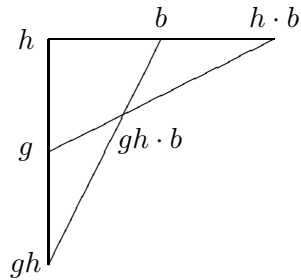
\[
\groupconf{h}{g}{gh}{b}{h\cdot b}{gh\cdot b}
\]
where $\cdot$ denotes the action of $\mathbb G_m \ltimes \mathbb G_a$ on $\mathbb
G_a$, is readily verified to be a field
configuration in the field $K$ (we will prove a slightly more general statement in Lemma \ref{lem-field} allowing us to construct field configurations form certain families of plane curves).
%Ultimately, the main result of the present paper will be obtained by constructing a field configuration in $\CD$.
%To achieve this we point out the following lemma.

When constructing a field configuration in Section \ref{sec-last}, we will need the lemma below. Given an algebraically closed field $K$,  denote by  $\mathrm{AGL}_1(K)$ the group of its affine transformations. Let $\CM$ be an o-minimal expansion of a real closed field, and $\CD$ a  $2$-dimensional definable strongly minimal structure. Here and below, we follow the conventions mentioned in Section \ref{sec-setting}. Namely, notions such as definability, genericity, $\dim$ and $\acl$ refer to $\CM$, unless indexed otherwise.
\begin{lemma}\label{L:Mcon2Dcon}
Let $K$ be a definable algebraically closed field and  $h,g\in \mathrm{AGL}_1(K)$ independent generics. Let $b\in K$ be a generic independent from $g,h$. Let  $\mathcal Y=\{h',g',k', b', c',d'\}\sub D^n$ be such that
    \begin{itemize}
        \item $h',g',b'$  are interalgebraic over $\0$ with  $h$, $g$, $b$ respectively, and
        \item $k',c',d'$ are interalgerbaic over $\0$ with $gh, h\cdot b, gh\cdot b$ respectively.
    \end{itemize}
  Then $\mathcal Y$ is a field configuration in $\CD$ if and only if it satisfies (3) of Definition \ref{D:gpconf}.
\end{lemma}
\begin{proof}
    Because $\CD$ is $2$-dimensional, if $S$ is a $\CD$-definable set, then by what we have already explained in  Section \ref{sec-setting}, $\dim(S)=2\mr(S)$. Since o-minimal dimension is preserved under interalgebraicity, it will suffice to show that (1), (2) and (4) of Definition \ref{D:gpconf} hold with $\mr$ replaced by $\frac{1}{2}\dim$.

    Because $\dim(K)=2$ we get that $\dim(\mathrm{AGL_1}(K))=4$. By exchange (in $\CM$), we get that (1) and (2) above hold. So it remains to verify (4). For that we note that, by genericity of $b$, for example, the point $(b,h\cdot b)$ is generic on the affine $K$-line $(x,h_1x+h_2)$ where $h=(h_1,h_2)$. Since any two distinct affine lines intersect in at most one point, any automorphism fixing the affine line $(x,h_1x+h_2)$ setwise must also fix $h$ (pointwise). So $h$ is a canonical base for $\mathrm{tp}_K(b,h\cdot b/h)$. Using the interalgebraicity it follows that $h'$ is
$\CD$-interalgebraic with $\mathrm{Cb}(b',c'/h')$. Similarly, the rest of clause (4) carries over from $K$ to $\CD$.
\end{proof}

\subsection{Notation}\label{sec-notation} If $S$ is a set in a topological space, its closure, interior, boundary and frontier are denoted by $\cl(S)$, $\intr(S)$, $\bdry(S):=\cl(S)\sm \intr(S)$ and
$\fr(S):=\cl(S)\sm S$, respectively.
Given a group $\la G, +\ra$ and sets $A, B\sub G$, we denote by $A-B$ the Minkowski
difference of the two sets, $A-B=\{x-y : x\in A,\, y\in B\}$. Given a set $X$ and
$S\sub X^2$, we denote $S^{\textup{op}}=\{(y,x)\in X^2 : (x,y)\in S\}$. The graph of a
function $f$ is denoted by $\Gamma_f$. If $\gamma: (a,b)\to R^n$ is a definable curve we will let $\gamma$ denote the image of $\gamma$ in $R^n$. Thus, if for some definable function $f$ we have $\gamma(t)\in \dom(f)$ for all $t$, we may write $f(\gamma)$ instead of $f(\mathrm{Im}(\gamma))$. For $\CM=\la R;+,\cdot,<,\ldots \ra$ as above and  $x=(x_1, \dots, x_n)\in R^n$, we write  $|x|=\sqrt{x_1^2+\dots+x_n^2}$.

%When we say that a property \emph{almost} holds, we mean that it holds except possibly for finitely many points. For example, a set $X$ almost contains (or coincides with) another set $Y$, if $Y\sm X$  (or $Y\triangle X$) is finite. Exception to this rule is the use of term ``almost faithful" (Definition \ref{def-almfaith}) and ``almost complex/$K$-structure" (Sections \ref{almost-complex}, \ref{sec-concl}).

%avoid \ni throughout
%ball, box
%definably connected

\section{Plane curves}\label{sec-prelim}
\label{sec-curves} In this section, we work in a strongly minimal structure $\CD$
and prove some lemmas about the central objects of our study, plane curves. When
$\CD$ expands a group $G$ and is non-locally modular, we construct in Sections
\ref{sec-dividing} and \ref{sec-ample} two special definable families of plane
curves which will be used in the subsequent sections.

%  we recall some basic facts about local modularity

\subsection{Some basic definitions and notations}\label{sec-basic curves}
Let $\CD$ be a strongly minimal structure.

%\noindent\textbf{Even though definability in this section is always taken with
%respect to $\CD$, we keep the index `$\CD$-' in notions such as `$\CD$-definable',
%`$\CD$-generic', etc, in order to be in accordance with the notation of the
%subsequent sections.}

\begin{defn} A \emph{$\CD$-plane curve} (or just \emph{plane curve}) is a $\CD$-definable subset
of $G^2$ of Morley rank $1$.
\end{defn}

%A definable family of plane curves is given by a definable set $C\sub D^2\times T$, with $T$ definable, such that for every $t\in T$, the set $C_t=\{(a,b)\in D^2:(a,b,t)\in C\}$ is
%a $\CD$-plane curve.\end{defn} We shall be mostly interested in $\CD$-definable families.

\begin{definition}\label{def-almfaith}
For two plane curves $C_1, C_2$, we write $C_1\sim C_2$ if $|C_1\triangle
C_2|<\infty$. Note that this gives a $\CD$-definable equivalence relation on any
$\CD$-definable family of plane curves.
 A $\CD$-definable family of plane curves $\CF=\{C_t : t\in T\}$ is
   \emph{faithful} if for $t_1\neq t_2$ in $T$, $C_{t_1} \triangle C_{t_2}$ is finite (i.e., $C_1\not\sim C_2$). It
   is \emph{almost faithful} if all $\sim$-equivalence classes are finite.
\end{definition}

Note that if $\CF=\{C_t : t\in T\}$ is a faithful family of plane curves, then $t$
is a canonical parameter for $C_t$. If $\CF$ is almost faithful, then $t$ is
inter-algerbaic with a canonical parameter of $C_t$.

Given a $\CD$-definable family of plane curves $\CF$, there exists
a $\CD$-definable almost faithful family of plane curves $\CF'=\{C'_t:t\in T'\}$ (possibly over additional parameters), such that every curve in $\CF$ has an equivalent curve in $ \CF'$ and vice versa (see for example, \cite{Hr2}, p.137)\footnote{Allowing imaginary elements, we can always
obtain faithful families of plane curves. The point here is to work in the real sort
only.}. It is not hard to see that $\mr(T')$ is independent of the choice of $\CF'$.
Thus, we can make the following definition.

\begin{defn}  A $\CD$-definable family of plane curves $\CF$ as above is said to be
{\em $n$-dimensional}, written also as $\mr(\CF)=n$, if in the corresponding almost
faithful family $\CF'$ as above, we have $\mr(T')=n$. We call $\CF$
\emph{stationary} if $\md (T')=1$.

\medskip We call $\CD$ a {\em non-locally modular structure} if there exists a $\CD$-definable family
 of plane curves $\CF$ with $\mr(\CF)\ge 2$.
\end{defn}

In fact, \cite[Proposition 5.3.2]{Pi}, if $\CD$ is non-locally modular, then for
 every $n$ there exists an $n$-dimensional $\CD$-definable family of plane curves.
  We will sketch a proof of a slightly stronger result in Proposition \ref{familyL} below.
  \smallskip

The following terminology is inspired by \cite{HrZil}.

%\textbf{Do we want to define Morley degree?}

%\textbf{There is an awkward switch of notation from $\CF$ to $\CL$}

%\textbf{In the last definition, why $\CD$-generic over $\es$ and not over
% the parameters used for $\CL$?}

\begin{defn}\label{def-ample}
Let $\mathcal F=\{C_t : t\in T\}$ be a $\CD$-definable family of plane curves. For every
$p\in G^2$, denote
$$T(p)= \{t\in T: p\in C_t\}\,\,;\,\, \CF(p)=\{C_t:p\in C_t\}.$$

We say that $\mathcal F$ is {\em (generically) very ample} if for every $p\neq q \in
G^2$ (each $\CD$-generic over the parameters defining $\CF$),
\[\mr(T(p)\cap T(q))<\mr(T(p)).\]
\end{defn}

$ $\\
\noindent\textbf{In the rest of this section,  $\CD=\la G; +, \dots\ra$ denotes a
strongly minimal expansion of a group $G$.}

\subsection{Local modularity} Here we recall some basic facts about local modularity.

\begin{definition}
    A $\CD$-definable set is {\em $G$-affine} if it
    is a finite boolean combination of cosets of $\CD$-definable subgroups of $G$.
\end{definition}

We use the following simple observation without further reference (see, for example,
\cite[Lemma 7.2.5, Corollary 7.1.6 and Corollary 7.2.4]{MaBook} for details).
\begin{remark}
    If $S\subseteq G^2$ is $\CD$-definable and strongly minimal, then $S$ is $G$-affine
     if and only if $S\sim H+ a$ for some $\CD$-definable strongly minimal subgroup $H\sub G^2$.
\end{remark}

\begin{defn} Given a strongly minimal plane curve $C$,
    the \emph{stabilizer of $C$} % (more precisely, the stabilizer of the generic type of $C$),
    is the set $$\stab^*(C)=\{g\in G^2: C\sim C+g\}.$$
\end{defn}
The stabilizer of $C$ is easily seen to be a $\CD$-definable subgroup of $G^2$. The
next properties are easy to verify.

\begin{lemma}\label{stabproperties} For $C$ a strongly minimal $\CD$-plane curve, and $p,q \in G^2$,
\begin{enumerate}
\item $C+p\sim C+q$ if and only if $p-q\in \Stab^*(C)$.

\item $\Stab^*(C)$ is trivial if and only if $\{C+p:p\in G^2\}$ is a faithful family.

\item $\Stab^*(C)$ is finite if and only if $\{C+p:p\in G^2\}$ is almost faithful.

\item $\Stab^*(C)$ is infinite if and only if $C$ is $G$-affine.
\end{enumerate}
\end{lemma}

 We only use here the following characterization of non-local
modularity in expansions of groups, which  follows from \cite{HrPi87}.

\begin{fact} \label{Hr-Pi} If $\CD$ is a strongly minimal expansion of a group $G$, then $\CD$ is
non-locally modular if and only if there exists a $\CD$-plane curve which is not
$G$-affine.
\end{fact}

%WE NEVER FIX SUCH CURVE $C$. Our plan is to fix one such curve $C$ which is not $G$-affine and then to consider the almost faithful family of curves $\{C+p:p\in G^2\}$. This family of curves -- subject to various manipulations -- will be our main object of investigation.

%\end{itemize}

Note that if $\CF=\{C+p:p\in T\}$ is a $\CD$-definable family of plane curves, with
$C$ strongly minimal and $T=G^2$, then for every $p\in G^2$,
$$T(p)=\{q\in T:p\in C+q\}=p-C.$$ In particular $T(p)$ is strongly minimal so that, in fact, if $\mr(T(p)\cap T(q))=\mr(T(p))$, then $T(p)\sim T(q)$. We thus have the following lemma.

\begin{lemma}\label{L:ample} If  $C$ is a strongly minimal plane curve and $\CF=\{C+p:p\in G^2\}$, then the following are
equivalent: \begin{enumerate} \item $\CF$ is very ample \item $\CF$ is faithful
\item $\Stab^*(C)$ is trivial.
\end{enumerate}
\end{lemma}

Finally, we will need the following definition.

\begin{defn}
   Let $\CF=\{C_t:t\in T\}$ be a $\CD$-definable family of plane curves. We call $C_t$ \emph{a $\CD$-generic curve in $\CF$} over $A$,  if  $t$ is $\CD$-generic in $T$ over $A$. We say that $\CF$ is {\em generically strongly
  minimal} if every $\CD$-generic curve in $\CF$ is strongly minimal.
\end{defn}

%\begin{lemma} \label{lemma1} Let $C\sub G^2$ be a $\CD$-plane curve.
%Then the following are equivalent: \begin{enumerate} \item there exist
%$(a_1,b_1)\neq (a_2,b_2) \in G^2$ such that
%$$\{(c,d)\in G^2 :(a_1,b_1)\in C+(c,d)\}\sim \{(c,d)\in G^2: (a_2,b_2)\in
%C+(c,d)\}.$$ \item   $\Stab^*(C)$ is not trivial. \end{enumerate}
%\end{lemma}
%\begin{proof}
%(1) is equivalent to :
%$$\{(a_1-c,b_1-d)\in C:(c,d)\in G^2\}\sim \{(a_2-c,b_2-d)\in
%C:(c,d)\in G^2\}, $$ and this is equivalent to  $(a_1-a_2,b_1-b_2)$ is in
%$\Stab^*(C)$.
%\end{proof}

%\noindent\textbf{Note:} The converse of the above statement is also easily seen to
%be true, but we are not going to make use of it here.

%\begin{lemma}  Let $C\sub G^2$ be a definable set of dimension $2$, such that $\Stab^*(C)$ is trivial. For every
% $t\in G^2$, let $C_t= C+t$.
%Then the family $\mathcal L=\{C_t : t\in G^2\}$  is faithful and very ample.
%\end{lemma}
%\begin{proof} Let $T=G^2$, and for every $(a, b)\in G$, $T(a,b)$ as in Definition \ref{def-ample}. By Lemma \ref{lemma1},
%\begin{equation}\label{eq1}
%\mbox{ for any }(a,b)\neq (a',b')\in G^2, \,\, \dim(T(a,b)\cap \mathcal T(a',b'))=0,
%\text{ and hence $T(a,b)\not\sim T(a',b')$}.\notag\end{equation} Hence $\mathcal L$ is
%very ample. It is easy to see that it is faithful.
%\end{proof}

\subsection{Dividing by a finite subgroup of $G$}\label{sec-dividing}

 The main goal of this subsection is to prove Lemma \ref{G and G/F} below, which will be
  used in the proof of Theorem \ref{main} in Section \ref{sec-last}. It will also allow us to assume, without
  loss of generality,
the existence of a $\CD$-definable faithful, very ample family of strongly minimal
plane curves of Morley rank
  two (Proposition \ref{familyL0} below).

Given a strongly minimal plane curve $C$ which is not $G$-affine, we plan to work
with the family $\CF=\{C+p:p\in G^2\}$. We know that $\Stab^*(C)$ cannot be infinite
but it can be a finite, non-trivial, group in which case $\CF$ is neither  faithful
nor very ample. We prove below that dividing the structure $\CD$ by a finite
group is harmless.

Given a finite subgroup $F\sub G$, $\CD$-definable over $\es$, we consider the map
$\pi_F:G\to G/F$, and still use $\pi_F:G^n\to (G/F)^n$ to denote the map
$\pi_F(g_1,\ldots, g_n)=(\pi_F(g_1),\ldots, \pi_F(g_n))$.

 We let $\CD_F$ be the structure whose universe is $G/F$ and whose atomic relations
 are all sets of the form $\pi_F(S)$ for $S\sub G^n$ a $\0$-definable set in $\CD$. The structure
$\CD_F$ is again an expansion of a group.

The following result implies that for the purpose of our main theorem we may work
with $\CD_F$ instead of $\CD$.

\begin{lemma}\label{G and G/F} Assume that the group $G$ has unbounded exponent.
Then the structures $\CD$ and $\CD_F$ are bi-interpretable, without parameters. In
particular, $\CD$ is bi-interpretable with an algebraically closed field if and only
if $\CD_F$ is.
\end{lemma}
\begin{proof} Because $F$ and $\pi_F$ are $\0$-definable in $\CD$,
the structure $\CD_F$ is interpretable, with no additional parameters, in $\CD$, via
the identity interpretation $\alpha(g+F)=g+F$.

%\begin{claim} If $S\ub G^n$ is definable in $\CD$ over $a\in G^k$ then $\pi_F(S)$ is
%definable in $\CD_F$ over $(aF)^k$. \proof We assume that $S=X_a$, where $X\sub
%G^n\times G^k$ is $\0$-definable in $\CD$.

% Clearly, every $\CD_F$-definable set is definable in $\CD$, so $\CD_F$ is
%interpretable in $\CD$. Also,  note that every $\CD$-definable subset of $(G/F)^n$
%is definable in $\CD_F$ (because its pre-image in $G^n$ is definable).

Next, let us see how we interpret $\CD$ in $\CD_F$.  Let $n=|F|$, and let
$\pi_F^*:G/F\to G$ be the map defined as follows: given $y\in G/F$, and $x\in G$ for
which $\pi_F(x)=y$, let
$$\pi_F^*(y)=nx.$$
Because $G$ is commutative, if $\pi(x)=\pi(x')=y$, then $nx=nx'+ng$ for some $g\in
F$. Since $ng=0$, this
 proves that $\pi_F^*$ is a well-defined group homomorphism with kernel $\pi_F(G[n])$,
where $G[n]=\{x\in G:nx=0\}$.

Since $G$ is strongly minimal and has unbounded exponent, the group $G[n]$ is finite
 and hence $\ker (\pi_F^*)$ is finite, so $\dim \mathrm{Im}(\pi_F^*)=\dim G/F=\dim
G$. Because $G$ is definably connected, $\pi_F^*$ is surjective. Thus the
homomorphism $\pi_F^*$ induces an isomorphism of $(G/F)/\pi_F(G[n])$ with $G$. Its
inverse $\beta:G\to (G/F)/\pi_F(G[n])$ is given by
$$\beta(g)=(\frac{g}{n}+F) +\pi_F(G[n]),$$ where $g/n$ is any element $h\in G$ such that
$nh =g$ (note that a strongly minimal group of unbounded exponent is divisible, \cite[\S3.3]{PoiGroups}).

By our assumptions,  $\pi_F(G[n])$ is $\0$-definable in $\CD_F$, and therefore the
quotient $(G/F)/\pi_F(G[n])$ is $\0$-definable in $\CD_F$. Now, given any
$\0$-definable $X\sub G^k$ in $\CD$, the set $\{(g_i/n)_{i=1}^k\in G^k:g\in X\}$ is
also $\0$-definable in $\CD$, and hence its image in $(G/F)^k/\pi_F(G[n])^k$ is
$\0$-definable in $\CD_F$. We therefore showed that $\CD$ is interpretable, without
parameters, in $\CD_F$ via $\beta$.

% It follows that the set $G$ is internal to the set $G/F$ in
%the structure $\CD$. Because every subset of $(G/F)^n$ that is definable in $\CD$
%over $\0$ is also definable in $\CD_F$ over $\0$ it follows that $\CD$ is
%interpretable in $\CD_F$, without parameters.

To see that this is indeed bi-interpretation,  we first note that the isomorphism
between $\CD$ and its interpretation in $\CD_F$ is $\alpha\circ \beta$, which equals
$\beta$. It is  clearly definable in $\CD$.

Let us examine the map  induced on  $G/F$ by $\beta\circ \alpha$ and prove that it
is definable in $\CD_F$. We denote by $F/n$ the preimage of $F$ in $G$ under the map
$g\mapsto ng$. It is not hard to see that the image of $F$ inside $\beta(G)$ is the
group $\pi_F(F/n)+\pi_F(G[n])=\pi_F(F/n)$, and hence the isomorphism which
$\beta\circ \alpha$ induces on $G/F$ is
$$g+F\mapsto g/n+\pi_F(F/n).$$ This map is definable in the group $G/F$ by sending
$g+F$ to the unique  coset $h+F/n$ such that $nh+F=g+F$.

 This completes the proof that $\CD$ and $\CD_F$ are bi-interpretable over $\es$.\end{proof}

 Note that in our case, when the group $G$ is abelian and definable in an o-minimal
 structure, then by
\cite{Strz}, the group $G$ has unbounded exponent, so the above result holds.

For the rest of this subsection, assume that $\CD$ is non-locally modular, and fix
(after possibly absorbing into the language a finite set of parameters) a strongly
minimal plane curve $C\sub G^2$ which is $\CD$-definable over $\es$ and not
$G$-affine. By Lemma \ref{stabproperties}(4), $F'=\Stab^*(C)\subseteq G^2$ is a finite subgroup and let
 $F\sub G$ be a $\CD$-$\0$-definable  subgroup such that $F'\sub F\times F$. Consider the
  structure $\CD_F$ expanding $\la G/F,+\ra $
 as above.

 \begin{claim}\label{Stabtrivial} $\pi_F(C)$ is strongly minimal in $\CD_F$ and $\Stab^*( \pi_F(C))$
  in $(G/F)^2$ is trivial.\end{claim}
 \begin{proof} The strong minimality of $\pi_F(C)$ is immediate from the strong minimality
 of $C$ in $\CD$.

 Assume that $q\in \Stab^*(\pi_F(C))\sub (G/F)^2$, namely
 $q+\pi_F(C)\sim \pi_F(C)$. Let $\tilde F=F\times F\sub G^2$ and fix $p\in G^2$ such that $\pi_F(p)=q$.
  Then $p+C+\tilde F\cap C+\tilde F$ is infinite and since $\tilde F$ is finite, there
  exist $g,h\in \tilde F$ such that $C+p+g\cap C+h$ is infinite. But then $p+g-h\in
  \Stab^*(C)\sub \tilde F$, implying that $p\in \tilde F$, and hence $0=\pi_F(p)=q$.

  We thus showed that $\Stab^*(\pi_F(C))$ is trivial.\end{proof}

Combining  Lemmas \ref{L:ample}, \ref{G and G/F} and Claim \ref{Stabtrivial}, we can
conclude the following statement.

%there is a strongly minimal plane curve $C$ with $\Stab^*(C)=\{0\}$. By Lemma
%\ref{L:ample}, the family  $\{C+p:p\in G^2\}$ is  faithful and very ample.
% We can thus conclude the following proposition.

\begin{proposition}\label{familyL0}
Assume $\CD$ is non-locally modular, expanding a group $G$ of unbounded exponent.
Then there exists a finite group $F\sub G$, possibly trivial, and in  the structure
$\CD_F$ defined above there exists a definable family  $\CL=\{l_t:t\in Q\}$,  of
strongly minimal plane curves, which is faithful, very ample, and
$\mr(Q)=2$.

The structures $\CD$ and $\CD_F$ are bi-interpretable, over the parameters defining
$F$.
\end{proposition}

\noindent{\bf Assumption: for the rest of the article, we replace the structure
$\CD$ with the structure $\CD_F$, and thus assume that a family $\CL$ as above is
definable  in $\CD$.}

%\noindent We fix $\CL$ as above.

% This family will play, in several applications a role similar to that of the family of affine
%  lines in $\mathbb C^2$.
%  (\textbf{not really several, only in the Frontier section. My suggestion is not to fix it.}\\

%For $p\in G^2$ we let $\CL(p)$ denote the subfamily of curves passing through $p$. Namely, $\{q\in G^2: p\in C+q\}$, equivalently, $\{q\in G^2: -q\in C-p\}$. So $\CL(p)=p-C$, and in particular, $\CL(p)$ is strongly minimal. In particular
%Very ampleness of $\CL$ is implies the stronger condition that $\CL(p)\cap \CL(q)$ is non-generic in neither $\CL(p)$ or $\CL(q)$, provided $p\neq q$.

\subsection{Very ample families of high dimension}\label{sec-ample}

  The goal of this subsection is to construct a larger family $\CL'$ of plane curves which still has the
geometric properties of the family $\CL$ from Proposition \ref{familyL0}. The main
method is to use composition of binary relations and families of plane curves.
Recall the notion of a composition of binary relations, extending composition of
functions: given $S_1,S_2\sub G^2$, we let
$$S_1\circ S_2=\{(x,z)\in G^2:\exists y (x,y)\in S_2\,\text{ and }\, (y,z)\in S_1\}.$$
Clearly, if $S_1,S_2$ are $\CD$-definable, then so is $S_1\circ S_2$. We will be
mostly interested in the composition
 of plane curves, and even more so, in the composition of families of plane curves: if $\CL_1,\CL_2$ are $\CD$-definable families
of plane curves, we let $\CL_1\circ \CL_2:=\{C_1\circ C_2: C_1\in \CL_1, C_2 \in
\CL_2\}$.

\begin{definition}
%\stodo{12.8-changed slightly def}
A plane curve $S\sub G^2$ is \emph{a straight line} if there exists $a\in G^2$, such that either $S \sim \{a\}\times G$ or $S\sim G\times \{a\}$.
\end{definition}

As a rule, geometric properties are not preserved under compositions of (families
of) curves.
 The composition of two strongly minimal curves, which are not both straight lines, has, indeed, Morley rank 1, but it need not be strongly minimal. More
  generally, a $\CD$-generic curve of $\CL_1\circ \CL_2$ need not be strongly minimal, and even if it were,
   $\CL_1\circ \CL_2$ need not be faithful. In fact, although the dimension of $\CL_1\circ \CL_2$
cannot decrease, it need not be greater than that of $\CL_1$ or $\CL_2$. For
example, if both families are the family of affine lines in $\mathbb A^2$, then
$\CL_1\circ \CL_2=\CL_1$.\\

We will need a series of lemmas to address these issues. We start with the following
easy observation.

\begin{lemma}\label{lem-comp} Assume that $\CL_1=\{C_t:t\in T\}$ and
    $\CL_2=\{D_r:r\in R\}$ are two $\CD$-definable almost faithful families of plane curves, none of which is a straight line, and let $\CL=\CL_1\circ\CL_2$.
    \begin{enumerate} \item For every $\CD$-generic $p$ in $G^2$, we have $\mr(\CL(p))=\mr(R)+\mr(T)-1$.
    \item If $\CL_1$ and $\CL_2$ are generically very ample, then so is $\CL$.
    \end{enumerate}
\end{lemma}
%\stodo{12.8-changed below. Also, in the statement of 3.17, don't we want to add that the family $E\circ L$ is %generaically strongly minimal?}

\begin{proof}
(1) Let $\CL:=\CL_1\circ \CL_2$, $C\in \CL$ a $\CD$-generic curve and $(a,b)\in C$ a $\CD$-generic point. So $[C]$ forks over $(a,b)$ and therefore $\mr([C]/(a,b))\le \mr(R)+\mr(T)-1$. Let us see that equality holds, or equivalently $\mr(\mathcal L(a,b))=\mr(R)+\mr(T)-1$. Fix some $\CD$-generic $e\in G$. Then $(e,b)$ is $\CD$-generic in $G^2$. So $\CL_1(e,b)$ has Morley rank $\mr(T)-1$. Similarly $\CL_2(a,e)$ has Morley rank $\mr(R)-1$. So the set
$$\CL(a,b)_e:=\{(t,r)\in T\times R: (a,e)\in D_r\land  (e,b) \in C_t\}$$ has rank $\mr(R)+\mr(T)-2$. But because for $\CD$-independent generics $e,e'$, the sets $\CL(a,b)_e$ and $\CL(a,b)_{e'}$ are disjoint up to a set of lower rank, $\mathcal L(a,b)$ has rank $\mr(R)+\mr(T)-1$.

(2) In order to show that $\CL$ is generically very ample it will suffice to show that $\CL(a,b)\cap \CL(c,d)$ has rank at most $\mr(R)+\mr(T)-2$ for $(a,b),(c,d)$ distinct $\CD$-generics. Fix some $r\in R$.  Then for $t\in T$ we have that $(t,r)\in \CL(a,b)$ only if for some $e$ such that $(a,e)\in D_r$ we also have $(e,b)\in C_t$. If, in addition $(t,r)\in \CL(c,d)$, then there exists $e'$ such that $(a,e')\in D_r$ and $(e',b)\in C_t$. But as there are only finitely many $e$ such that $(a,e)\in D_r$ and only finitely many $e'$ such that $(e',b)\in C_t$ it follows that there is a $\CD$-generic (over $a,b,c,d$) element of $\CL_1(e,b)$ that is also an element of $\CL_1(e',d)$. Unless $b=d$, this contradicts generic very ampleness of $\CL_1$. So we are reduced to the case where $b=d$, in which case a symmetric argument will show that unless also $a=c$ we get a similar contradiction. But since $(a,b)\neq (c,d)$, we are done.
\end{proof}

% I HAVE NO IDEA WHY WE ARE FIXING SUCH $C$ AND WHERE THAT $C$ IS EVEN USED. Next, we fix As we have already mentioned $C$ need not be strongly minimal. In order to sort out this situation we need a few more technical observations. First, we need:

\begin{definition}
    Given two $\CD$-definable families of plane curves, $\CL$ and $\CL'$, we say that {\em $\CL$ extends $\CL'$}  if for every
     $C'\in \CL'$ there exists
     $C\in \CL$ such that $C'\subseteq C$.
\end{definition}

In the next couple of lemmas we show that although the composition of two families
of curves need not preserve the properties of the original families (as already
discussed), it extends a family of curves that does.

\begin{lemma}\label{L:almost-faithful-comp}
    Let $\CL$ be a $k$-dimensional almost faithful $\CD$-definable family of plane curves.
     Let $E$ be a plane curve. Assume neither $E$ nor any $\CD$-generic plane curve is a straight line.
    Then $E\circ \CL$ extends a $k$-dimensional almost faithful $\CD$-definable family of plane curves whose $\CD$-generic members are strongly minimal. In fact, if $C\in \CL$ is
     $\CD$-generic over $[E]$, then for any strongly minimal $C_E\subseteq E\circ C$ we have $\mr([C_E]/[E])=k$.
\end{lemma}
\begin{proof}
    Fix some $C\in \CL$ which is $\CD$-generic over $[E]$ and $C_E\subseteq E\circ C$ strongly minimal. Note that $(E^{-1}\circ C_E)\cap C$
    is infinite, and since $C$ is strongly minimal, $E^{-1}\circ C_E$ is a set of Morley rank 1, containing the set
     $C$,
     up to a
    finite set. It follows that $[C]\in \acl_\CD([C_E][E])$. Since $\mr([C]/[E])=\mr([C])$, we get, by exchange,
     that $\mr([C])=\mr([C]/[E])=\mr([C_E]/[E])$.

    Absorbing $[E]$ into the language,  we can find $\bar c\in \acl_\CD([C])$ and a formula $\phi(x,\bar c)$ defining $C_E$.
     By compactness, there is a formula $\theta\in \tp(\bar c)$ such that whenever $\bar c'\models \theta$ there is some
      $C'\in \CL$ such that $\phi(x,\bar c')\subseteq E\circ C'$, and for all $\CD$-generic $\bar c'\models \theta$ the formula
       $\phi(x,\bar c')$ is strongly minimal. We may further require -- by compactness, again -- that if
        $\phi(x,\bar c')\land \phi(x,\bar c'')$ is infinite, then the symmetric difference
         $\phi(x,\bar c')\triangle \phi(x,\bar c'')$ is finite for all $\bar c',\bar c''\models \theta$. %\stodo{12.8-changed the last line of proof}
          By rank considerations, the family $\{\phi(G^2,\bar c'):\theta(\bar c')\}$ is almost faithful of rank $k$.
\end{proof}

As an immediate application (since the only families of straight lines are $1$-dimensional), we get the following statement.
\begin{corollary}\label{cor-compfam}
    Let $\CL_1, \CL_2$ be almost faithful $k$-dimensional $\CD$-definable families of plane curves, $k>1$. Then $\CL_1\circ \CL_2$ extends an almost faithful, stationary, generically strongly
    minimal family of plane curves of dimension at least $k$.
\end{corollary}

We can now show that a $2$-dimensional family of plane curves closed under composition (such as the family of affine lines in a field) gives rise to a field configuration.

\begin{lemma}\label{lem-field}
    Let $\CL_1, \CL_2$ be almost faithful $2$-dimensional families of plane curves. Assume that $\CL_1, \CL_2$ are $\CD$-definable over $\es$. Let $X\in \CL_1$ and $Y\in \CL_2$ be $\CD$-independent generic curves, and $E\sub X\circ Y$ strongly minimal.

    \begin{enumerate} \item If $\mr([E]/\emptyset)=2$,
  then  $\CD$ interprets an infinite field.
  \item If $\mr([E]/\emptyset)=k>2$ and $\CL_1, \CL_2$ are generically very ample, then $\CL_1\circ\CL_2$ extends a  $k$-dimensional almost faithful,   generically strongly minimal, stationary  and generically very ample family of curves.
 \end{enumerate}
\end{lemma}
\begin{proof}
 (1) As we note at the beginning of the proof of Lemma \ref{L:almost-faithful-comp}, each one of $[X],[Y],[E]$ is in the algebraic closure of the other two.  Since $\CL_1$ and $\CL_2$ are almost faithful and $2$-dimensional, we have $\mr([X]/\0)=\mr([Y]/\0)=\mr([E]/\0)=2$, and the Morley rank of each two of $[X],[Y],[E]$ is $4$.

    Now choose a $\CD$-generic $(x,y)\in X$, and $z$ so that $(y,z)\in Y$, and hence  $(x,z)\in E$.
    We claim that $\{[X],[Y],[E],x,y,z\}$ is a field configuration as in Definition \ref{D:gpconf}.
    We have
    $$\mr([X],x,y/\0)=\mr([Y],y,z/\0)=\mr([E],x,z/\0)=3 $$ and  $$\mr(x/\0)=\mr(y/\0)=\mr(z/\0)=1.$$
    Also, the Morley rank of the whole configuration over $\0$ is $5$.  It is thus left to verify (4) of Definition \ref{D:gpconf}.

    Because $(x,y)\in X$ and $\mr([X]/\0)=2$, we have $\mr(\mathrm{Cb}(x,y/[X]))=2$. We similarly verify the other conditions and therefore  $\{[X],[Y],[E],x,y,z\}$ is indeed a field configuration. By Fact \ref{Hrushovski}, an infinite field is interpretable in $\CD$.

\smallskip (2) Let $\CL:=\CL_1\circ \CL_2$. Let $\CF$ be an $\acl_\CD(\0)$-definable, generically strongly minimal, almost faithful and stationary family of plane curves, so that $E$ is contained, up to finitely many points, in a $\CD$-generic member of $\CL$ (such a family always exists). The family $\CL$ extends $\CF$ and by Lemma \ref{lem-comp}(2) is generically very ample. We need to show that so is $\CF$.

    Since $\CF$ is $k$-dimensional and almost faithful, for every $\CD$-generic $p\in G^2$, we have $\mr(\CF(p))=k-1$.
    It is thus sufficient to prove that for $p,q$, each $\CD$-generic in $G^2$, we have $\mr(\CF(p)\cap \CF(q))<k-1$.

    We write $\CL_1=\{C_t:t\in T\}$ and $\CL_2=\{D_r:r\in R\}$. By our assumption on $E$, it is a strongly minimal
    subset of $C_t\circ D_r$, for $(t,r)$ $\CD$-generic in $T\times R$. That is, $\mr(t,r/\0)=4$. It follows that $\mr(t,r/[E])=4-k$ and thus for every $\CD$-generic $(t',r')\in T\times R$ there exists a strongly minimal $E'\sub^* C_{t'}\circ D_{r'}$ in $\CF$ with
    $\mr(t',r'/[E'])=4-k$. Here $\sub^*$ means ``contained up to finitely many points''.

    Now let $p,q$ be $\CD$-generic in $G^2$ and  assume towards contradiction that $\mr(\CF(p)\cap \CF(q))=k-1$. Take $E'$ $\CD$-generic in $\CF(p)\cap\CF(q)$ over $p,q$.
    Consider the set $$P=\{(t_1,r_1)\in T\times R: E' \subseteq^* C_{t_1}\circ D_{r_1}\}.$$ Since $\mr(t',r'/[E'])=4-k$ and $(t',r')\in P$, we have $\mr(P)\geq 4-k$.  Fix $(t_0,r_0)$ $\CD$-generic in $P$ over $[E']$, $p$ and $q$.
    We have \[\mr(t_0, r_0,[E']/p,q) =\mr(t_0,r_0/[E'],p,q)+\mr([E']/p,q) \geq (4-k )+ (k-1)=3.\]

    Finally, since $[E']\in \acl_\CD (t_0,r_0)$, we have $\mr(t_0,r_0/p,q)\geq 3$ and in addition $(t_0,r_0)\in (\CL_1\circ L_2)(p)\cap (\CL_1\circ\CL_2)(q)$ (because $E'\sub C_{t_0}\circ D_{r_0}$).  However, by Lemma \ref{lem-comp}(1),(2) we have $\mr(\CL_1\circ \CL_2)(p)\cap (\CL_1\circ\CL_2)(q))<3$, contradiction. Thus $\CF$ is indeed generically very ample.
    \end{proof}

%\begin{lemma} \label{lem-comp} Assume that $\CL_1=\{C_t:t\in T\}$ and
%$\CL_2=\{D_r:r\in R\}$ are two very ample families of plane curves. Then the family
%$$\CL_1\circ \CL_2=\{C_t\circ D_r: t\in T, r\in R\}$$ is also very ample.
%\end{lemma}
%\proof Denote, for every $(a, b)\in G^4$,
%$$S(a,b)=\{(t, r)\in T\times R: (a, b)\in C_t\circ D_r\}.$$
%It suffices to find, for every $(a',b')\neq (a,b)\in G^2$ a generic $(t,r)\in
%S(a,b)$ over $(a,b)$ such that $(a', b')\notin C_t\circ D_r$. Without loss of
%generality $a\neq a'$. Choose $y_0\in G$ generic over $a,b,a',b'$ and let $r$ be
%generic in $R(a,y_0)=\{r\in R: (a,y_0)\in D_r\}$ over $r,y_0$. Let $t\in T$ be
%generic in $T(y_0,b)$ over $\{r,a,b,a',b',y_0\}$. Then $C_t\circ D_r$ is a generic
%curve in $S(a,b)$. Notice that if $(D_r)_a=\{y_1,\ldots, y_k\}$ then by  genericity
%of $t$ over $r,a',b'$, none of the $(y_i,b')$ is in $C_t$ and therefore
%$(a',b')\notin C_t\circ D_r$. Thus, $\CL_1\circ\CL_2$ is very ample. \qed
%\\
%
%{\bf The following needs to be rewritten:
%
%
%We would like to reach next:
%
%$\bullet$ The family $\CL\circ \CL$, or an extract of it,  is very ample, faithful,
%with $\mr(\CL)\geq 3$.
%
%$\bullet$ Every generic curve curve in $\CL\circ CL$ has Morely degree one.
%
%
%}

Under our standing assumptions at the end of Sections \ref{sec-basic curves} and \ref{sec-dividing}, we can finally conclude the last result of this section.

\begin{proposition}\label{familyL}
There exists a $\CD$-definable almost faithful, stationary family of generically strongly
minimal plane curves, $\mathcal F=\{C_t : t\in T\}$,  which is generically very ample,
and $\mr(T)\ge 3$.
\end{proposition}
\begin{proof}
Let $\CL$ be an almost faithful family of rank $2$ as in Proposition \ref{familyL0} and consider the family $\CL\circ \CL$. Let $C\in \CL\circ \CL$ be $\CD$-generic. By Lemma \ref{L:almost-faithful-comp}, there exists a strongly minimal $E\subseteq C$ with $\mr([E])\ge 2$. Either $\rm([E])=2$ and by Lemma \ref{lem-field}(1) there is an infinite field interpretable in $\CD$, in which case a family as required exists (take the family of graphs of polynomials of degree $d>1$ over $K$) or $\mr([E])>2$ in which clause (2) of Lemma \ref{lem-field} gives a $\CD$-definable family of curves as required.
\end{proof}

$ $\\
\noindent\textbf{From now on, until the end of the paper, we fix a sufficiently
saturated o-minimal
 expansion of a real closed field  $\CM=\la R;+,\cdot,<,\ldots \ra$, and a $2$-dimensional group
 $G=\la G; \oplus\ra$ definable in $\CM$. We also fix a strongly minimal  non-locally modular structure
 $\CD=\la G;  \oplus, \dots\ra$ definable in $\CM$. As discussed in Section \ref{sec-setting}, we include the index
 $\CD$ when referring to definability, genericity and such in the structure $\CD$, and omit
  the index when referring to $\CM$. We also assume the existence of a
  $\CD$-definable, very ample stationary family of plane curves $\CL$, as noted
   after Proposition \ref{familyL0}.
  In Sections \ref{sec-frontier} - \ref{sec-topcor}, we denote $\oplus$ by $+$, for simplicity.}
$ $\\

\section{Frontiers of plane curves}

\label{sec-frontier}

%{\bf Notation (to put somwhere?)} $\dim(t):=\dim (t/\es)$.

%\textbf{Kobi/Assaf: in this section, is genericity always in $\CM$? Or have we forgotten to say $\CD$-generic
% sometimes?}

\subsection{Strategy}

Our goal is to show (Theorem \ref{main-frontier}) that if $S\sub G^2$ is a plane curve, then its frontier $\fr(S)$ is
 finite and
in fact contained in $\acl_{\mathcal D}([S])$.  The geometric idea originates in
\cite{PeStExpansions} and it is implemented in Lemma \ref{FrgoodS} below, as
follows. We consider the  family $\CL$ from the assumption following Proposition
\ref{familyL0}. We also fix $b\in \fr(S)$ and consider a curve $l_q\in \CL$ going
through $b$ with $q$ generic over $[S]$. If $l_q$ meets $S$ transversely at every
point of intersection and $b$ is sufficiently generic in $G^2$, then by moving $l_q$
to an appropriate $l_{q'}$ close to $l_q$, the curve $l_{q'}$ will intersect $S$
near all points of $l_q\cap S$, and in addition at a new  point near $b$. Since $b$
itself was not in $S$, it follows that a generic $l_q$ through $b$ intersects $S$ at
fewer points than a generic curve in $\CL$. Thus $b$ is $\CD$-algebraic over $[S]$
and in particular $\fr(S)$ is finite.

While this strategy works well when the curves in $\CL$  are complex lines in
$\mathbb C^2$, the problem becomes more difficult when they are arbitrary plane
curves and $b$ is not necessarily generic in $G^2$. To get around this problem, the
idea in \cite{HaKo} was to replace $S$ by its image under composition with a
``generic enough'' curve from a new ``large'' family $\mathcal L'$ (Proposition
\ref{familyL}). We carry out this replacement in
 Lemma \ref{massageS} below. An additional complication of this strategy in the current setting is that instead
 of the functional language in \cite{HaKo} we need to work with arbitrary curves, and control  their
composition.

\subsection{Two technical lemmas about $2$-dimensional sets in $G^2$}

The following lemmas will be used in the sequel.

\begin{lemma}\label{intersectionfinite}
  Assume that $\{Y_e:e\in E\}$ is a definable family of $2$-dimensional subsets of $G^2$, with $\dim E=k\ge 2$.
    Assume that for all $e\in E$ there are at most finitely many $e'\in E$, such that
    $|Y_e\cap Y_{e'}|=\infty$.
 Then
$\dim(\bigcup_{t\in E} Y_t)= 4$.
\end{lemma}
\begin{proof}
The set $$\{(e,s):e\in E\, , \, s\in Y_e\}$$ has dimension $k+2$. Therefore, if the
union of the $Y_e$ had dimension  smaller than 4, then for a generic $s$ in this
union, the dimension of $E(s)=\{e\in E: s\in Y_e\}$ is at least $k-1\ge 1$, and in
particular, is infinite. Hence, there are $e_1,e_2\in E(s)$, independent and generic
over $s$. Therefore, $\dim(e_1,e_2/s)=2k-2$ and hence $\dim(e_1,e_2,s)=2k-2+3=2k+1$.
But this is impossible since $\dim(e_1,e_2/\es)\leq 2k$ and,  by our assumption on
the family, the set $Y_{e_1}\cap Y_{e_2}$ is finite, so $s\in \acl(e_1,e_2)$.
\end{proof}

%\textbf{To be defined:} submersion

\begin{defn} We say that two $2$-dimensional sets $C_1$ and $C_2$ {\em intersect
transversely at $p\in C_1\cap C_2$} if $C_1$ and $C_2$ are both smooth at $p$, and
their tangent spaces at $p$ generate the full tangent space of $G^2$ at $p$, namely
$T_p C_1+T_p C_2=T_p G^2$.
\end{defn}

\begin{lemma}\label{transversal} Let $\mathcal L=\{l_q : q\in Q\}$ be a $\0$-definable family of $2$-dimensional subsets
of $G^2$, and $S\sub G^2$ a $\0$-definable $2$-dimensional set. Let $q$ be generic in $Q$ over $\emptyset$ and assume that
$l_q$ and $S$ intersects transversely at $s$. Then for every  neighborhood $U\sub
G^2$ of $s$,
 there exists a neighborhood $V\sub Q$ of $q$, such
that for every $q'\in V$, we have $l_{q'}\cap S\cap U\neq \emptyset$.
\end{lemma}

\begin{proof} Without loss of generality, $U$ is definable over $\0$ and $l_q\cap U$ is smooth (we can shrink
it so that $q$ is generic in $Q$ over the parameters defining it). Reducing $U$
further, if needed, we may -- by cell decomposition, and the assumption that $l_q$
is smooth at $s$ -- write $l_q\cap U$ as the zero set of a definable $C^1$-map
$F_q:U\to R^2$, and similarly write $S$ as the zero set of a $C^1$-map $G:U\to R^2$.
The transversal intersection of $l_q$ and $S$ implies that the joint map
$(F_q,G):U\to R^4$ is a diffeomorphism at $s$, so in particular there is $U_0\sub U$
such that $(F_q,G)$ is a diffeomorphism on $U_0$ and $\bar 0\in R^4$ is in its open
image. We may choose $U_0$ so $q$ is still generic over the parameters defining
$U_0$. It follows that there is a neighborhood $V\sub Q$ of $q$ such that for every
$q'\in V$, $l_{q'}\cap U=F_{q'}^{-1}(0)$, for some definable $F_{q'}:U_0\to R^2$,
and the map $(F_{q'},G)$ is still a diffeomorphism on $U_0\sub U$, with $\bar 0$ in
its image. But now, if $(F_{q'},G)(s')=\bar 0$, then $s'\in U_0\cap l_{q'}\cap S$.
\end{proof}

% \textbf{Needed?:} Clearly, for every $(a,b)\in G^2$, $\dim(\mathcal T_{a,b})=2$ (\textbf{why?}). \textbf{Seems we need $\mathcal T_{c_2, b_2}=6$ below!}

\subsection{Bad points}

Recall that $\CL=\{l_q:q\in Q\}$ is a faithful and generically very ample $\CD$-definable
family of strongly minimal plane curves, with $\mr(Q)=2$. Notice that for $b\in G^2$
generic, the set $Q(b)=\{q\in Q:b\in l_q\}$ has Morley rank $1$. As in Section \ref{sec-prelim}, we let
$\CL(b)=\{l_q: q\in Q(b)\}$.

\begin{defn}\label{D:fibers} Let $U\sub G^2$ be an open set and $b\in G^2$. We say that \emph{$\mathcal L(b)$ fibers $U$} if for every
$s\in U$ there exists a unique $q\in Q(b)$ such that $s\in l_{q}$, the set $Q(b)$ is
smooth at $q$  and furthermore the function $s\mapsto q: U\to Q(b)$ is a submersion
at $s$ (that is, the differential map between the tangent spaces is surjective).
\end{defn}

\begin{defn} For $b\in G^2$, we say that  a point $s=(s_1,s_2)\in G^2$
 is  {\em $b$-good} if
\begin{enumerate}

\item There exists an  open neighborhood $U\sub G^2$ of $s$ such that the family
$\CL(b)$ fibers $U$.

\item For all $q\in Q(b)$ such that $s\in l_q$, the curve $l_q$ is smooth at $s$.

\end{enumerate}
Otherwise, we say that  $s$ is a \emph{$b$-bad}. We denote by $\Bad(b)$ the set of all
$b$-bad points.
\end{defn}
Clearly, the set $\Bad(b)$ is definable over $b$.

\begin{lemma}\label{b-good} For every $b\in G^2$, the set $\Bad(b)$ has
dimension at most $3$.
\end{lemma}
\begin{proof} Note that since  $\CL(b)$ is faithful, it follows that
$\mr(G^2\setminus \bigcup_{q\in Q(b)}l_q)\leq 1$.

 By cell decomposition, for a fixed generic $q\in Q$, the set
of points $s\in l_q$ failing (2) is at most $1$-dimensional. So the set of all
points $s$ failing (2)  is at most $3$-dimensional.

We now fix $s\in G^2$ generic over $b$, and show that it satisfies (1).
 %, and hence that the set of points failing (3) is at most $3$-dimensional.
The set of singular points $q$ on $Q(b)$ has dimension one, and for every such $q$,
$l_q$ has dimension $2$. Thus, the union of all such $l_q$ has dimension at most
$3$, and does not contain $s$. So if $s\in l_q$ for some $q\in Q(b)$, then $q$ is a
smooth point on $Q(b)$.

Since $s$ is generic in $G^2$, there are at most finitely many curves in $\CL(b)$
containing $s$. Hence, there is an open neighborhood $W\sub Q(b)$ such that $W\cap
Q(b)\cap Q(s)=\{q\}$.  We may choose $W$ to be definable over generic parameters.
Hence the first-order property over $b$: ``$\varphi(s'):=|W\cap Q(s')\cap Q(b)|=1$''
must hold for all  $s'$ in a neighborhood $U\sub G^2$ of $s$. Let $g:U\to Q(b)$ be
the map sending $s'$ to the unique $q'\in W\cap Q(b)$ with $s'\in l_{q'}$. Note that
for every $q'\in g(U)$, $g^{-1}(q')=l_{q'}\cap U$. Since the family $\CL(b)$ is
faithful, we have $\dim g(U)=2=\dim Q(b)$, and by the genericity of $s$ in $\dom
(g)$, the function $g$ is a submersion at $s$, thus $s$ is a $b$-good point.
\end{proof}

\subsection{Finiteness of the frontier}\label{sec-ff}

%The following easy observation simplifies the exposition.
%\begin{lemma}\label{L: S sm}
   % If $\fr(S)\subseteq \acl_\CD([S])$ for any $\CD$-definable strongly minimal $S\subseteq G^2$ which is not a straight line, then $\fr(S)\subseteq \acl_\CD([S])$ for any  $\CD$-definable set $S\subseteq G^2$ of Morley rank one.
%\end{lemma}
%\begin{proof}
   % Note that $\fr(\bigcup\limits_{i=1}^k S_i)\subseteq \bigcup\limits_{i=1}^k \fr(S_i)$. Since any set $S\subseteq G^2$ of Morley rank one can be written as $\bigcup\limits_{i=1}^k S_i$ for some strongly minimal sets  $\CD$-definable over $\acl_\CD([S])$ it will suffice to prove the lemma for $S$ strongly minimal.  If $S$  is a straight line, say $\{a\}\times G$ then $\fr(S):=(\{a\}\times G)\setminus S$, so it is $\CD$-definable and finite. Otherwise $\fr(S)\subseteq \acl_\CD([S])$ by assumption. So the result follows.
%\end{proof}

%\noindent\textbf{Until the end of Section \ref{sec-ff}, unless explicitly stated otherwise, $S$ will denote a $\CD$-definable strongly minimal subset of $G^2$ not  coinciding with a straight line.}\\

The heart of the geometric argument is contained in the following lemma showing
that in a generic enough setting the frontier of $S$ is indeed contained in
$\acl_{\CD}([S])$.

%\textbf{In the statement we said we should include $[S]$, but I am not sure what $[S]$ is. Was it $[\mathcal S]$, or $[S_{t_0}]$?}

\begin{lemma}\label{FrgoodS}
Let $\mathcal F=\{S_t: t\in T\}$ be a $\CD$-definable stationary almost faithful family of plane
curves with $\mr(T)\geq 3$. Assume that  $b\in G^2$ with $\dim(b/\emptyset)=4$, and
$t_0\in T$ generic over $\emptyset$.  If $b\in \fr(S_{t_0})$, then $b\in \acl_{\mathcal
D}(t_0)$.
\end{lemma}
\begin{proof}
We may assume first that $S_{t_0}$ is strongly minimal. Indeed, $S_{t_0}$ is
a finite union of strongly minimal sets, each definable over $\acl_{\CD}(t_0)$ and $b$ is in the frontier of one of those so we may replace $S_{t_0}$ by this
strongly minimal set, and modify the family $\CF$ accordingly. %Faithfulness of $\CF$ and the Morley rank and degree of $T$ assure that $S_{t_0}$ does not  coincides with a straight line. \\

\smallskip

 Denote $S=S_{t_0}$ and  $B=\Bad(b)$.
\\

 \noindent\textbf{Claim 1.} {\em $\dim(S\cap B)\leq 1$.}
\begin{proof}[Proof of Claim 1]
Since $\dim(t_0/\emptyset)\geq 6$ and $\dim(b/\emptyset)\leq 4$, we obtain
$\dim(t_0/b)\geq 2$. Assume towards a contradiction that $\dim(S\cap B)=2$. Let
$$I=\{t\in T :  \dim (S_t\cap B)=2\}.$$
Notice that $I$ is defined over $b$ and $t_0\in I$, so $\dim I\geq 2$. Because $\CF$
is almost faithful,  $\{S_t\cap B: t\in I\}$ is a  definable family of
$2$-dimensional subsets of $G^2$ satisfying the assumptions of Lemma
\ref{intersectionfinite}. It follows that $\dim \bigcup_{t\in T}(S_t\cap B)= 4$. But
$\bigcup(S_t\cap B)\sub B$, contradicting Lemma \ref{b-good}.
\end{proof}

%Recall that, for every $s\in G^2$, we let $Q(s) = \{q\in Q: s\in l_q\}.$

%(but not over $a$, since $b_0\in \fr(S_a)$!)

\noindent\textbf{Claim 2.} {\em For every $q'\in Q$, $S\cap l_{q'}$ is finite.}
%(it actually must be nonempty but this will not play any role here):
\begin{proof}[Proof of Claim 2]
If not, then by   strong minimality of $S$, we would have $S\sim l_{q'}$ for some
$q'\in Q$, implying -- since $S=S_{t_0}$ and $\CF$ is almost faithful -- that
$t_0\in \acl(q')$. However, we assumed that $\dim(t_0/\emptyset)\geq 6$, while
$\dim(q'/\emptyset)\leq 4$, a contradiction.
\end{proof}

 We fix an
element $q\in Q(b)$ generic over $t_0$ and $b$. Since $\dim(b/\emptyset)=4$, $q$ is
generic in $Q$ over $\emptyset$, hence
we have $\dim(q/\es)=4$.\\

Since $\CL$ is very ample, no two points in $G^2$ belong to infinitely many curves
in $\CL$, and hence each $s\in S\cap l_q$ is inter-algebraic with $q$ over $t_0$
and $b$. Thus such an $s$ is generic in $S$ over $t_0$
 and $b$.  So in particular $S$ is
smooth at $s$. It is not hard to see now (using the fact that  $\CF$ is almost
faithful) that $\dim(s/b)=4$.
\\

For the rest of this proof, we fix an element $s\in S\cap l_q$.\\

\noindent{\bf Claim 3.} {\em The curve $l_q$ is smooth at $s$, and the intersection
of $S$ and $l_q$ is transversal at $s$.}

\begin{proof}[Proof of Claim 3]
Because $\dim(s/b)=4$, it follows from Claim 1 that $s$ is $b$-good, so in
particular $l_q$ is smooth at $s$ and there exist neighborhoods $U\sub G^2$ of $s$
and $W\sub Q(b)$ of $q$, and a $\CD$-definable parameter choice function $g_b:
U\to W$, such that $g(s')$ is the unique $q'\in W$ with $s'\in l_q\cap U$.
Restricting $U, W$ if needed we may assume that $l_q\cap U$ (which equals
$g_b^{-1}(q)$) is a $C^1$-submanifold of $G^2$.
    Thus, the tangent space to $l_q$ at $s$, $T_s(l_q)$,  equals $\ker(d_s(g_b))$, where $d_s(g_b)$ is the differential of $g_b$ at $s$ viewed as a linear map between the tangent spaces (see Definition \ref{D:fibers}).   If the intersection is not
    transversal,
     then $\dim(T_s(l_q)\cap T_s(S))\geq 1$. It follows that $\dim(d_s(g_b)(T_s(S)))\le 1$,
      and by genericity of $s$ in $S$ over $t_0,b$, the same is true of any $s'\in S$ in some open neighborhood  $U'\ni s$.
      Thus, the image of $g_b(S\cap U)$ is a 1-dimensional manifold (or finite), and it follows that for some
      $q'$ in this image, $l_{q'}\cap S$ is infinite. This contradicts Claim 2.
\end{proof}

\noindent{\bf Claim 4.} {\em For every neighborhood $V\sub Q$ of $q$, there exists a
neighborhood $U\sub G^2$ of $b$ such that for every $b'\in U$ there are infinitely
many $q'\in V$ with $b'\in l_{q'}$.}

\begin{proof}[Proof of Claim 4]
By assumptions, $b\in l_q$ is generic in $G^2$ over $\emptyset$. Thus, by shrinking
$V$ if needed, we may assume $b$ is still generic in $G^2$ over the parameters
defining $V$. Since $Q(b)\cap V$ is infinite, the first order statement
$$\phi(b'):=(\exists^{\infty} q'\in V  )(b'\in l_{q'}) $$ holds for $b$ and therefore there
is a neighborhood $U\ni b$ for which it holds.\end{proof}

 Let $N$ be the number of intersection points of a
curve from $\CL$, generic over $t_0$, with $S$ (recall that $\md(\CL)=1$,
so $\CL$ has a unique generic type).\\

\noindent{\bf Claim 5.} The curve $l_q$ intersects $S$ in less than $N$ points.

\begin{proof}[Proof of Claim 5]
We write $l_q\cap S=\{s_1,\ldots, s_n\}$ (note that $b$ is not among them). We first
fix some open disjoint neighbourhoods $U_1,\ldots, U_n\sub G^2$, of $s_1,\ldots,
s_n$, respectively. By Claim 3 and Lemma \ref{transversal}, applied to each of the
$s_i$, there is a neighbourhood $V\sub Q$ of $q$ such that for every $q'\in V$, the
curve $l_{q'}$ intersects $S$ at least $n$ times -- at least once in each of the
$U_i$, $i=1,\ldots, n$. Next, we apply Claim 4 to $V$ and find $U_0\ni b$, which we
may assume is disjoint from all the $U_i$, as in Claim 4.

Because $b$ is in $\cl(S)\setminus S$, we can find in $S\cap U_0$ some $s'$, an
element $\CD$-generic over $t_0$, and by Claim 4, we can find in $V$ some $q'\in
Q(s')$ generic over $s'$ and $t_0$. But now $l_{q'}$ intersects $S$ at least $n+1$ times: at
$s'$ and in each of $U_1,\ldots, U_n$. Since $S\cap l_{q'}$ is finite, the curve
$l_{q'}$ is generic in $\CL$ over $t_0$. So we have $N\geq n+1>n=|l_q\cap S|$.
\end{proof}

Finally, let us see that $b\in \acl_{\CD}(t_0)$. The set $Y$ of all $q_1\in Q$ such that $|l_{q_1}\cap S|<N$ is $\CD$-definable over $t_0$ and has Morley Rank at most $1$. Since $q$ is generic in $Q(b)$ over $t_0$ and $b$, it follows from Claim 5 that $\mr (Q(b) \cap Y)=1.$ Also, there are at most finitely many $b$'s such that $\mr(Q(b) \cap Y )=1$, for otherwise there would be $b_1\ne b_2$ such that $Q(b_1)\cap Q(b_2)$ is infinite, contradicting the very ampleness of $\mathcal L$. Thus $b\in \acl_{\CD}(t_0)$. This ends the proof of Lemma \ref{FrgoodS}.\end{proof}
% EARLIER:
%Finally, let us see that $b\in \acl_{\CD}(t_0)$. By Claim 5 no generic curve in
%$\CL(b)$ intersects $S_{t_0}$ in a generic
% number of points. So $b$ is contained in the set $Y$ of all those $b'\in G^2$
% such that for all but finitely many $q_1\in \CL(b') $, we have  $|l_{q_1}\cap
% S|<N$. The set $Y$ is $\CD$-definable over $\es$ and has
%  Morley rank at most $1$. Since $t_0$ is generic in $T$ over $\0$ and $\mr(T)\ge 3$
%  we get that $\mr(t_0/\es)\geq 3$, and hence $Y\cap S_{t_0}$ is finite.
%   Since $b\in Y\cap S_{t_0}$ it follows that $b\in \acl_\CD(t_0)$.

%$|l_{q_1}\cap S|<N$
%The set $Y$ of all $q_1\in Q$ such that
%$|l_{q_1}\cap S|<N$ is $\CD$-definable over $t_0$ and has Morley rank at most $1$.
%By Claim 5, $Q(b)\setminus Y$ is finite. It follows that $b\in \acl_{\CD}(t_0)$. \qed

In our next step we show that the assumptions of Lemma \ref{FrgoodS} can be met for
a $\CD$-definable set $S$ of $\mr(S)=1$, after replacing $S$ by its composition with
a generic enough curve in a family $\mathcal L'$ as in Proposition \ref{familyL}.

\begin{lemma}\label{massageS}
Let $S\sub G^2$ be a $\CD$-definable strongly minimal set which is not a straight line, and assume that $c$ is
generic in $G^2$  over $\0$ and belongs to $\fr(S)$. Then there are
 \begin{enumerate} \item An almost faithful
stationary family of plane curves $\mathcal S'=\{S'_t : t\in T\}$, $\CD$-definable over $[S]$, with
$\mr(T)\geq 3$.

\item $t_0$ generic in $T$ over $c \cup \acl_\CD([S])$
%({\bf what do we mean by that? what is generic over $[S]$?}).

\item $b$ which is $\CD$-interalgebraic with $c$ over $t_0 \cup [S]$.

\item $b\in \fr(S'_{t_0})$ and $\dim(b/\emptyset)=4$.
\end{enumerate}

\end{lemma}

\begin{proof}  Let $\mathcal L'=\{C_t : t\in T\}$ be a $\CD$-definable family of plane curves as in Proposition \ref{familyL}.
Recall that, for every $(a, b)\in G^2$,
$$T(a, b):= \{t\in T: (a, b)\in C_t\}.$$
If we write  $c=(c_1, c_2)$, then by assumption, $c_2$ is generic in $G$ over $\0$.
Fix an element $b_2\in G$, which is generic over $c_2 \cup \acl_\CD([S])$ (abusing
notation, in the present proof we will write $[S]$ for $\acl_\CD([S])$), and let
$t_0$ be generic in
 $T(c_2, b_2)$ over $c_1, c_2, b_2$ and $[S]$. Note that $(c_2,b_2)\in G^2$ is generic and $C_{t_0}$ is generic
 through it. So $\dim(t_0c_2b_2)=\dim(T)+2$, whereas $\dim(t_0/c_2b_2)=\dim(T)-2$. Since $b_2\in \acl_\CD(t_0c_2)$,
 we get that $\dim(t_0/c_2)=\dim(T)$. Because $t_0$ was chosen generic over $c_1,[S]$ too, we get
 $\dim(t_0/c_1c_2[S])=\dim(T)$.
%  $c_1$ is generic over $c_2$ Thus one can show ({\bf ???}) that $\dim(t_0/c_1,c_2,[S])=\dim T$.

%Consider the composition $\CL\circ S$. By Lemma \ref{L:almost-faithful-comp} it extends a $\CD$-definable almost faithful family $\CS'$ of the same dimension.
%Replacing $\CS'$
%by a $\CD$-conjugate (over $[S]$) composition subfamily, if needed, we may assume that the strongly minimal subset of $C_{t_0}\circ S$ is contained in $\CS'$. We denote it (abusing
%notation) $S_{t_0}'$.

%We define
%$$\mathcal S'=\{C_t\circ S:t\in T\},$$ clearly defined in $\CD$ over $[S]$. Because
%$\mathcal L'$ was almost faithful, so is the family $\mathcal S'$, and we still have
%$\mr(T)\geq 3$. We  consider the curve $S'_{t_0}=C_{t_0}\circ S$.

   We set
  $b:=(c_1, b_2)$.  Since $(c_2, b_2)\in C_{t_0}$ and $\mr(C_{t_0})=1$, $b_2$ and $c_2$
  are inter-algebraic in $\CD$ over $t_0$ and $[S]$,
  and hence so are $(c_1, b_2)$ and $(c_1, c_2)$. \\
  %This  establishes property (1)-(3), together with $\dim(b/\emptyset)=4$.
  %  It is therefore sufficient to prove:\\

\noindent\textbf{Claim.} $b\in \fr(C_{t_0}\circ S)$.
\begin{proof}[Proof of Claim] Since
$c_2$ is generic in $G$ over $\emptyset$, $(c_2,b_2)$ is generic in $G^2$ over
$\emptyset$ and therefore, by our choice of $t_0$, the point $(c_2,b_2)$ is also
 generic in $C_{t_0}$ over $t_0$. Hence, the curve $C_{t_0}$ is a homeomorphism at $(c_2,b_2)$.
 Denote this local map by $f_{0}$.   It follows that  the map $(x,y)\mapsto (x,f_{0}(y))$ is a local homeomorphism
on a neighborhood $W$ of $(c_1,c_2)$, sending $(c_1,c_2)$ to $(c_1,b_2)$. It is easy
to verify that it sends every point in $S\cap W$ to a point in $C_{t_0}\circ S$, and
therefore sends every point in $\cl(S)\cap W$ to a point in $\cl(C_{t_0}\circ S)$.
We conclude that $(c_1,b_2)\in \cl(C_{t_0}\circ S)$.

It remains to see that $(c_1,b_2)\not\in C_{t_0}\circ S$. Let
$$S_{c_1}=\{y\in G:(c_1,y)\in S\}=\{d_1,\ldots, d_k\}.$$ Note that, since  $(c_1,c_2)\notin S$, we
 have $c_2\notin S_{c_1}$. Also, $(c_1,b_2)\in C_{t_0}\circ S$ if and only if there
 is some $i=1,\ldots, k$ for which $(d_i,b_2)\in C_{t_0}$.

 Since $\CL'$ is very ample, for every $i=1,\ldots, k$, $\dim(T(c_2,b_2)\cap T(d_i,b_2))<\dim T$. But $t_0$
  is generic in $ T(c_2,b_2)$  over
$\{c_1,c_2,d_1,\ldots, d_k, b_2, [S]\}$ and therefore, $t_0\notin T(c_2,b_2)\cap
T(d_i,b_2)$. That is, none of the points
 $(d_i,b_2)$ are in $C_{t_0}$. It follows that $(c_1,b_2)\notin C_{t_0}\circ S$, so
 we may conclude that $b=(c_1,b_2)\in \fr(C_{t_0}\circ S)$.
\end{proof}

Since $S$ is not a straight line, we have $\mr(C_{t_0}\circ S)=1$, and hence there is a strongly minimal $C\subseteq C_{t_0}\circ
S$ such that $b\in \fr(C)$. By Lemma \ref{L:almost-faithful-comp},
$\mr[C]=\mr[C_{t_0}]=\mr(T)$ and is contained, therefore, in an almost faithful
family $\CS'$ of the same rank. This gives condition (1) of the lemma, (2) is by the
choice of $t_0$, (3) is the line before the above claim, and (4) is what we just
showed. So the lemma is proved.
%This ends the proof of Lemma \ref{massageS}.
\end{proof}

We can now conclude the main result of this section.

\begin{theorem}\label{main-frontier}
Let $S\sub G^2$ be a $\CD$-definable set with $\mr(S)=1$. Then $\fr(S)\sub
\acl_{\mathcal D}([S])$ and hence $\fr(S)$ is finite. In particular, $S$ is locally
closed, namely every $p\in S$ has a neighborhood $U\ni p$ in $G^2$ such that $S\cap
U$ is closed in $U$.
\end{theorem}
\begin{proof}   Since $\mr(S)=1$, $S$ can be written as $\bigcup\limits_{i=1}^k S_i$ for some strongly minimal sets  $\CD$-definable over $\acl_\CD([S])$. Since $\fr(\bigcup\limits_{i=1}^k S_i)\subseteq \bigcup\limits_{i=1}^k \fr(S_i)$, it suffices to prove the theorem for $S$ strongly minimal. Moreover, if $S$ is a straight line, then clearly its frontier is contained in finitely many points, which are in $\acl_\CD([S])$. So we may assume that $S$ is strongly minimal not coinciding with any straight line.

Fix $c\in \fr(S)$. Replacing $S$ by $S+p$ for $p$ generic in $G^2$ over $c$ and $[S]$,
 we may assume that $\dim(c/\emptyset)=4$. We can now apply Lemma \ref{massageS} and obtain $t_0$, $S'_{t_0}$ and $b\in
\fr(S'_{t_0})$ as in the lemma. Working first in a richer language where $[S]$ is
$\0$-definable, we may apply Lemma \ref{FrgoodS}, and  then conclude that $b\in
\acl_{\CD}(t_0,[S])$.

By Lemma
 \ref{massageS}, $c$ is interalgebraic with $b$ over $t_0$ and $[S]$ hence, $c\in \acl_\CD(t_0,[S])$.
  Since $\dim(t_0/c,[S])=\dim(t_0/c)$, we obtain that $c\in \acl_\CD([S])$.

For $p\in S$, let $U\ni p$ be any neighborhood such that $U\cap \fr(S)=\0$, and then
$S\cap U$ is closed in $U$.
\end{proof}

\subsection{Two structural corollaries on plane curves}

The first corollary will be used in the next subsection.

\begin{corollary}\label{C:uniform-cl}
    Let $\CL$ be a family of plane curves. Assume $\CL$ is $\CD$-definable over $\0$. Then there exists a family of plane curves $\CL'$, also $\CD$-definable over $\0$, such that:
    \begin{enumerate}
        \item Every curve in $\CL'$ is closed.
        \item For every curve $X_s\in \CL$, there exists a curve $X_s'$, defined over the same parameters, such that $X_s\sim X_s'$.
        \item For every $X_s'\in \CL'$, there exists $X_s\in \CL$, defined over the same parameters, such that $X_s'\sim X_s$.
    \end{enumerate}
\end{corollary}
\begin{proof}
Let $\chi(x,y)$ define $\CL$ and $\psi(y):=(\exists^\infty x)\chi(x,y)$. We prove the corollary by induction on  $(\mr(\psi), \md(\psi))$. For $\mr(\psi)=0$, the corollary is Theorem \ref{main-frontier}.

In the general case, fix $s\models \psi(y)$ generic. By definition, $[X_s]\in \dcl_\CD(s)$. By Theorem \ref{main-frontier}, there is a finite set $R_s$, $\CD$-definable over $[X_s]$, so a fortiori also over $s$, such that $\fr(X_s)\subseteq R_s$. Let $\phi(x,s)$ define $R_s$. By compactness, there is a formula $\theta(y)\in \tp(s)$  such that for all $r\models \theta$ the formula $\chi(x,r)$ is algebraic and, if not empty, its set of realisations contains $\fr(X_r)$. So, for all $r\models \theta$, the formula $\phi(x,r)\lor \chi(x,r)$ defines a closed plane curve $\sim$-equivalent to $X_r$. Because $\theta(y)\in \tp(s)$ and $X_s$ was generic in $\CL$, we get that $\psi(y)\land \neg \theta(y)$ has smaller $(\mr, \md)$ (in the lexicographic order) than $\psi(y)$. So we are done by the induction hypothesis.
\end{proof}

The second corollary below  will be used several times in the rest of the paper.
\begin{definition}
    Let $S\sub G^2$  and $a=(a_1,a_2)\in S$.
    We say that {\em $S$ is injective at $a$ over $a_1$} if there is an open
    neighborhood $U_1\times U_2\sub G\times G$ of $a$ such that for every $y\in U_2$ there exists at most
    one $x\in U_1$ such that $(x,y)\in S$. Namely, $S\cap (U_1\times U_2)$ is the graph of a function from a subset of $U_2$ into $U_1$.
        We say that {\em $a$ is an injective point of $S$} if $S$ is injective at $a$
    over $a_1$ and $S^{\textup{op}}$ is injective at $(a_2, a_1)$ over $a_2$. Otherwise, we say that $a$ is {\em a non-injective point of $S$}.

    Let $S\sub G^2$ and $a_1\in G$. We say that {\em $S$ is injective over $a_1$} if for every $a=(a_1,a_2)\in S$, the set $S$ is injective at $a$ over $a_1$.

    % Say that $x\in \pi(S)$ is a
    % {\em an injective point of $S$} if for all $y\in S_x:=\{z: (x,z)\in S\}$,
    % there is an open $U\ni (x,y)$ such that $S\cap U$ is the graph of an injective function in the $x$ variable
    % (we do not specify its domain).

    % We call $x\in \pi(S)$ {\em a non-injective point of $S$} otherwise.
\end{definition}

Note that $S$ is injective at every isolated point. Also, we cannot yet rule out the
possibility that $a$ is an injective point of $S$ belonging to a $1$-dimensional component of $S$.

%We now prove the following corollary to Theorem  \ref{main-frontier}.

\begin{corollary}\label{injective}
    Let $S\subseteq G^2$ be a $\CD$-definable strongly minimal set. If $S$ is not $\sim$-equivalent to any fiber
    $G\times \{a\}$, then the set of $x\in G$ such that $S$ is non-injective over $x$  is finite and contained in
    $\acl_{\CD}([S])$.  If $S$ is not a straight line, then the set of non-injective points of $S$ is finite and contained in $\acl_{\CD}([S])$.
\end{corollary}
\begin{proof} By Theorem \ref{main-frontier}, we may assume that $S$ is closed.
Let
\[S_1=\{(x_1,x_2)\in G^2: x_1\neq x_2 \,\, \&\,\, \exists y ((x_1,y)\in S \wedge
 (x_2,y) \in S)\}=(S^{\textup{op}}\circ S)\setminus \Delta.\]
The set $S_1$ is $\CD$-definable over the same parameters as $S$. Since     $S$ is not $\sim$-equivalent to any fiber
    $G\times \{a\}$, we have  $\mr(S_1)\leq 1$. Note
 that $(x,x)\notin \fr(S_1)$ if and only if there exists an open $U\ni x$ such that
for all $y\in G$ there exists at most one $x'\in U$ such that $(x',y) \in S$.
It follows that $(x,x)\notin \fr(S_1)$ if and only if $S$ is injective over $x$.
 By Theorem \ref{main-frontier}, $\fr(S_1)\sub \acl_{\CD}([S])$ thus the set of  $x\in G$ such that $S$ is non-injective over $x$  is finite and contained in
    $\acl_{\CD}([S])$.

The second clause follows immediately by applying the first one also to $S^{\textup{op}}$.
\end{proof}
%Similarly,  let
%\[S_2:=\{(x_1,x_2)\in G^2:x_1\neq x_2\,\, \& \,\, \exists y ((x_1,y)\in S \wedge (x_2,y)\in S)\}=(S^{\textup{op}}\circ S)
%\setminus \Delta.\] Then $(b,b)\notin \fr(S_2)$ if and only if there exists an open
%$V\ni b$ such that for all $y\in G$ there exists at most one $x\in V$ such that
%$(x,y)\in S$. Equivalently, $S^{\textup{op}}\cap (G\times V)$ is the graph of a function.

%Thus, for $(a,b)\in S$, if $(a,a)\notin \fr(S_1)$ and $(b,b)\notin \fr(S_2)$ then
%there is an open $U\times V\ni (a,b)$ such that $S\cap U\times V$ is the graph of an
%injective function (notice that we cannot determine yet the domain of the function).
%It follows that if $(a,b)\in S$ is a non-injective point then either $(a,a)\in
%\fr(S_1)$ or  $(b,b)\in \fr(S_2)$. By Theorem \ref{main-frontier}, these frontiers
%are finite and contained $\acl_{\CD}([S])$.

%Finally,  by our assumptions on $S$, for every $a,b\in G$ the sets $S_a$ and $S^b$
%are finite  and therefore the set of $(a,b)$ which are non-injective is finite and
%contained in $\acl_{\CD}([S])$.

\subsection{On $\CD$-functions}\label{D-functions}

%\textbf{Let us put some introductory paragraph.}

Every plane curve $S\sub G^2$ gives rise to  a definable partial function from $G$
into $G$   around almost every point in $S$ (except when $S$ is  contained in
finitely many straight lines $\{a\}\times G$). The goal of this subsection is to establish the
basic theory of such functions.

\begin{defn} Let $U\sub G$ be a definable open set and  $f:U\to G$ be a definable
continuous function. \begin{enumerate}
        \item We say that $f$ is {\em a $\CD$-function} if there exists a plane curve $S\sub G^2$ such that
        $\Gamma_f\sub S$. We say in this case that {\em $S$ represents $f$}.
        \item We say that $f$ is \emph{$\CD$-represented over $A$} if is there exists
        $S$ representing $f$ which is $\CD$-definable over $A$.
        \item We say that a plane curve $S$ {\em represents the germ of  $f$ at $x_0\in U$} if there exists an
        open neighborhood $W\ni x_0$, $W\sub \dom(f)$, such that $\Gamma_{f|_W}\sub S$.
    \end{enumerate}

\end{defn}

Note that our definition does not require that $S$ is, locally at $(x_0, f(x_0))$,
the graph of a function, but only that it {\em contains} the graph of $f$. Indeed,
at least for some of the $\CD$-functions we need to consider we do not know whether
this stronger property can be achieved as well.

\begin{lemma} \label{rep-strongly} Let $U\sub G$ be a definably connected open set and
$f:U\to G$ a continuous $\CD$-function, $\CD$-represented over $A$. Then $f$ can be
$\CD$-represented over $\acl_{\CD}(A)$ by a strongly minimal set.
%$S\sub G^2$, definable in $\CD$ over

\end{lemma}
 \begin{proof}
  Assume that $f:U\to G$ is $\CD$-represented over $A$ by $S$.  We let $S=S_1\cup\cdots \cup S_r$ be
  a decomposition of $S$ into strongly minimal
sets, definable in $\CD$ over $\acl(A)$. By Theorem \ref{main-frontier}, we may
assume, by adding finitely many points in $\acl_{\CD}(A)$, that each $S_i$ is closed
in $G^2$, but now the intersection $S_i\cap S_j$ for $i\neq j$ may be non-empty and
finite. We claim that one of the $S_i$ must contain $\Gamma_f$. Indeed, for each
$i=1,\ldots r$, let $C_i=\pi(S_i\cap \Gamma_f)\sub U$, where $\pi:G^2\to G$ is the
projection on the first coordinate. By the continuity of $f$,  these are definable,
relatively closed subsets of $U$, whose pair-wise intersection is at most finite.

Let $U':=U\sm \bigcup\limits_{i\neq j}C_i\cap C_j$. Because $U$ is open and
definably connected so is $U'$. For $i=1,\ldots, r$ let $C_i'=C_i\cap U'$. The
$C_i'$'s are pairwise disjoint and still relatively closed in $U'$. So each $C_i'$ is clopen (having a closed complement) in $U'$ so for some $j$, $C_j'=U'$. Because $C_j$ is closed in
$U$ it follows that $C_j=U$.\end{proof}

%\begin{defn} We let $\mathfrak F$ be the collection of all smooth $\CD$-functions
%$f$ in a neighborhood of $0$, with $f(0)=0$.

%We say that a strongly minimal set {\em represents $f\in \mathfrak F$ } if it
%represents the germ of $f$ at $0$.
%\end{defn}

\begin{proposition}\label{rep-Dfunction} Let $\{S_t:t\in T\}$ be a family of plane curves. Assume that this family is $\CD$-definable over $A$, and that for every $t\in T$, $(0,0)\in S_t$.
 Then there exists a  family $\CF=\{f_s: s\in T_0\}$, definable (in $\CM$) over $A$, of functions in $\mathfrak F$ (defined  in Section \ref{sec-overview}),  such
that:
\begin{enumerate}
\item For every $t\in T$, if $S_t$ represents the germ at $0$ of a $\CD$-function
$f\in\mathfrak F$, then there exists $s\in T_0$ and an open $W\ni 0$ such that
$f|_W=f_s|_W$. \item For every $s\in T_0$ there exists $t\in T$ such that $S_t$
represents the germ at $0$ of $f_s$.
\end{enumerate}
\end{proposition}
\begin{proof}  By Corollary \ref{C:uniform-cl}, there exists a $\CD$-definable family $\CL'$ of closed plane curves such that each curve in $\CL$ is $\sim$-equivalent to one in $\CL'$ and vice versa.
%
%finite sets $\{R_t:t\in T\}$,  such that for every $t\in T$,  $\cl(S_t)\sub S_t\cup
%R_{t}$. Since $R_{t}$ is finite it follows that $S_t\cup R_{t}$ is closed, so we can
%replace the original family with $\{S_t\cup R_t:t\in T\}$ (note that we do not claim
%that $R_t=\fr(S_t)$).
Note that whenever $S_t$ represents the germ of a $\CD$-function $f_t$ at $0$,
%(so, since $\CD$-functions are continuous, $S_t$ is locally closed at $(0,0)$),
if $S'_t\in \CL'$ is $\sim$-equivalent to $S_t$, then it also represents the germ of $f$ at $0$.  Thus, we may replace $\CL$ with $\CL'$ and assume that every curve $S_t$ is closed.

By fixing a coordinate system near $0$, we can identify some neighbourhood  $W\ni 0$
in $G$ with an open subset of $R^2$. For each $r>0$, we consider the disc $B_{r}$
centered at $0$, and let $S_t^{r}=S_t\cap (B_{r}\times W)$. By o-minimality, there
exists a uniform cell decomposition of the sets $\{S_t^{r}:t\in T,r>0\}$. In
particular, there is a bound $k\in \mathbb N$ such that every such decomposition
contains at most $k$ cells. By allowing cells to be empty, we obtain a definable
collection of cells $\{C^{r}_{t,i}:t\in T\, ,\, r>0, i=1,\ldots, k\}$, such that for
every $t\in T$, $r>0$,
$$S_t^{r}=\bigcup_{i=1}^k C^{r}_{t,i}.$$
Recall that the notion of a decomposition implies that for $C^{r}_{t,i},
C^{r}_{t,j}$, if $\pi:G^2\to G$ is the projection onto the first coordinate, then
either $\pi(C^{r}_{t,i})=\pi(C^{r}_{j,j})$ or $\pi(C^{r}_{t,i})\cap
\pi(C^{r}_{j,j})=\emptyset$.\\

\noindent\textbf{Claim.} \emph{For every $t\in T$, and a $\CD$-function $f\in
\mathfrak F$, the following are equivalent:
\begin{enumerate}
\item $S_t$ represents the germ of $f$ at $0$. \item There exist $r>0$, and $A\sub
\{1,\ldots, k\}$, such that
 $$\Gamma_{f|_{B_{r}}}=\bigcup_{i\in A}C^{r}_{t,i}.$$
\end{enumerate}}
\begin{proof}[Proof of Claim] $(1)\Rightarrow (2)$. We  assume that $S_t\cap (B_r\times G)$ contains the
graph of  $f|_{B_{r}}$, for $r>0$. To simplify notation we omit $r$ and consider the
cell decomposition $S_t=C_{t,1}\cup\cdots \cup C_{t,k}$.

We let $A\sub \{1,\ldots, k\}$ be all $i$ such that $C_i\cap \Gamma_f\neq
\emptyset$. We fix a cell $C=C_i$ with $i\in A$ and claim that $C\sub \Gamma_f$.
Without loss of generality $\dim C>0$, and since $C\sub \Gamma_f$, the projection
$\pi:C\to G$ is injective. Since $C$ is definably connected, it is sufficient to
prove that $C\cap \Gamma_f$ is clopen inside $C$. Because $C$ is locally closed and
$f$ is continuous, it follows that $C\cap \Gamma_f$ is closed in $C$, so we need to
prove that it is also open in $C$.

%Let $U\sub B_{r}$ be an open neighborhood of $x_0$, and consider $U\cap \pi(C)$.
Fix some $x_0\in \pi(C\cap \Gamma_f)$. Since $\Gamma_f\sub S_t$, there exists a cell $C'$ in the decomposition of $S_t$
containing $\Gamma_f\cap [(U\cap \pi(C))\times G]$ for some open set $U\ni x_0$.
But then $\pi(C)\cap \pi(C')\neq \emptyset$ and therefore $\pi(C)=\pi(C')$. By the
continuity of $f$, it follows that $(x_0,f(x_0))\in C'$, forcing $C'=C$. It follows
that $C\cap \Gamma_f$ is clopen in $C$, and therefore $C\sub \Gamma_f$.

We showed that for each $i\in A$, $C_i\sub \Gamma_f$ and hence
$\Gamma_f=\bigcup_{i\in A}C_i.$

$(2)\Rightarrow (1)$. This is immediate, since $\Gamma_{f|_{B_r}}\sub S_t$.\end{proof}

We now return to the proof of Proposition \ref{rep-Dfunction} and consider the
uniform decomposition
$$S_t^{r}=\bigcup_{i=1}^k C^{r}_{t,i}.$$
%For each $t\in T$ and $r>0$, we let $A^{r}_t\sub \{i=1,\ldots,k\}$ be all $i$ such that $\dim(C^{r}_{t,i})=2$.

For each $A\sub \{1,\ldots,k\}$, we consider
$$G_{t,A}^{r}=\bigcup_{i\in A} C^{r}_{t,i}.$$
The family
$$\CF=\{G_{t,A}^{r}: G_{t,A}^{r} \mbox{ is the graph of a continuous function on $B_{r}$}\}$$
is definable in $\CM$, as $t$ varies in $T$, $A$ varies among subsets of
$\{1,\ldots,k\}$ and $r>0$. By the above claim, this family satisfies our
requirements.
\end{proof}

\begin{remark}
\begin{enumerate}
    \item Note that in the above family $\CF$ of $\CD$-functions, each germ of a function
    appears infinitely often since we allow arbitrarily small $r$. One can divide the
    family, definably in $\CM$,  by the equivalence of germs at $0$ and then, using
    Definable Choice in o-minimal structures, obtain a unique $\CD$-function in the
    family representing each germ. Thus, if $f\in \mathfrak F$ is represented by the plane
    curve $S_t$, then there exists  $g\in \mathfrak F$ which has the same germ as $f$
    at $0$ and is definable in $\CM$ over $t$.
    \item It follows from the above that if $S_t$ represents $f\in \mathfrak F$, then $J_0(f)$ is in $\dcl(t)$ (recall from Section \ref{sec-overview} that $J_0f$ is the Jacobian of $f$ at $0$ with respect to some fixed differential structure on $G$).
\end{enumerate}
\end{remark}

\noindent{\bf Notation.} For a $\CD$-function $f$, we reserve the notation $S_f$ for
a strongly minimal set representing $f$. Note that $S_f$ is unique only up to $\sim$-equivalence. \\

We conclude this section with an open mapping theorem for $\CD$-functions.

\begin{theorem}\label{C:continuous}
        Let $U\sub G$ be an open definably connected set and $f:U\to G$ a continuous non-constant $\CD$-function. Then $f$ is an open map.
\end{theorem}
    \proof  By Lemma \ref{rep-strongly}, there exists a $\CD$-definable strongly minimal $S_f\sub G^2$
    representing $f$. Because $f$ is not constant, the projection of $S_f$ onto both
    coordinates is finite-to-one, so it is not a straight line. By Corollary \ref{injective}, $S_f$ is injective at co-finitely many points, and therefore so is also $f$. By the o-minimal version of Brouwer's invariance of domain \cite{Johns}, it follows that $f$ is open at every injective point of its domain. So $f$ is open after possibly removing finitely many points from its domain. It is easy to check (see, for example, the proof of \cite[Proposition 4.7]{HaKo} for  details) that a function which is continuous on a disc and open on the punctured disc is open on the whole disc. So $f$ is open.

%
%
%   By o-minimality, since $f$ is continuous and not open, there are infinitely many
%   points $a\in U$, such that $f$ is not open at $a$. This contradicts Lemma \ref{open-rel1}.
%} Because $f$ is not constant, the projection of $S_f$ onto both
%coordinates is finite-to-one.
%By o-minimality, since $f$ is continuous and not open, there are infinitely many
%points $a\in U$, such that $f$ is not open at $a$. This contradicts Lemma \ref{open-rel1}.

\section{Poles of plane curves} \label{sec-poles}

Recall that we assume that $G$ is a definable closed subset of some $R^n$,
 equipped with
the subspace topology, making it a topological group.

The goal of this section is to prove that just like affine algebraic curves in
$\mathbb C^2$, every plane curve has at most finitely many poles. We may assume that
$0_G=0\in R^n$.  For $x\in G$ and $\epsilon>0$ in $R$, we write
$$B(x;\epsilon)=\{g\in G: |x-g|<\epsilon\}, $$
and $B_{\epsilon}$ for $B(0;\epsilon)$.
 For $A\sub G$, and $\epsilon>0$, we let
$$B(A;\epsilon)=\{y\in G:\exists x\in A\,\, , \, y\in B(x;\epsilon)\}.$$

In this section, we will also consider definable curves, that is, definable maps $\gamma:(0,1)\to U\sub R^n$, which we will denote, for simplicity, by $\gamma(t)\in U$. We recall from   \cite[\S6.1]{vdDries} that if $x\in \cl(X)$, for some definable $X\sub R^n$, then by curve selection for $\mathcal M$, there is a definable path $\gamma(t)\in X$ with $\lim\limits_{t\to 0} \gamma(t)=x$.

\begin{defn}
Let $S\sub G^2$ be a definable set. We call $a\in G$  a {\em pole of $S$} if for
every open $U\sub R^n$ containing $a$,   the set $(U\times G)\cap S$ is an unbounded
subset of $R^n$. We denote the set of poles of $S$ by $\UB{S}$.
\end{defn}

Given $S\sub G^2$ and $U\sub G$, we define
$$S(U):=\{y\in G:\exists x\in U\,\, (x,y)\in S\}\sub G.$$
Note  that  then $a\in \UB{S}$ if and only if for every open $U\sub G$ containing
$a$, $S(U)$ is  unbounded.
%(equivalently, the closure of $S(U):=\{y\in G:\exists x\in U\,\, (x,y)\in S\}$ in $G$ is not definably compact with respect to the group topology).
Another remark is that if $S$ is $G$-affine, then $\UB{S}=\es$. Indeed, if $S$ is a
subgroup of $G^2$ or its coset,  then its projection onto the first coordinate is a
finite-to-one topological covering map, and hence $S$ has no poles.

\smallskip The main result of this section is the following.

\begin{theorem} \label{inf-thm} If $S\sub G^2$ is a $\CD$-definable set and $\mr(S)=1$, then $\UB{S}$ is finite.\end{theorem}

%{\bf Question: can we show that each point in $\UB{S}$ is in $acl_{\CD}([S])$?}
%\proof

Notice that if $G$ is a definably compact group (for example, a complex elliptic
curve), then $G$ is a closed and bounded subset of $R^n$, and hence
$\UB{S}=\emptyset$. So the theorem is of interest for those $G$ which are not definably
compact.

Let us first introduce the key notion of ``approximated points'' and then discuss the strategy of our proof. %Given sets $S\sub G^2$ and $I\sub G$, we introduce the set of
%all points in $G$ ``approximated by $S$ near $I$'', denoted by $A(S, I)$.
Recall that for $S\sub G^2$ and $x\in G$, we let $S_x=\{y\in G:(x,y)\in S\}$.

\begin{defn}\label{A(S)} Let $S\sub G^2$, $b,x_1,x_2\in G$, and $I\sub G$. We say that
\begin{enumerate}
\item \begin{itemize} \item {\em $b$ is $S$-attained at $(x_1,x_2)$ } if $b\in
S_{x_1}-S_{x_2}$. \item {\em $b$ is $S$-attained in $I$} if it is $S$-attained at some $(x_1,x_2)\in I$. %there are $x_1,x_2\in I$
%such that $b\in S_{x_1}-S_{x_2}$.
\end{itemize}
\item \begin{itemize} \item {\em $b$ is $S$-attained near $(x_1,x_2)$} if for every
$\epsilon>0$ $b$ is $S$-attained in $B((x_1,x_2); {\epsilon})$.
%there are $x_1'\in B(x_1;\epsilon)$, $x_2'\in B(x_2;\epsilon)$ such
%that $b\in S_{x'_1}-S_{x'_2}$.
 \item {\em $b$ is $S$-attained near $I$} if for every $\epsilon>0$, $b$ is $S$-attained
at $B(I;\epsilon)$.
\end{itemize}
\item \begin{itemize} \item {\em $b$ is $S$-approximated near $(x_1,x_2)$} if for
every $\epsilon>0$, some $b'\in B(b;\epsilon)$ is $S$-attained at $(x_1,x_2)$.
%are $x_1'\in B(x_1;\epsilon)$, $x_2'\in B(x_2;\epsilon)$
%such that $B(b;\epsilon)\cap (S_{x_1'}-S_{x_2'})\neq \emptyset$.
\item {\em $b$ is
$S$-approximated near $I$} if for every $\epsilon >0$, there are $x_1, x_2\in
B(I;\epsilon)$ such that $B(b;\epsilon)\cap (S_{x_1}-S_{x_2})\neq \emptyset$. The
set of such points $b$ is denoted by $A(S, I)$.
\end{itemize}
\end{enumerate}
We omit $S$ from the above notation whenever it is clear from the context.
\end{defn}

The following claim is immediate from the definitions.
\begin{claim} For any $S\sub G^2$ and $I\sub G$,
$$\mbox{$b$ attained at $I$}\,\Rightarrow\, \mbox{$b$ is attained near $I$}\,
\Rightarrow\, \mbox{$b$ is approximated near $I$}.$$ If, in addition, $S$ and $I$ are
closed and bounded, then the above notions are  equivalent and $A(S,I)=S(I)-S(I)$.
\end{claim}

Here is a simple example.

\begin{example} Let $G=\la \mathbb
C,+\ra$ and consider the complex algebraic curve $$S=\{(z,w)\in \mathbb
C^2:zw=1\}.$$ The following are easy to verify: $\UB{S}=\{0\}$, every $b\in \mathbb
C$ is attained near $0$, and thus $A(S, \{0\})=\bb C$.
\end{example}

 The strategy of the proof of Theorem
\ref{inf-thm} is as follows. Assume towards a contradiction that the theorem fails.
It is easy to see that we may assume that  $S$ is closed,  strongly minimal and not
$G$-affine. Now, for any such $\CD$-definable set $S$ and infinite definable $I\sub
G$, we first find an infinite definable set $I_0\sub I$ and an open bounded ball
$B\sub R^n$, such that the set $A(S, I_0)\sm B$ is at most $1$-dimensional
(Proposition \ref{attained3}(1)). Then, using further that $\UB{S}$ is infinite, we
construct (Proposition \ref{lower}) another $\CD$-definable set $\hat S$, again
closed, strongly minimal and not $G$-affine, and an infinite definable $\hat I\sub
G$, such that for every infinite definable set $T\sub \hat I$ and open bounded ball
$B$, the set $A(\hat S, T)\sm B$ is $2$-dimensional. This gives the desired contradiction.

\subsection{An upper bound on the dimension of the set of approximated points}

The goal of this subsection is to prove the following proposition.
\begin{proposition}\label{attained3}
Assume that $S\sub G^2$ is a $\CD$-definable strongly minimal closed set which is
not $G$-affine, and let $I\sub G$ be an infinite definable set. Then there is
 a definable $1$-dimensional  $I_0\sub I$,
 such that
\begin{enumerate}
\item there exists a bounded $B\sub G$, such that the set $A(S, I_0)\sm B$ is at
most $1$-dimensional,

\item for every definable open $V\ni 0$  in $G$ there exist $\epsilon>0$ and a bounded ball
$B'\ni 0$ such that for all $x\in G\setminus B'$,
$$x+V\nsubseteq S(B(I_0,\epsilon)).$$

\end{enumerate}
\end{proposition}

The rest of this subsection is devoted to the proof of Proposition \ref{attained3}.
We fix throughout $S$ as in its assumptions. Since $\dim S=2  \mr(S)=2$, it follows easily from cell decomposition for $\mathcal M$ that $\dim \UB{S}\le
1$.  Absorbing $[S]$ into the language, we assume that $S$ is $\CD$-definable over $\0$.

We begin with an observation regarding the notions of Definition \ref{A(S)}.

\begin{lemma}\label{attained} Let $I\sub G$ be a definable bounded set over $\0$.
\begin{enumerate}

\item If $b\in G$ is attained near $I$, then there are $x_1,x_2\in \cl(I)$ such that
$b$ is attained near $(x_1,x_2)$.

\item If $b\in G$ is generic over $\0$ and $b$ is approximated near $I$, then $b$
is attained near $I$.
\end{enumerate}
\end{lemma}
 \begin{proof} (1) By assumption, and by curve selection in $\mathcal M$, there are definable curves $x_1(\epsilon)$, $x_2(\epsilon)$, $y_1(\epsilon)$, $y_2(\epsilon)\in G$, such that for every $\epsilon$, we have $x_1(\epsilon), x_2(\epsilon)\in B(I;\epsilon)$, $(x_i(\epsilon),y_i(\epsilon))\in S$,
$i=1,2$, and $b=y_1(\epsilon)-y_2(\epsilon)$. Since $I$ is bounded, the curves
$x_i(\epsilon)$ have limits $x_1,x_2\in \cl(I)$, so $b$ is attained near $(x_1,x_2)$.
%Note that if $y_1(t), y_2(t)$ are bounded, then their limits $y_1,y_2$ will satisfy $y_1-y_2=b$ and since $S$ is closed, also $(x_1,y_1),(x_2,y_2)\in S$, and therefore $b$ is
%attained at $x_1,x_2\in X$.

   (2) Fix $b$ generic in $G$ over $\0$, and assume that it is approximated near $I$.
It follows from the definition that for every $\epsilon>0$,  the element $b$ is in
the closure of
$$Y_{\epsilon}=\{y_2-y_1:\exists x_1,  x_2\in B(I;\epsilon)\,\, (x_1,y_1)\,,\,
(x_2,y_2)\in S\}.$$

Notice that the collection of $Y_{\epsilon}$ forms a definable chain of definable
sets decreasing with $\epsilon$.  We may now take $\epsilon$ sufficiently small, so
that $b$ is still generic in $G$ over $\epsilon$, and therefore $b$ is generic in
$\cl(Y_\epsilon)$ over $\epsilon$. Hence,  $b\notin \fr(Y_{\epsilon})$, a set of
dimension at most $1$. It follows that  $b\in Y_{\epsilon}$ for all sufficiently
small $\epsilon$, and so  $b$ is attained near $I$.\end{proof}

The following technical claim about definable and $\CD$-definable sets will be used
in the subsequent lemma.
\begin{claim} \label{dimension1} Let $P\sub G^2\times G$ be a $\CD$-definable set of Morley rank $1$
whose projection on the $G^2$-coordinate is finite-to-one. Then for any definable
sets $I,J\sub G$ of dimension at most $1$,  $$\dim(P\cap (I\times J\times G))\leq
1.$$
\end{claim}
\begin{proof}

Since the projection $\pi:P\to G^2$ is finite-to-one,
\[
\dim(P\cap (I\times J\times G))=\dim(\pi(P\cap (I\times J\times G)))\le \dim(I\times J),
\]
so if one of $I$ and $J$ is finite, then $\dim \pi(P)\le 1$ and we are done.

Suppose now that $\dim I =\dim J=1$. Since $P$ is $\CD$-definable and infinite, the projection of $P$ on one of the coordinates of $G^2$ has infinite image. Let us assume it is the projection on the first coordinate. Hence, since $\mr(P)=1$,
for every $\CD$-generic $a\in G$, the set $\{(w,z)\in G\times G: (a,w,z)\in P\}$ is
finite. Since $I\sub G$ is infinite every generic of $I$ is also $\CD$-generic in
$G$. But then, for such an $a\in I$ the set $\{(w,z)\in J\times G:(a,w,z)\in P\}$ is
finite. It follows that $\dim (P\cap (I\times
J\times G))=\dim I= 1$.\end{proof}

We proceed with a lemma towards the proof of Proposition \ref{attained3}.
\begin{lemma} \label{attained1} There exists a finite set $F\sub G$, with $F\sub \acl_{\CD}(\0)$,
 and a definable set $X\sub G$, with $\dim X\leq 1$,
 such that
for every  $b\in G\setminus X$ and for every $(x_1,x_2)\in G^2\setminus F^2$, if $b$
is attained near $(x_1,x_2)$, then $b$ is attained at $(x_1,x_2)$.
\end{lemma}
\begin{proof} Consider the $\CD$-definable set
$$T=\{(x_1,x_2,b)\in G^3:b\in S_{x_1}-S_{x_2}\}.$$ Since every generic fiber $S_x$ is finite,
$\mr(T)=2$.  Also, by fixing $x_1$ and letting $x_2$ vary, it is easy to see that the projection of $T$ on the last
coordinate is infinite and hence for every $\CD$-generic $b\in G$
 the set
$$T^b=\{(x_1,x_2)\in G^2: (x_1,x_2,b)\in T\} $$ has Morley rank $1$. Note that $(x_1,x_2)\in T^b$ if and only if $b$ is attained at $(x_1,x_2)$, and
$(x_1,x_2)\in \cl(T^b)$ if and only if $b$ is attained near $(x_1,x_2)$. We also
note, although we will not use this, that
$(x_1,x_2,b)\in \cl(T)$ if and only if $b$ is approximated near $(x_1,x_2)$.\\

\noindent\textbf{Claim 1.} {\em For $b\in G$ and $x_1,x_2\in G$, the following are
equivalent:
\begin{enumerate}
\item $b$ is attained near $(x_1,x_2)$ but not attained at $(x_1,x_2)$.

\item $(x_1,x_2)\in \fr(T^b) $ and $x_1,x_2\in \UB{S}$.

\item $(x_1,x_2)\in \fr(T^b).$
\end{enumerate}}

\begin{proof}[Proof of Claim 1] $(1)\Rightarrow (2)$. The fact  $(x_1,x_2)\in \fr(T^b) $ is immediate from the notes just above the claim. Since $b$ is attained near $(x_1,x_2)$, by curve selection in $\mathcal M$, we can find
definable curves $(x_1(t),y_1(t))\in S$ and $(x_2(t),y_2(t))\in S$ such that
$x_i(t)\to x_i$, for $i=1,2$, and $y_1(t)-y_2(t)=b$. Notice that $y_1(t)$ is bounded if and only if $y_2(t)$ is bounded, in which case,
since $S$ is closed, their limit points $y_1,y_2$ satisfy $(x_1,y_1), (x_2,y_2)\in
S$ and $y_2-y_1=b$, so $b$ is attained at $(x_1,x_2)$. Because we assumed that this
is not the case,  $y_1(t)$ and $y_2(t)$ are unbounded, hence $x_1,x_2$ are both in
$\UB{S}$.

The other implications are easy, thus ending the proof of Claim 1.\end{proof}

By Theorem \ref{main-frontier}, for each $b\in G$, $\fr(T^b)\sub \acl_{\CD}(b)$ (recall that $[S]$ was absorbed into the language). By compactness, we may therefore find a set $P\sub G^2\times
G$, $\CD$-definable over $\0$, such that for every $b\in G$ the set $P^b$ is finite
and contains $\fr(T^b)$. It follows that $\mr(P)=1$. Note however that we do not
claim that for every $b\in G$, we have $P^b=\fr(T^b)$. Thus, for example, we allow
at this stage the possibility that the set of $b$ for which $T^b$ is not closed is
$1$-dimensional.

Now, by Claim 1, if $b$ is attained near $(x_1,x_2)$ and not attained at $(x_1,x_2)$
then $(x_1,x_2)\in P^b$.

Assume first that the image of $P$ under the projection onto the $G^2$-coordinates,
call it $F_1$, is finite, and let $F\sub G$ be a finite set, $\CD$-definable over
$\acl_{\CD}(\0)$, such that $F_1\sub F^2$.  We may take $X=\es$ and complete the
proof of the lemma in this case. Assume then that $F_1$ is infinite.

Let $F_0\sub G^2$ be the set of all $p\in G^2$ such that  $P_p\sub G$ is
infinite. This is a finite set, $\CD$-definable over $\acl_{\CD}(\0)$, and because
we assumed that $F_1$ is infinite, the set $P^*:=(G^2\sm F_0)\times G$ still has
Morley rank $1$, and the projection map from  $P^*$ onto the $G^2$-coordinate is
finite-to-one.

 Set
$$X=\{b\in G: \fr(T^b)\setminus F_0\ne \es\}, \,\, \mbox{ a definable set in $\CM$
}.$$

\noindent\textbf{Claim 2.} $\dim(X)\leq 1$.

\begin{proof}[Proof of Claim 2] Assume towards contradiction that $\dim X=2$. For
every $b\in X$ there exists $(x_1,x_2)\in \fr(T^b)\sm F_0$.  By Claim 1 and our
choice of $P$,  $(x_1,x_2)\in (P^*)^b\cap (\UB{S}\times \UB{S})$, so since $\dim
X=2$, it follows that
$$\dim(P^*\cap (\UB{S}\times \UB{S}\times X))\geq 2.$$ This contradicts Claim
\ref{dimension1}. \end{proof}

  By  Claim 1, for every $b\in G$ and for every
$(x_1,x_2)\in G^2$, if $b$ is attained near $(x_1,x_2)$ and not at $(x_1,x_2)$, then
$(x_1,x_2)\in \fr(T^b)\sub P^b$. Now, either $(x_1,x_2)\in F_0$, or $b\in X$.
Thus, we may take any finite set $F\sub \acl_{\CD}(\0)$ with $F_0\sub F^2$ to complete the proof of  Lemma
\ref{attained1}.\end{proof}

%We are now ready to conclude the first part of Proposition \ref{attained3}

We now fix a definable $1$-dimensional $I\subseteq G$. Fix also a finite $F\sub G$ as in  Lemma \ref{attained1}, and   a definable
$1$-dimensional closed set $I_0\sub I$,
  such that $I_0\cap F=\emptyset$.

\begin{proof}[Proof of  Proposition \ref{attained3} ~(1)] Because
$S\cap (I_0\times G$) is a $1$-dimensional subset of $G\times G$, we may  shrink
$I_0$ further and assume that the set $S\cap (I_0\times G)$ is closed and
  bounded. Thus, the set $$B=\{b\in G:b\mbox{ is attained at }
   I_0\}=S(I_0)-S(I_0)$$ is a closed and bounded subset of $G$.
 By Lemma \ref{attained1} and the choice of $I_0$, there is a definable $X\sub
G$ with $\dim X\leq 1$ such that
 for every $b\in G\setminus X$, if $b$ is attained near
 $(x_1,x_2)\in I_0^2$, then $b$ is attained at $(x_1,x_2)$.
 Assume towards a contradiction that the set $A(S, I_0)\setminus B$ has dimension
$2$. By Lemma \ref{attained} (2), the set $L$ of all $b\in G\setminus B$ which are
attained near $I_0$ has dimension $2$, and therefore there is some $b\in L$ which is
not in $X$. By Lemma \ref{attained} (1), $b$ is attained near some $(x_1,x_2)\in
\cl(I_0)=I_0$, and since $b\notin X$, it is attained at $(x_1,x_2)$. Namely, $b\in
S(I_0)-S(I_0)=B$, a contradiction. \end{proof}

The rest of this subsection is devoted to the proof of Proposition
\ref{attained3}(2). Fix an open $V\sub G$ containing $0$. We may assume that $V$ is
bounded and symmetric, namely $-V=V$. Given $r>0$, let
  $P_r=\cl(B_r)\cap G$ and $S_r=\fr(B_r)\cap G$,
where $B_r$ is defined in the beginning of this section. Let $B$ be as in Proposition \ref{attained3}(1).
\\

\noindent{\bf Claim 1.} {\em There are $r_1>r_0>0$ sufficiently large such that
$B\sub P_{r_0}\sub P_{r_1}$ and  $S_{r_0}+V\sub P_{r_1}\setminus B.$}

  \begin{proof}[Proof of Claim 1] Since  $B+V$ is bounded,  there exists $r_0>0$ such that $B\sub P_{r_0}$
and $B+V$ does not intersect $S_{r_0}$. Since $V$ is symmetric, it follows that
$(S_{r_0}+V)\cap B=\emptyset$. Because $S_{r_0}+V$ is bounded  there exists
$r_1>r_0$ such that $S_{r_0}+V\sub P_{r_1}$. It follows that $S_{r_0}+V\sub
P_{r_1}\setminus B$.\end{proof}

Fix such $r_0,r_1$. For $\epsilon>0$ let as in Lemma \ref{attained} $$Y_{\epsilon}=\{y_2-y_1:\exists x_1,
x_2\in B(I_0;\epsilon)\,\, (x_1,y_1)\,,\, (x_2,y_2)\in S\}.$$

\noindent{\bf Claim 2}. {\em There exists $\epsilon_0>0$, such that no translate of
$V$ is contained in $(P_{r_1}\setminus B)\cap Y_{\epsilon_0}.$}

\begin{proof}[Proof of Claim 2]   The family of $Y_\epsilon$
decreases with $\epsilon$, and it is immediate from the definitions that
 $$A(S, I_0)=\bigcap_{\epsilon} \cl(Y_\epsilon).$$

We restrict our attention to the definably compact set $P_{r_1}\setminus \intr(B)$
and let
$$\bar Y_\epsilon^{r_1}=\cl(Y_\epsilon)\cap (P_{r_1}\setminus \intr(B)) \, \mbox{ and } \,
A_{r_1}(S,I_0)=A(S,I_0)\cap (P_{r_1}\setminus \intr(B)).$$ Thus, we have
$A_{r_1}(S,I_0)=\bigcap_{\epsilon>0}\bar Y_{\epsilon}^{r_1}.$ Each $\bar
Y_\epsilon^{r_1}$ is definably compact, and hence $A_{r_1}(S,I_0)$ is also definably
compact.

By the choice of $B$, Proposition \ref{attained3}(1) implies that
$\dim(A(S,I_0)\setminus B)\leq 1$ and hence, since the boundary of $B$ is at most
1-dimensional, also $\dim(A(S,I_0)\setminus \intr(B))\leq 1$. It follows that
$A_{r_1}(S,I_0)$ is a definably compact set which is at most $1$-dimensional. Using
that, it is not hard to see that for sufficiently small open $W\ni 0$ the set
$A_{r_1}(S,I_0)+W$ does not contain any translate of our open set $V$. Fix such a
set $W$.

Because $A_{r_1}(S,I_0)=\bigcap_{\epsilon}\bar Y_\epsilon^{r_1}$ it is not hard to
see that there exists $\epsilon_0>0$, such that $\bar Y_{\epsilon_0}^{r_1}\sub
A_{r_1}(S,I_0)+W$.
 It follows that the set $\bar Y_{\epsilon_0}^{r_1}$ does not contain any translate of $V$, thus proving Claim 2.\end{proof}

It is left to show that setting $\epsilon:=\epsilon_0$ for $\epsilon_0$ as in Claim 2, the requirements of Proposition \ref{attained3}~(2) are satisfied.
\\

\noindent{\bf Claim 3.} {\em There exists $r>0$ such that for all $x\in G\setminus
P_r$, $x+V\nsubseteq S(B(I_0,\epsilon_0)).$}

 \begin{proof}[Proof of Claim 3] Assume towards a contradiction that no such $r$ exists. Then we can find an
unbounded, definably connected curve $\Gamma\sub G$ such that $\Gamma+V\sub
S(B(I_0,\epsilon_0))$. It follows from the definition of $Y_{\epsilon_0}$ that
$(\Gamma+V)-(\Gamma+V)\sub Y_{\epsilon_0}$.

Fix any $\gamma_0\in \Gamma$ and let $\Gamma_0=\Gamma-\gamma_0$. The curve
$\Gamma_0$ is unbounded, definably connected, with $0\in \Gamma_0$ and in addition
$\Gamma_0+V\sub (\Gamma+V)-(\Gamma+V)\sub Y_{\epsilon_0}$. Note that $\Gamma_0\cap
S_{r_0}\neq \emptyset$, where $r_0$ as in Claim 1. Indeed, although $S_{r_0}=\fr(B_{r_0})\cap G$ need not be
definably connected, $\Gamma_0 \cap B_{r_0}\neq 0$  because $\Gamma_0$ is unbounded, definably connected and contains $0$. Fix $x_0\in \Gamma_0\cap S_{r_0}$. This intersection  point necessarily lies in $S_{r_0}$.

By our choice of $\Gamma_0$,  $x_0+V\sub \Gamma_0+V\sub Y_{\epsilon_0}$ and by our
choice of $r_0$ in Claim 1, $x_0+V\sub P_{r_1}\setminus B$. However, by Claim 2, no
translate of $V$ is contained in $Y_{\epsilon_0}\cap (P_{r_1}\setminus B)$,
contradiction.\end{proof}

Choose $r$ as in Claim 3.  Setting $B'=P_r$ and $\epsilon=\epsilon_0$ finishes the
proof of Proposition \ref{attained3}~(2).
\qed

\subsection{A lower bound  on the   dimension of the set of approximated points}

In this subsection,  assuming that $S_{\pol}$ is infinite, we modify the set $S$ from Proposition \ref{attained3} to a set $\hat S$ as in the next proposition, using an
idea from \cite[Section 4]{HaKo}.  The proof of Theorem \ref{inf-thm} in the next
subsection is by
 contradiction, and towards that we need this proposition.

%In this subsection, assuming that $S_{\pol}$ is infinite, we modify $S$, using an idea from \cite[Section 4]{HaKo}, to obtain a strongly minimal $\hat S$ as in the proposition below. Combined with Proposition \ref{attained3} this will quickly lead to a contradiction, allowing us to finally prove Theorem \ref{inf-thm} in the next subsection.

\begin{proposition}\label{lower}
Let $S\sub G^2$ be a $\CD$-definable strongly minimal, closed set which is not
$G$-affine, and assume that $\UB{S}$ is infinite. Then there is a  strongly minimal
closed set $\hat S\sub G^2$ which is not $G$-affine, definable in $\CD$ (over
additional parameters), and there exists  an infinite definable $\hat I\sub G$, such
that for every infinite definable set $T\sub \hat I$ and any  bounded ball $B$, the
set  $A(\hat S, T)\sm B$ is $2$-dimensional.
\end{proposition}

The rest of this subsection is devoted to the proof of Proposition \ref{lower}. We fix the sets $S$  and $\UB{S}$ as in its assumptions. Applying Proposition \ref{attained3} to $S$ and $S_{\pol}$ (in the role of $I$ there) we fix a definable $1$-dimensional $I_0\sub \UB{S}$ satisfying  Clauses (1) and (2) as in that proposition.

\begin{lemma}\label{attained2.9}  There is a definable smooth $1$-dimensional $I_1\sub I_0$ and
    \begin{enumerate}
      \item a  definably connected bounded open $U\sub G$,
       \item a definable continuous function $f:U\to G$ with $\Gamma_f\sub S$, and
           \item a definable family $\{\gamma_x : x\in I_1\}$ of curves $\gamma_x:(0,1)\to U$ with
            $\lim\limits_{t\to 0} \gamma_x(t)=x$, \linebreak $\lim\limits_{t\to 0}f(\gamma_x(t))=\infty$, and for  every $x_1, x_2\in I_1$,
    $$\lim\limits_{t\to 0} f(\gamma_{x_1}(t))-f(\gamma_{x_2}(t))=0.$$
    \end{enumerate}
\end{lemma}
\begin{proof} Using o-minimality and the fact that the projection of $S$ onto $G$ is
finite-to-one, we may partition $S$ and $I_0$ into finitely many cells and reach the
following situation. There is a definable, definably connected bounded open $U\sub
G$ and a definable $1$-dimensional smooth $I_1\sub I_0$, with $I_1$ on the boundary
of $U$ and $U\cup I_1$ a manifold with a boundary. We may assume that $\cl(U)\cap
\UB{S}=\cl(I_1)$. Furthermore, there is a definable, injective, continuous function
$f:U\to G$ whose graph is contained in $S$, such that for every $x_0\in I_1$ and
every curve $\gamma:(0,1)\to U$ tending to $x_0$ at $0$, the image of $\gamma$ under
$f$ is unbounded.

After applying a definable local diffeomorphism, we may assume that $I_1=(a,b)\times
\{0\}\sub R^2$ and $U=(a,b)\times (0,1)\sub R^2$. By shrinking $I_1$ if needed we
may assume that $f$ is defined on the box $[a,b]\times (0,1]$. For $\epsilon\leq 1$,
let
$$U_\epsilon=(a,b)\times (0,\epsilon)\sub U$$
and
 $$C_\epsilon= f([a,b]\times \{\epsilon\})\,, \,
\Gamma_{\epsilon,1}=f(\{a\}\times (0,\epsilon))\, , \,
\Gamma_{\epsilon,2}=f(\{b\}\times (0,\epsilon)).$$ When $\epsilon=1$, we denote
$C_1,\Gamma_{1,1}$ and $\Gamma_{1,2}$ by $C,\Gamma_1$ and $\Gamma_2$, respectively.
For every $\epsilon\leq 1$, the set $C_\epsilon$ is bounded and
$\Gamma_{\epsilon, i}$ are unbounded curves for $i=1,2$. Recall that $\partial
f(U_\epsilon)$ denotes the boundary of $f(U_\epsilon)$ (which is contained in $G$, since $G\subseteq R^n$ is closed).
Because $f:U\to G$ is continuous and injective it is in fact a homeomorphism, by
\cite{Johns}, hence
$$\partial f(U_\epsilon)= \Gamma_{\epsilon, 1} \cup \Gamma_{\epsilon, 2} \cup C_\epsilon$$
(we use here the fact that the limit of $|f(x)|$ as $x$ tends to any point in $I_1$
is $\infty$).

The next claim roughly says that for an infinitesimal $\epsilon$, the set
$f(U_\epsilon)$ is contained in  two infinitesimal tubes around  $\Gamma_1$ and $\Gamma_2$.\\

\noindent\textbf{Claim 1.} {\em For every $\epsilon_1>0$ there exists $\epsilon_2>0$
such that
$$f(U_{\epsilon_2})\sub \bigcup_{i=1}^2 \Gamma_i+B_{\epsilon_1}.$$}
\begin{proof}[Proof of Claim 1]  We fix $\epsilon_1>0$. Using Proposition \ref{attained3}(2), we can find
$\epsilon>0$ and a bounded neighborhood, $B'\sub G$, of $0$ such that for every
$y\in G\sm B'$, $y+B_{\epsilon_1}\not \subseteq f(U_{\epsilon})$. Next, choose
$0<\epsilon_2<\min\{\epsilon, \epsilon_1\}$, such that $f(U_{\epsilon_2})$ does not
intersect the bounded sets $B'$ and $C_{\epsilon}+B_{\epsilon_1}$. This can be done
since the limit of $|f(x)|$ is $\infty $ as $x$ tends in $U$ to any point in $I_1$.
We claim that this $\epsilon_2$ satisfies our requirements.

Indeed, given $x\in U_{\epsilon_2}$, we have $f(x)\notin B'$ and hence
$f(x)+B_{\epsilon_1}\not \subseteq f(U_{\epsilon})$. However, clearly $f(x)\in
f(U_\epsilon)$ (since $\epsilon_2<\epsilon$) and so,  because $f(x)+B_{\epsilon_1}$
is definably connected, we must have $(f(x)+B_{\epsilon_1})\cap \partial
f(U_\epsilon)\neq \emptyset$. Since $f(x)\not\in C_{\epsilon}+B_{\epsilon_1} $, we
have $(f(x)+B_{\epsilon_1})\cap C_\epsilon=\emptyset$, and therefore
$f(x)+B_{\epsilon_1}$ must intersect $\Gamma_{\epsilon,1}\cup \Gamma_{\epsilon,2}$,
and hence also $\Gamma_1\cup \Gamma_2$. It now follows that for some $i=1,2$,
$f(x)\in \Gamma_i+B_{\epsilon_1}$.\end{proof}

\noindent\textbf{Claim 2.} {\em There is a definable $1$-dimensional subset $I_2\sub
I_1$, and a definable family $\{\gamma_x:x\in I_2\}$ of
 curves $\gamma_x:(0,1)\to U$ with $\lim\limits_{t\to 0}\gamma_x(t)=x$, such that for every $x_1,x_2\in
I_2$, $$\lim\limits_{t\to 0} f( \gamma_{x_1}(t))-f(\gamma_{x_2}(t))=0.$$}

\begin{proof}[Proof of Claim 2] Consider the unbounded curves $\Gamma_1,\Gamma_2\sub \partial f(U)$,  and for each
 $i=1,2$  fix a definable
parametrization $\gamma_i(t):(0,1)\to G$ for $\Gamma_i$, such that $\lim_{t\to
0}|\gamma_i(t)|=\infty$.

Now fix a definable family $\{\gamma_x:x\in I_1\}$ of curves $\gamma_x:(0,1)\to U$
with $\lim\limits_{t\to 0}\gamma_x(t)=x$. By Claim 1, for each $x\in I_1$, the curve
$f(\gamma_x(t))$ approaches one of the $\Gamma_i$ as $t$ tends to $0$, and
therefore, after possibly re-parameterizing $\gamma_x$, we can find $\gamma_i$,
$i=1,2$, such that $\lim\limits_{t\to 0} f(\gamma_x(t))-\gamma_i(t)=0.$ The
re-parametrization can be done uniformly in $x$. We can now find an infinite
subinterval $I_2\sub I_1$ and $i\in \{1,2\}$ such that  if $x\in I_2$, then
$\lim\limits_{t\to 0}f(\gamma_x(t))-\gamma_i(t)=0$.
\end{proof}
Replacing $I_1$ by $I_2$ finishes the proof of Lemma \ref{attained2.9}.
\end{proof}

The rest of this subsection is devoted to the proof of Proposition \ref{lower}. We fix
$I_1\sub I_0, U, f, \{\gamma_x : x\in I_1\}$ as in Lemma  \ref{attained2.9}.
%Without loss of generality, $S$ is defined over $\es$.

\smallskip  Because $I_0$ is smooth on the boundary of $U$, we can  find an infinite sub-cell
$\hat I\sub I_0$ and $c\in G$ generic over $\0$
  such that $\cl(\hat I+c)$ is contained in $U$. We fix such $\hat I$ and  $c$. We say that two definable sets $X,Y$ have \emph{the same germ at $0$} if there is some open neighbourhood $W\ni 0$ such that $X\cap W=Y\cap W$. With this in hand, the key initial observation is the following.

\begin{claim} \label{inf5} For any infinite definable set $T\sub \hat I$, the set $V_c=f(T+c)-f(T+c)$ is
 a $2$-dimensional bounded set.\end{claim}
 \begin{proof} Since $f$ is continuous and $\cl(\hat I+c)\sub U$, it follows that $V_c$ is bounded.
  Assume now towards contradiction that $\dim V_c=1$. By shrinking $T$ further, if needed, we get from \cite[Lemma 2.7]{OtPe} that $f(T+c)$ is $G$-linear, that is, the sets $f(T+c)-g$ and $f(T+c)-h$ have the same germ at $0$ for all $h,g\in f(T+c)$.

%  ,  the one dimensional
%set $f(T+c)$ is a translate of a local subgroup $H_c$ of $G$ (by that  we mean
%that for all $x,y,z\in H_c$ sufficiently close to each other,  $x-y+z\in H_c$).

 By
shrinking $T$ if needed, we may assume that $c$ is still generic in $G$ over the parameters defining $T$. It follows that there is an open neighborhood $W\ni c$, such that for all $c'\in W$ the set $f(T+c')$ is $G$-linear. By definable choice, there is a definable function $g:W\to G$ such that $g(c')\in f(T+c')$ for all $c'\in W$. Denote $H(c'):=f(T+c')-g(c')$ and define an equivalence relation $E$ on $W$ by $E(c_1,c_2)$ if $H(c_1)$ and $H(c_2)$ have the same germ at $0$. Since $f(T+c')$ was $G$-linear, we easily get (see \cite{OtPe} for details) that $H(c')$ is a local subgroup.

We claim that there is a generic $E$-equivalence class that is infinite. Since $W$ is two dimensional, it will suffice to show that the class of germs at $0$ of the sets $H(c')$ is at most $1$-dimensional as $c'$ varies on $W$. Since the tangent space to $H(c)$ at $0$ is a subspace of the 2-dimensional tangent space to $G$ at $0$, our claim will follow from the fact that $H(c)$ and $H(c')$ have the same germ at $0$ if and only if they have the same tangent space at $0$. This latter fact is \cite[Claim 2.20]{PeStDefSimple} (note that the argument given there for definable subgroups goes through verbatim for germs of definable local subgroups).

If we now fix generic and independent
$x,y,z\in T$ sufficiently close to each other, then there is $w\in T$ and there are
infinitely many $E$-equivalent $c'$ such that
$$f(x+c')-f(y+c')+f(z+c')=f(w+c').$$
Since $\Gamma_f\subseteq S$, it follows readily from the above that $\stab^*(S)$ is infinite and therefore, by Lemma \ref{stabproperties}(4), that $S$ is $G$-affine, a contradiction.\end{proof}

Consider now the  $\CD$-definable set
\[
S'=\{(x,y_1-y_2) : (x+c,y_1) , (x,y_2) \in S\}.
\]
and the continuous function $\hat f:U\to G$
$$\hat f(x)=f(x+c)-f(x).$$
Clearly, $\mr(S')=1$ and $\Gamma(\hat f)\sub S'$. By Lemma \ref{rep-strongly}, there is a
$\CD$-definable strongly minimal set $\hat S\sub S'$ containing   $\Gamma(\hat f)$.
  Clearly, $\UB{\Gamma(\hat f)}\sub \UB{\hat S}$. Since $\fr(\hat S)$ is finite, we
  may assume that $\hat S$ is closed.

\begin{claim}  $\hat I\sub \UB{\hat S}$.\end{claim}
\begin{proof}  It suffices to prove $\hat I\sub\UB{\Gamma(\hat f)}$. Take $x\in \hat I$, and denote by $\gamma$ our fixed $\gamma_x:(0,1)\to U$. Then $\lim\limits_{t\to 0}\gamma(t)=x$. Also
$\hat f(\gamma(t))=f(\gamma(t)+c)-f(\gamma(t))$. Since
 $\lim\limits_{t\to 0}\gamma(t)+c=x+c$, it follows that
$\lim\limits_{t\to 0} f(\gamma(t)+c)=f(x+c)$, and because $\lim\limits_{t\to
0}|f(\gamma(t))|=\infty $,  also $\lim\limits_{t\to 0}|\hat
f(\gamma(t))|=\infty$, so $x$ is a pole of $\Gamma(\hat f)$.\end{proof}

Since $\UB{\hat S}\ne \es$, it follows that $\hat S$ not $G$-affine.\\

We can now proceed with the proof of Proposition \ref{lower}.  Let $T$ be any infinite
definable subset of $\hat I$, and $B$
 any open bounded ball.
We want to prove that $A(\hat S, T)\sm B$ has dimension $2$. \\

 \noindent\textbf{Claim 1.} {\em There is a definable unbounded
$1$-dimensional subgroup $H\sub G$, such that
 for every $x\in T$ and $h\in H$, there is a definable $\pi:(0,1)\to (0,1)$, with $\pi(0^+)=0^+$
  and
$$\lim\limits_{t\to 0} f(\gamma_x(\pi(t)))-f(\gamma_x(t))=h.$$ }
\begin{proof}[Proof of Claim 1] We first recall a theorem from \cite{PeSte}: given a definable curve
$\sigma:(0,1)\to G$ with $\lim\limits_{t\to 0}|\sigma(t)|=\infty$, the set of all
limit points of $\sigma(t)-\sigma(s)$, as $s$ and $t$ tend to $0$, forms an
$1$-dimensional torsion-free unbounded subgroup $H_\sigma\sub G$. In particular, for
each $h\in H_\sigma$ there is a definable function $\pi_h:(0,1)\to (0,1)$ with
$\pi_h(0^+)=0^+$ such that $\lim\limits_{t\to 0} \sigma(\pi_h(t))-\sigma(t)=h.$ It
follows from the definition of $H_{\sigma}$ that for every other definable curve
$\sigma':(0,1)\to G$, if $\lim\limits_{t\to 0} \sigma'(t)-\sigma(t)=0$, then
$H_{\sigma}=H_{\sigma'}$. We now apply this result to the unbounded curves
$f(\gamma_x(t))$, $x\in T$, and obtain the desired $H$.
\end{proof}

%For each $x\in T$, fix $H_x$ as in Claim 1 and observe that by its definition
%$H_x\sub A((S, T)\sub A(S, \hat I)$. We next claim that  the family $\{H_x : x\in T\}$ is finite.
%Indeed, if it is infinite then for every bounded $B\sub R^n$, the
%set $\bigcup_{x\in T} (H_x\sm B)$   is a $2$-dimensional set. Hence, $\dim A(S,
%\hat I)\sm B=2$. Since $\hat I\sub I_0$, $\dim A(S, I_0)\sm B)=2$, contradicting our choice of $I_0$ in Proposition \ref{attained3}(1).

%We fix an infinite definable set $T_1\sub T$ and a subgroup $H$ of $G$, such that
%for every $x\in T_1$,
%$H_x=H$. Denote $V_1=f(T_1+c)-f(T_1+c)$.\\

\noindent{\bf Claim 2.}  {\em For every $b\in V_1:=f(T+c)-f(T+c)$ and $h\in H$, we
have $b+h\in A(\hat S, T)$.}

  \begin{proof}[Proof of Claim 2] Let  $b=f(x_1+c)-f(x_2+c)\in V_1$, where $x_1, x_2\in T$, and let
   $\pi$ be as in Claim 1, for $x=x_2$ and $h$. Hence $h=\lim\limits_{t\to 0} f(\gamma_{x_2}(\pi(t)))-f(\gamma_{x_2}(t))$. We have:
\begin{align*}
&\hat f(\gamma_{x_1}(t))- \hat f(\gamma_{x_2}(\pi(t)))  =  \hat
f(\gamma_{x_1}(t))-\hat f(\gamma_{x_2}(\pi(t))) +
f(\gamma_{x_2}(t)) - f(\gamma_{x_2}(t)) \notag\\
&=[f(\gamma_{x_1}(t)+c)-f(\gamma_{x_1}(t)))]-
[f(\gamma_{x_2}(\pi(t))+c)-f(\gamma_{x_2}(\pi(t)))]+f(\gamma_{x_2}(t)) - f(\gamma_{x_2}(t))\notag \\
&=
[f(\gamma_{x_1}(t)+c)-f(\gamma_{x_2}(\pi(t))+c)]+[f(\gamma_{x_2}(t))-f(\gamma_{x_1}(t))]+
[f(\gamma_{x_2}(\pi(t)))-f(\gamma_{x_2}(t))].
\end{align*}

 As $t$ tends to $0$, for $i=1,2$, the curve $\gamma_{x_i}(\pi(t))+c$ still tends to
$x_i+c$, since $\pi(0^+)=0^+$, so its image under $f$ tends to $f(x_i+c)$. By Lemma \ref{attained2.9}(3),
$\lim\limits_{t\to 0} f(\gamma_{x_2}(t))-f(\gamma_{x_1}(t))=0$. Thus, by Claim 1, the above
expression tends to $f(x_1+c)-f(x_2+c)+h=b+h$, proving that $b+h$ can be
approximated near $T$.
\end{proof}

We can now conclude the proof of Proposition \ref{lower}, as follows. Because $V_1$ and
$B$ are bounded, we can find $r_0>0$ such that for every $h\in G\setminus B_{r_0}$,
the set $V_1+h\sub G\setminus  B$. In particular, for every $b\in V_1$, $b+(H\sm
B_{r_0})\sub G\sm B$.
 Moreover, since $H$ is unbounded, $H\sm B_{r_0}$ has dimension $1$. Hence, by Claim 2, the $2$-dimensional
  set $V_1+(H\sm B_{r_0})$ is contained in $A(\hat S, T)\sm B$, as needed.

\subsection{Proof of Theorem \ref{inf-thm}} Assume towards a contradiction that $\UB{S}$ is infinite.
%By Lemma \ref{L: S sm}
Since for any $S_1, S_2\sub G^2$, $\UB{(S_1\cup S_2)}=\UB{S_1}\cup\UB{S_2}$, and $\UB{S}=\UB{\cl(S)}$,
we may assume that $S$ is
strongly minimal and closed. Since $\UB{S}\ne \es$, we have that $S$ is not
$G$-affine. By Proposition \ref{lower}, there is a $\CD$-definable set $\hat S$ which is
closed, strongly minimal and not $G$-affine, and an infinite definable $\hat I\sub
G$, such that for every infinite set $T\sub \hat I$
and open bounded ball $B$, $A(\hat S, T)\sm B$ is $2$-dimensional. This contradicts  Proposition \ref{attained3}(1) for $\hat S$ and $\hat I$.\\

\begin{example}  One of the difficulties in the above proof was the
need to replace the initial set $S$ with a set $\hat S$, in order to reach a
situation where $\dim(A(\hat S, T)\sm B)=2$, for every infinite definable $T\sub \hat I\sub
\UB{\hat S}$ and any open bounded ball $B$. The following example shows that the initial $S$ can indeed have
infinitely many poles and yet $\dim A(S,I_0)=1$ for some (in fact, any bounded)
infinite $I_0\sub \UB{S}$. Consider the graph of the function $f:\mathbb R^2\to \mathbb
R^2$ defined by
$$f(x,y)=\begin{cases}
  (x,0) &  \text{ if } y=0 \\
   (xy,1/y) & \text{ if } y\neq 0
\end{cases},$$ with $G=\la \mathbb C,+\ra$.
The function $f$ is a bijection of $\mathbb C$ which is its own inverse. Its set of
poles is the  $x$-axis. For every $x\in \mathbb R$, as $(x,y)\to (x,0)$, $f(x,y)$
approaches the $y$-axis, with $|f(x,y)|\to \infty$. Thus, for any bounded $I_0\sub
\mathbb R\times \{0\}$, $A(S, I_0)= \text{$y$-axis}$. After  moving to $\hat S$ as
in the proof of Proposition \ref{lower}, we can see that $\dim (A(\hat S,T)\sm B)=2$, for
any infinite $T\sub \UB{\hat S}$ and bounded ball $B$.
\end{example}

\section{Topological corollaries}\label{sec-topcor}

\label{sec-topology} We establish here several topological properties of plane
curves, typically true for complex algebraic plane curves.  These properties are
used later on in our proof of the main theorem. Our first definition generalizes the notion of a function being open at a point.

\begin{defn} Let $S\sub G^2$ and $a=(a_1,a_2)\in S$. We say that {\em $S$ is open  at
$a$ over $a_1$} if for every open  neighborhood $U$ of $a$, $a_1$ is in the interior of $\pi_1(U\cap S)$. We say that {\em $S$ is open at $a$} if $S$ is open at $a$ over $a_1$ and $S^{\textup{op}}$ is open at $a$ over $a_2$.

Let $S\sub G^2$ and $a_1\in G$. We say that {\em $S$ is open over $a_1\in \pi_1(S)$} if for every $(a_1,a_2)\in S$, $S$ is open at $a$ over $a_1$.
\end{defn}

We note  that if  $B\ni a=(a_1,a_2)$ is an open box such that $a_1\notin
\intr(\pi_1(B\cap S))$, then the same remains true for all smaller open boxes.

\begin{lemma}\label{open-rel1}
    Assume that  $S\subseteq G^2$ is a plane curve. Then there are
    at most finitely many $a_1\in \pi_1(S)$ such that $S$ is not open over $a_1$. In
    particular, $S$ does not contain any $1$-dimensional components.

    If $S$ does not  contain any straight line, then there are at most finitely points $a\in S$ such that $S$ is not open at $a$.
\end{lemma}
\begin{proof} First note that if $S=S_1\cup S_2$ and $S$ is not open at $a$ over $a_1\in G$, then
either $S_1$ or $S_2$ is not open at $a$ over $a_1$. Thus we may assume
that $S$ is strongly minimal. Without loss of generality, $S$ is $\CD$-definable
over $\emptyset$.

Assume towards a contradiction that the set $N$ of all  $x$ in $\pi_1(S)$ over which
$S$ is not open is infinite. Pick $a_1$  generic in $N$ over $\0$. Because
$\mr(\pi_1(S))=1$, the point $a_1$ is  $\CD$-generic in $\pi_1(S)$ over $\emptyset$.

Fix $a=(a_1,a_2)\in S$ and  $B=B_1\times B_2\ni a$ such that $a_1\notin
\intr(\pi_1(S\cap B))$. Let $B_S=S\cap B$ and write $\bar B_S:=\cl(B_S)$. Note that
$a$ is $\CD$-generic in $S$ over $\emptyset$.

 By
Theorem \ref{inf-thm}, $S$ has  finitely many poles and since
$\dim(a_1/\emptyset)\geq 1$, the point $a_1$ is not a pole of $S$. By Corollary
\ref{injective}, there are at most finitely many points in $\pi_1(S)$ over which $S$
is  non-injective and each one of those is in $\acl_{\CD}([S])=\acl_{\CD}(\0)$. Thus $S$ is
injective over $a_1$. By Theorem \ref{main-frontier},  $\fr(S)\sub
\acl_{\CD}(\0)$ and hence we have $(\{a_1\}\times G)\cap \fr(S)=\emptyset$.

    %By our above assumption, the projection of $U_S$ on both coordinates is infinite.
    %By genericity of $a$ we may assume that $\pi_1(B)$, the projections of $B$ onto the first $G$-coordinate, contains only
    %$\CD$-generic elements.

Since $a_1\notin \intr(\pi_1(B_S))$, there exists
    a definable curve $\gamma:(0,1)\to B_1\setminus
    \pi_1(B_S)$ such that $\lim\limits_{t\to 0} \gamma(t)=a_1$. Notice that for $t$
    small enough $\gamma(t)$ must be $\CD$-generic in $G$, and therefore, because $\pi_1(S)$
    is co-finite in $G$,
    $\gamma(t)$ is $\CD$-generic in $\pi_1(S)$ over $\emptyset$. So, we may assume that
    the fiber $S_{\gamma(t)}$ has constant  size $n\geq 1$.
For each $t$,
     let $y_1(t),\dots, y_n(t)\in G$ be distinct such that $(\gamma(t),y_i(t))\in
     S$. Because $\gamma(t)\notin \pi_1(B_S)$, none of the $y_i(t)$ is in $B_2$.

Since $a_1\notin \UB{S}$, each of the curves $\gamma_i(t)$ is bounded, and hence has
a limit $y_i\in G\setminus B_2$. Since $(\{a_1\}\times G)\cap \fr(S)=\emptyset$, each
of the limit points $(a_1,y_i)$ is in $S$ and in addition $(a_1,a_2)\in S$, with
$a_2\neq y_i$ for all $i$. However, since $a_1$ is $\CD$-generic, we must have
$|S_{a_1}|=n$. This implies that for some $i\neq j$, we have $y_i=y_j$, so $S$ is
non-injective at $(a_1,y_i)$, contradiction.

Assume now that the intersection of $S$ with any straight line is finite. We apply the above to both $S$ and $S^{\textup{op}}$, and then by
removing from
$\pi_1(S)$ and $\pi_1(S^{\textup{op}})$ finitely many points, we remain, by our assumption on
$S$, with a co-finite subset $S'$ of $S$, such that $S$ is open at each point of $S'$. \end{proof}

%We are ready to deduce two useful topological corollaries:

\begin{corollary} \label{open-rel2} Assume that $S\sub G^2$ is strongly minimal and $a=(a_1,a_2)$ is a non-isolated point of $S$.

\begin{enumerate}

\item If $S$ is not  a straight line, then $S$ is open at $a$.

\item If there is $y\in G$ such that $S\sim G\times \{y\}$, then $y=a_2$
and there exists an open $U\ni a_1$ such that $U\times \{a_2\}\sub S$. In particular, $S$ is open at $a$ over $a_1$.

% \item If $S\sim \{x\}\times G$ then the dual result holds for $S^{\textup{op}}$.

\item  If $S^{\textup{op}}$ is injective at $(a_2,a_1)$ over $a_2$, then either $S\sim G\times \{a_2\}$  or
there exists an open $B=B_1\times B_2\ni a$ such that $S\cap B$ is
the graph of an open continuous map from $B_1$ into $B_2$.

% \item If $S^{\textup{op}}$ is injective at $(a_2,a_1)$ over $a_2$, then we are either under clause (2) or
% there exists an open $B=B_1\times B_2\ni a$ such that $S\cap B$ is the graph of an open continuous map from
%$B_1$ into $B_2$

\item If $S$ is not  a straight line and $a$ is $\CD$-generic in $S$, then there exists an open $U\ni a$ such that $S\cap U$
is the graph of a homeomorphism from $\pi_1(U)$ onto $\pi_2(U)$.

\end{enumerate}

\end{corollary}
\begin{proof}
(1) We assume that $S$ is not  a straight line and  show that $S$ is open at $a$. Assume towards a contradiction that $S$ is not  open at $a$ over $a_1$. In order to reach a contradiction it is sufficient, by Lemma
\ref{open-rel1}, to conclude that there are
infinitely many points in $\pi_1(S)$ over which $S$ is not open.

By Theorem \ref{main-frontier},  we can find an open box $B=B_1\times B_2$
containing $a$ such that $S\cap \cl(B)$ is closed, and $a_1\notin \intr(\pi_1(B\cap
S))$. Let $B_S=B\cap S$ and denote $\bar B_S=\cl(B_S)$. Repeating the argument with
a smaller box, we  see that we also have $a_1\notin \intr(\pi_1( \bar
B_S))=\intr(\pi_1(S\cap \cl(B)))$. Because $S\cap (\{a_1\}\times G)$ is finite, we
may also assume that $ S\cap (\{a_1\}\times \cl(B_2))=\{a\}$.

The set $\pi_1(\bar B_S)$ is closed in $G$ and since, by  Lemma \ref{open-rel1}, $S$ has no $1$-dimensional components and $\pi_1$ is finite-to-one, it is $2$-dimensional. The point $a_1$ belongs to the boundary of $\pi_1(\bar B_S)$, so by
o-minimality, there exists a definable curve $\gamma_1:(0,1)\to
\partial(\pi_1(\bar B_S))$,  with $a_1=\lim\limits_{t\to 0}\gamma_1(t)$. Since $\bar B_S=S\cap \cl(B)$, there exists a
definable curve $\gamma_2:(0,1)\to \cl(B_2)$ such that for every $t$,
$(\gamma_1(t),\gamma_2(t))\in S\cap \cl(B)$. Let $b=\lim\limits_{t\to
0}\gamma_2(t)\in \cl(B_2)$.

Since $S\cap \cl(B)$ is closed, it follows that $(a_1,b)\in S$, and therefore by our
choice of $B_2$, $b=a_2$. But then the curve $\gamma(t)=(\gamma_1(t),\gamma_2(t))$ tends
to $a$, so for small enough $t$, it must belong to the open set $B$, and its
projection is not in $\intr(\pi_1(B_S))$. Therefore  $S$ is not open over every
$\gamma_1(t)$ for $t$ small enough. This contradicts Lemma \ref{open-rel1} and ends
the proof of (1).

(2)  Assume that $S_1:=S\cap G \times \{y\}$ is infinite. Because $S$ is strongly minimal and $a$ is non-isolated
we must have $(a_1,a_2)\in S_1$, so $y=a_2$. The set $S_1$ is strongly minimal thus its projection on the first coordinate is co-finite and so $S_1$ (and therefore $S$) is locally
near $a$ the graph of a constant function. In particular, $S$ is open at $a$ over $a_1$.

%(3) This is identical to (2).

(3) Assume now that $S^{\textup{op}}$ is injective at $(a_2,a_1)$ over $a_2$ and that we are not under Clause (2).
Namely, the intersection of $S$ with any line $G\times \{y\}$ is finite.

By definition of injectivity, there exists an open box $B=B_1\times B_2$ such that
$B\cap S$ is the graph of a function, call it $f_S$, from a subset of $B_1$ into $B_2$, so the intersection of
each $\{x\}\times G$ with $S$ is finite. We may also assume that $B\cap S$
has no isolated point (by o-minimality, there are only finitely many).
By (1), we may shrink $B$ so
that $B_S$ is open over every point in $\pi_1(B_S)$ and $B_S^{-1}$ is open over every point in $\pi_1(B_S)$. It follows that the domain of $f_S$ is the whole of $B_1$
and in addition $f_S$ is continuous and open.

(4) By Corollary \ref{injective}, $S$ is injective at $a$ over $a_1$ and $S^{\textup{op}}$ is injective at $(a_2,a_1)$
over $a_2$. The result follows from (3), applied to $S$ and to $S^{\textup{op}}$. \end{proof}

Notice that even though, by o-minimality, the set of isolated points of any plane curve
$S$ is finite we do not know yet  that it is contained in
$\acl_\CD([S])$.

% However, by Corollary \ref{injective}, the non-injective points of
%$S$ are in $\acl_{\CD}([S])$ so we may conclude from the above:

%\begin{corollary}\label{open-rel4} Let $S\sub G^2$ be a strongly minimal
%whose projection on both coordinates is finite-to-one. If $a\in S$ is a non-isolated
%point which is $\CD$-generic in $S$ over $[S]$, then there exists an open box
%$B=B_1\times B_2\ni a$ such that $B\cap S$ is the graph of a homeomorphism $f:B_1\to
%B_2$.
%\end{corollary}

%\subsection{To include at an earlier point?}

%$\bullet$ Let  $X\sub G^2$  be  a $\CD$-definable plane curve and $p\in X$ is
%non-isolated point of $X$ such that $p\notin acl_{\CD}([X])$. Then there exists an
%open neighborhood $U=U_1\times U_2$ of $p$ such that $X\cap U$ is the graph of a
%continuous, open, injective function $f:U_1\to U_2$.

%$$\bullet$ Notation: For plane curves $X,Y\sub G^2$, let $$X\boxplus Y=\{(x,y_1\oplus
%$y_2):(x,y_1)\in X\,\, ,\,\, (x,y_2)\in Y\}.$$ We similarly define $\boxminus$.
%\section{More intersection theory}
%{\bf This section should be included at the end of the "non-injectity'' section.}

\section{The ring of Jacobian matrices}

\label{sec-jacobians}

\subsection{The ring $\mathfrak R$}
Our next goal is to show that if $f$ is a $\CD$-function, then its Jacobian matrix
vanishes at $0$ if and only if $f$ is not locally invertible at $0$. This will be
done in this and the next section. In the present section we prove that similarly to
a complex analytic function, the Jacobian matrix is non-zero if and only if it is an
invertible matrix.

Throughout this section we fix a definable local coordinate system for $G$ near
$0_G$, identifying $0_G$ with $0\in R^2$. From now on we identify $G$ locally with
an open subset of $R^2$. For a differentiable $\CD$-function $f$ in a neighborhood
of $0$, with $f(0)=0$, the Jacobian matrix at $x$, denoted by $J_xf$, is computed
with respect to this fixed coordinate system, and we denote by $|J_xf|$ its
determinant. We use $d_xf$ to denote the differential of $f$, viewed as a map from
the tangent space of $G$ at $x$, denoted by $T_x(G)$, to $T_{f(x)}(G)$. As we soon
observe, the collection of all matrices $J_0f$ is a subring of $M_2(R)$, and the main goal of
this section is to show that it is in fact a field (thus every nonzero matrix is
invertible).

We first observe the following statement.

\begin{lemma}\label{C:sgn}
    Let $f:U\to G$ be a non-constant $\CD$-function.
Then \begin{enumerate} \item The set of $a\in U$ at which $|J_af|=0$ is at most
$1$-dimensional. \item The set of $a\in U$ at which $J_af=0$ is finite.
\end{enumerate}
\end{lemma}
\proof (1) By strong minimality, for every open $V\sub U$, we have $\dim f(V)=2$,
for otherwise the pre-image of some point is infinite and co-infinite. By the
o-minimal version of Sard's Theorem (\cite[Theorem 2.7]{Wil99}), it follows that the set of singular points of $f$ is at
most $1$-dimensional.  For (2), note that if $J_af=0$ on a definably connected path
then $f$ must be constant there, which by strong minimality implies that $f$ is
constant on $U$.\qed

\begin{definition}\label{mfr}
Recall from Section \ref{sec-overview} that $\mathfrak F$ is the collection of all $\CD$-functions $f$ which are
$C^1$ in a neighborhood of $0$, with $f(0)=0$. We let $$\mathfrak R=\{J_0f\in
M_{2}(R): f\in \mathfrak F\}.$$
\end{definition}

It is important here to distinguish between the group operation in $G$ and the usual
ring operations in $M_2(R)$. Thus we reserve the additive notation $\pm$ for matrix
addition, and let $\oplus, \ominus$ denote the group operations in $G$.

\begin{lemma}\label{the ring}
 The set   $\mathfrak R$  is
 a subring with $1$ of $M_2(R)$ and for every $A\in \mathfrak R$ which is invertible, $A^{-1}\in \mathfrak R$.
\end{lemma}
\begin{proof}
  We first note that the collection of germs of functions in $\mathfrak F$ is
closed under $\oplus$ and functional composition. Indeed, if $S_f $ and $S_g$
represent $\CD$-functions $f$ and $g$ in $\mathfrak F$, then the plane curve
$S_f\circ S_g$ represents $f\circ g$ and the plane curve
$$S_f\boxplus  S_g=\{(x,y_1\oplus y_2):(x,y_1)\in S_f\,\, ,\,\ (x,y_2)\in S_g\}$$
represents $f\oplus g$.

Using the chain rule it is easy to verify that for $f,g\in \mathfrak F$,
$J_0(f\oplus g)=J_0f+J_0g,$ and $J_0(f\circ g)=J_0f\cdot J_0g.$ Since the germs in
$\mathfrak F$ are closed under $\oplus$ and functional composition, it follows
that $\mathfrak R$ is a ring. If $J_0 f$ is invertible, then $f$ is a locally
invertible function in which case it is clear that $f^{-1}$ is also in $\mathfrak
F$, and therefore $(J_0f)^{-1}\in \mathfrak R$.
\end{proof}

\noindent{\bf Note.} Given $\CD$-functions $f,g\in \mathfrak F$ it seems possible that every strongly minimal set representing $f \circ g$ (or every set representing $f\oplus g$) will have nodal singularity at $(0,0)$ and thus will not be locally at $(0,0)$ the graph of a function.

%It seems possible that even if $S_f$ and $S_g$ are, locally at $(0,0)$, the graphs of
%functions, the sets $S_f\circ S_g$ and $S_f\boxplus S_g$ need not be such, but they
%still represent $f\circ g$ and $f\oplus g$, respectively. It is in fact possible
%that no strongly minimal set representing $f\circ g$ will be locally at $(0,0)$ the
%graph of a function.

\subsection{Definability and  dimension of $\mathfrak R$}

Our aim is to show that $\mathfrak R$ is a definable field isomorphic to $R(\sqrt{-1})$. This is
achieved in several steps. We first show (Theorem \ref{P:ClassifyingmfR}) that
$\mfR$ is a definable ring of one of two kinds, and  then -- by eliminating one of
these possibilities -- we deduce the desired result.

 We are going to use extensively the following operation.
\begin{defn}  For a $\CD$-function $f$ which is
$C^1$ in a neighborhood of some $a\in G$, we let
 $$\tilde J_a f =J_0(f(x\op a)\om f(a)) \,\,; \,\, \tilde d_af=d_0(f(x\op a)\om f(a)).$$
\end{defn}
Note that $f(x\oplus a)\om f(a)$ is in $\mathfrak F$ and thus $\tilde J_a f\in
\mathfrak R$. We let $\ell_a(x)=x\op a$.

\begin{lemma}\label{basics}
\begin{enumerate}
\item \label{eq-jac}$\tilde d_a f=(d_{0}\ell_{f(a)})^{-1}\circ d_a f\circ
d_0\ell_a$.

\item For every $a\in dom(f)$, $J_af$ is invertible if and only if $\tilde J_af$ is
invertible, and $J_af=0 \Leftrightarrow \tilde J_af=0$.

\item  For any two differentiable $\CD$-functions $f,g:U\to G$ and $x_0\in U$,
$\tilde J_{x_0}(f\om g)=\tilde J_{x_0}f-\tilde J_{x_0} g.$

\end{enumerate}
\end{lemma}
\proof (1) is easy to verify and (2) follows, so we prove (3). Note that
\[\tilde J_{x_0}(f\om g)=J_0[(f\om g)(x_0\op x)\om (f\om g)(x_0)],\] which equals
\[
J_0[((f(x_0\op x)\om f(x_0))\om (g(x_0\op x)\om g(x_0))].
\]

As we noted in the proof of  Lemma \ref{the ring}, $J_0(h_1\om h_2)=J_0h_1-J_0h_2$,
therefore the above equals
\[
J_0(f(x_0\op x)\om f(x_0))- J_0(g(x_0\op x)\om g(x_0))=\tilde J_{x_0}(f)- \tilde
J_{x_0}(g).
\]
\qed

\begin{defn}
We say that a $\CD$-function $f:U\to G$ is {\em $G$-affine} if there exist non-empty
open sets $V\sub U$  and $W\ni 0$ such that for every $x_1,x_2\in V$ and $x\in W$,
$$f(x+ x_1)- f(x_1)=f(x+ x_2)- f(x_2).$$

%Given $S\sub G^2$ a strongly minimal set, we let
%$$\Stab^*(S)=\{g\in G^2:(g+ S)\Delta S \mbox{ is finite }\}.$$
\end{defn}

We note that a $\CD$-function  $f$ is $G$-affine if and only if $S_f$ is $G$-affine if and only if $\Stab^*(S_f)$ is infinite. Indeed, by Lemma \ref{stabproperties}(4) $S_f$ is $G$-affine if and only if $\stab^*(S_f)$ is infinite, and strong minimality implies that if $f$ is $G$-affine, then so is $S_f$. Furthermore, since $S_f$ is unique up to $\sim$-equivalence (as noted in the concluding paragraph of Section \ref{D-functions}) this does not depend on the choice of $S_f$.

\begin{remark}\label{L:vanishes} If $f$ is $G$-affine and $f(0)=0$, then $f$ is a partial group homomorphism, in
a neighborhood of $0$. It follows that for all $a$ in some open $V\ni 0$ we have $\tilde J_a f=J_0(f(x\oplus a)\ominus f(a))=J_0(f)=0$. Since $f$ is represented by some strongly minimal $S_f$, if it vanishes on some infinite set, $f$ vanishes on its domain.
\end{remark}

 As we already saw in Fact \ref{Hr-Pi}, since $\CD$ is not locally modular, there
exists at least one $\CD$-function which is not $G$-affine. \\

We are going to need the following lemma.
\begin{lemma}\label{invertible}
There are  invertible matrices in $\mathfrak R$ arbitrarily close to the $0$ matrix.
\end{lemma}
\begin{proof}
 We go via the following claim which is also used later in the text.
\begin{claim} \label{last claim} There exists $g\in \mathfrak F$ which is not $G$-affine, with
$J_0g=0$.
\end{claim}
\begin{proof} For $y\in G$ and $n\in \mathbb N$ we write
$ny:=\overbrace{y\op \cdots \op y}^{\mbox{ $n$-times}}.$ Fix $f:U\to G$ in
$\mathfrak F$ which is not $G$-affine, and for $n\in \mathbb N$, let
$g_n(x)=f(nx)-nf(x)$. It is easy to see that $g_n\in \mathfrak F$ and $J_0
g_n=nJ_0f-nJ_0f=0$. We want to show that for some $n\in \mathbb N$, the function
$g_n$ is not $G$-affine, so gives the desired $g$.

Notice that if $g_n$ is $G$-affine, then since $J_0g_n=0$, the function $g_n$ must
vanish on its domain (by Remark \ref{L:vanishes}). Assume towards contradiction that for every $n\in \mathbb N$,
the function $g_n$ vanishes on its domain, namely $f(nx)=nf(x)$ whenever $nx\in U$.
Pick a $\CD$-generic $x\in U$ sufficiently close to $0$ so that for all $n$,
$nx\in U$ and $nf(x)\in U$ (we can do it by saturation). For all $n$ we have
$$f(x+nx)=f((n+1)x)=(1+n)f(x)=f(x)+nf(x)=f(x)+f(nx).$$

Thus, since $x$ is generic, it is not torsion, and hence there are infinitely many $y\in G$ such that $f(x+y)=f(x)+f(y)$. Because
$f$ is a $\CD$-function it follows that for almost all $y$ with $x+y\in U$,
$f(x+y)=f(x)+f(y)$. Since $x$ is $\CD$-generic, the function $f$ must be
$G$-affine, a contradiction.\end{proof}

We return to the proof of Lemma \ref{invertible}. Take the function $g:V\to G$ from Claim \ref{last claim}. By Lemma \ref{C:sgn}, for every $x\in V$ generic, $J_x g$ and hence $\tilde
J_x g$ is invertible. Because $g$ is smooth and $J_0g=0$ there are invertible
matrices of the form $\tilde J_xg\in \mathfrak R$ arbitrarily close to the $0$
matrix.\end{proof}

\begin{defn}\label{def-Wrealized} Given a set $W\sub \mathfrak R$ and a family $\CF=\{f_t:t\in T_0\}$ of
$\CD$-functions, we say that {\em $W$ is realized by $\CF$} if
$$W=\{J_0 f_t:t\in T_0\}.$$

%We that $\{J_0 f_t:t\in T_0\}$ is {\em a continuous family} if the map $(t,x)\mapsto
%J_0f_t$ is continuous at $(t,0)$ for every $t\in T_0$.
\end{defn}
\begin{proposition}\label{definability}
The ring $\mathfrak R$ is a definable subring of $M_2(R)$ which is also an $R$-vector subspace.
\end{proposition}
\proof We first show that $\mathfrak R$ is a $\bigvee$-definable
subring of $M_2(R)$,  namely that $\mathfrak R$ is a bounded union of definable subsets of $M_2(R)$. Let $M\in \mathfrak R$. By definition, there exists some $\CD$-function $f\in \mathfrak F$ such that $J_0f=M$. Let $S_f$ represent $f$. Let $\phi(x,a)$ $\CD$-define $S_f$  such that $\phi(x,y)$ is a family of plane curves all passing through $(0,0)$. By Proposition \ref{rep-Dfunction}, there is a $\0$-definable family  $\mathcal F$ of $\CD$-functions in $\mathfrak{F}$ such that the germ of $f$ at $0$ is represented in $\CF$.  Since $J_0f$ only depends on the germ of $f$ at $0$, we get that $M$ is realised as the Jacobian at $0$ of some $\CD$-function in $\CF$. Let $T_\CF$ be the set of all jacobians of  $\CD$-functions in $\CF$, where $\CF$ is a $\0$-definable family of $\CD$-functions in $\mathfrak{F}$.  We have thus seen that $\mathfrak R$ can be covered by all the sets $T_\CF$. There is a bounded number of such sets, where the bound is given by the cardinality of the language of $\CD$.

% By
%Proposition \ref{rep-Dfunction} and the remarks following its proof, any family $\{S_t:t\in T\}$ of plane curves all containing $(0,0)$ $\CD$-definable over $\emptyset$
%gives rise to a $\0$-definable family of $\CD$-functions in $\mathfrak F$,
%$\CF_0=\{f_t:t\in T_0\}$ containing precisely all functions in $\mathfrak F$ whose
%germ at $0$ is represented by $S_t$ for some $t\in T$. The family $\CF_0$ realizes the definable set $J=\{J_0f_t:t\in T_0\}\sub \mathfrak R$ in the sense of Definition \ref{def-Wrealized}.
%
%As pointed out in the proof of Lemma \ref{the ring}, the families $$\{f_t\oplus
%f_s:t,s\in T_0\}\,\, ,\,\, \{f_s\circ f_t:s,t\in T_0\}\,\, ,\,\, \{f_t^{-1}:t\in T_0
%\mbox{ and $f_t$ invertible near $0$}\}$$ are also families of $\CD$-functions.
%Hence, the sets of matrices $J+J$, $J\cdot J$, and $J^{-1}=\{A^{-1}: A\in J \mbox{
%invertible}\}$ can be realized by definable families of $\CD$-functions. Thus
%$\mathfrak R$ can be first identified as a bounded $\bigvee$-definable subring of
%$M_2(R)$, where the set of disjuncts is bounded in size by the cardinality of the
%language of $\CD$.
It follows that there is a definable open neighborhood $U\sub
M_2(R)$ of the zero matrix, such that $U\cap \mathfrak R$ is definable (for more on
$\bigvee$-definable groups and rings see \cite{PeSt4}). More precisely, there exists a
$\emptyset$-definable family of $\CD$-functions which realizes $U\cap \mathfrak R$.

% $\CD$-definable
%family of plane curves $\{S_t:t\in T\}$ and a definable $T_0\sub T$, such that
%$U\cap \mathfrak R$ is realized by $\{f_t:t\in T_0\}$.

We now proceed to show that $\mathfrak R$ is actually a \emph{definable} subset of
$M_2(R)$.  Let $U\ni 0$ be a neighborhood of $0$ in $M_2(R)$ such that $U\cap
\mathfrak R$ is definable as above.
 We  claim that
$$\mathfrak R=\{AB^{-1}:A, B\in U\cap \mathfrak R \, \, , B \mbox{ is invertible}\}.$$
Indeed, for every $C\in \mathfrak R$ we can find, by Lemma \ref{invertible}, an
invertible matrix $B\in U\cap \mathfrak R$, sufficiently close to $0$, such that
$CB\in U\cap \mathfrak{R}$. It follows that $\mathfrak R$ is definable.

Finally, the subring of scalar matrices in $\mathfrak R$ is, in particular, a subgroup of $(R,+)$, and it is non-trivial since it contains $1$. By o-minimality, the only non-trivial definable subgroup of $(R,+)$ is $(R,+)$ itself. So $\mathfrak R$ contains all diagonal matrices, and is therefore an $R$-vector subspace of $M_2(R)$.
\qed

\begin{proposition}\label{thm-dim}
  Let $U\sub G$ be an open neighborhood of $0$, and  assume that $f:U\to G$ is a continuously differentiable
$\CD$-function which is not $G$-affine.  Then the set $\tilde J(U)=\{\tilde J_a f\in
M_2(R): a\in U\}$ has dimension $2$.  In particular, $\dim \mathfrak R\geq 2$.
\end{proposition}
\proof Since $\dim U=2$, we have $\dim \tilde J(U)\leq 2$. Assume towards a
contradiction, that $\dim \tilde J(U)\leq 1$.
\\

\noindent {\bf Claim.} {\em There exists $g_0\in G$, $g_0\notin \dcl(\emptyset)$,
and infinitely many $a\in G$ such that $\tilde J_af=\tilde J_{a}(f(x\op g_0))$.}

\begin{proof}[Proof of Claim]
  For every matrix $A\in \tilde J(U)$ let  $C_A:=\{x\in U:\tilde J_xf=A\}$. By
our assumptions, there exists $A\in \tilde J(U)$ such that $\dim C_A\geq 1$, and by
possibly shrinking $U$, we may assume that $C_A$ is definably connected. Consider
$B_A=C_A\ominus C_A\sub G$. There are two cases to consider:
\\

\noindent{\bf Case 1.} There exists  $A\in \tilde J(U)$ such that $\dim B_A=1$.
\\

We may apply \cite[Lemma 2.7]{OtPe} and conclude that the set $C_A$ consists of a
subset of a coset of a $\bigvee$-definable one dimensional subgroup $\CH$ of $G$. It
follows that for $g_0\in \CH$ sufficiently close to $0$,  there are infinitely many $a\in
C_A$ such that $a\op g_0\in C_A$, and thus $\tilde J_af=\tilde J_{a}f(x\op
g_0)=\tilde J_{a+g_0}f$.
\\

\noindent{\bf Case 2.} For all $A\in \tilde J(U)$, $\dim B_A=2$, so $B_A$ contains an
open subset of $G$.
\\

Given $A$ generic in $\tilde J(U)$  we may find an open set $W\sub G$ in $B_A$ such
that $A$ is still generic in $\tilde J(U)$ over the parameters defining $W$. Thus
there are infinitely many $A\in \tilde J(U)$ for which $W\sub B_A$. Pick $g_0$
generic in $W$ and then for each $A$ such that $g_0\in B_A$ there are $a,b\in C_A$
such that $a\om b=g_0$, so $a=b\op g_0$. By definition of $C_A$, we know that for
every such pair $(a,b)$ we have $\tilde J_a f =\tilde J_b f$, so
 $\tilde J_{b\op g_0}f=\tilde J_b f$. We get:
\begin{align*}
\tilde J_b f(x\oplus g_0)=J_0(f(x\op b\op g_0)\om f(b\op g_0)) = \tilde J_{b\op g_0}
f =\tilde J_b f.
\end{align*} Since there are infinitely many such
pairs $b, b\op g_0$, as $A$ varies, we are done.
\end{proof}

To conclude the proof, fix $g_0$ as in the claim and infinitely many $a$
such that $\tilde J_af=\tilde J_af(x\oplus g_0)$. By Lemma \ref{basics} (3), for
each such $a$,  $\tilde J_a(f(x\op g_0)\om f(x))=0$. But then, by Lemma \ref{basics} (2), for the $\CD$-function $k(x)=f(x\op g_0)\om f(x)$  there are infinitely many $a$, such that
$J_a k=0$, so $k(x)$ is constant on
its domain, say of value $d$. By strong minimality of $\CD$,  $(g_0,d)$ is in $\Stab^*(S_f)$.
Since $g_0$ is not in $\dcl(\0)$, it is not a torsion-element so  $\Stab^*(S_f)$ is
infinite and therefore $f$ is $G$-affine, contradiction.\qed

\subsection{The structure of $\mathfrak R$} The main result of this section is the following theorem.
\begin{theorem}\label{P:ClassifyingmfR} There exists a fixed invertible matrix $M\in
M_2(R)$ such that one of the following two holds:
    \begin{enumerate}
        \item $$\mathfrak R=\{M^{-1}\begin{pmatrix}
        a & -b\\
        b & a
        \end{pmatrix}M: a,b\in R\}.$$ In particular, $\mathfrak R$ is a field which is definably
isomorphic to $R(\sqrt{-1})$. Or,

        \item $$\mathfrak R=\{M^{-1}\begin{pmatrix}
        a & 0\\
        b & a
        \end{pmatrix}M: a,b\in R\}.$$
\end{enumerate}
\end{theorem}

 We  need some preliminaries.
\begin{lemma}\label{L:ConstSgn}
    Let $U\subseteq G$ be a definably connected open neighbourhood of $0$. Let $f:U\to G$ be
     a non-constant $\CD$-function. Then $|J_xf|$ has constant sign at all $x\in U$ where $f$ is differentiable
      and $J_xf$ is invertible.
\end{lemma}
\begin{proof}
     By Corollary \ref{injective}, we may
      assume -- possibly removing finitely many points from $U$ -- that $f$ is locally injective.
      The result now follows from  \cite[Theorem 3.2]{PetStaTopInv}.
\end{proof}

  Now, for $f\in \mathfrak F$ non-constant  we denote by $\sigma(f)$ the sign of
$|J_x(f)|$ for all $x$ sufficiently close to $0$ at which $J_xf$ is invertible.

%Let $\{f_t:t\in T\}$ be a definable family of functions into $G$, each defined and
%differentiable on an open neighborhood $U_t\sub G$ of $0$. We say that {\em the
%family is continuous} if the function $F(t,x)=f_t(x)$ is continuous at every
%$(t,0)$. If $T$ is definably connected then by \cite[Theorem 3.2]{PetStaTopInv},
%$\sigma(f_t)$ is constant as $t$ varies in $T$.

\begin{proposition}\label{C:no-neg}
    Every invertible $A\in \mathfrak R$ has positive determinant.
\end{proposition}
\begin{proof}
Fix $A_0\in \mathfrak R$ generic over $\0$, and $W\sub \mathfrak R$ a definable open
neighborhood of $A_0$. Fix also a definable family of $\CD$-functions, $\{f_t:t\in
T\}$ realizing $W$, provided by Proposition \ref{definability}. Let $a_0$ be generic in $T$, such that $J_0f_{a_0}=A_0$. We may
assume that $T$ is a cell in some $R^k$, and by definable choice in o-minimal
structures further assume that the map $t\mapsto J_0f_t$ is a homeomorphism of $T$
and $W$. By Proposition \ref{thm-dim}, $\dim T=\dim W=\dim \mathfrak R\geq 2$.

For every $t\in T$, let $U_t\sub G$ be the domain of $f_t$ (containing $0$). We can
find a definably connected neighborhood $U_0\ni 0$ and a definably connected
neighborhood $T\supseteq T_0\ni a_0$, such that for every $t\in T_0$, $U_0\sub U_t$. The definition of the
sets $U_0$ and $T_0$ may use additional parameters but we may choose them so that
$a_0$ is still generic in $T_0$ over those parameters. Let $W_0=\{J_0 f_t:t\in
T_0\}$ be the corresponding neighborhood of $A_0$ in $\mathfrak R$.

Consider now the set of matrices $\hat W_0=W_0-A_0\sub \mathfrak R$. It is an open
neighborhood of $0$ in $\mathfrak R$, which is realized by the family $\{f_t\ominus
f_{a_0}:t\in T_0\}.$

Our goal is to show that every invertible matrix in $\hat W_0$ has positive
determinant. Let us first see that they all have the same determinant sign. Note
that for all $t\in T_0\setminus \{ a_0\}$, the function $f_t\ominus f_{a_0}$ is
non-constant on $U_0$, thus by Theorem \ref{C:continuous} it is an open map. We now
show, using our above notation, that $\sigma(f_t\ominus f_{a_0})$ is constant as $t$
varies in a punctured neighborhood of $a_0$.

Fix $x_0\in U_0$ which is generic over $a_0$.  Since $(a_0,x_0)$ is generic in
$T_0\times U_0$, there exist an open $T_0'\ni a_0$ inside $T_0$ and an open $\hat
U_0\ni x_0$ inside $U_0$ such that the map $F(t,x)=f_t(x)\ominus f_{a_0}(x)$ is
continuous on $T_0'\times U_0$. Because $\dim T_0=\dim W\geq 2$, the set $\hat
T_0=T_0'\setminus \{a_0\}$ is still definably connected, and for each $t \in \hat
T_0$, the function $f_t\ominus f_{a_0}$ is open on $\hat U_0$. Given, $t_1\neq
t_2\in \hat T_0$, there exists a definable path $p:[0,1]\to \hat T_0$ connecting
$t_1$ and $t_2$, and by possibly shrinking $\hat U_0$, the induced map $(s,x)\mapsto
F(p(s),x)$ is a definable proper homotopy (see \cite[Section  3.5.1]{PetStaTopInv}) of $f_{t_1}\om f_{a_0}$ and $f_{t_2}\om f_{a_0}$,
hence by \cite[Theorem 3.19]{PetStaTopInv}, for every $x$ generic in $\hat U_0$,
$|J_{x}(f_{t_1}\ominus f_{a_0})|$ and $|J_{x}(f_{t_2}\ominus f_{a_0})|$ have the
same sign. It follows that
$$\sigma(f_{t_1}\ominus f_{a_0})=\sigma(f_{t_2}\ominus f_{a_0}).$$ Thus, every
invertible matrix in $\hat W_0$ has the same determinant sign.

Next, note that for every invertible $A\in \hat W_0$ sufficiently close to $0$, the
matrix $A^2$ is also in $\hat W_0$ and clearly has positive determinant. Thus all
invertible matrices in $\hat W_0$  have positive determinant.

 Finally, as we saw in the proof of Proposition \ref{definability}, $\mathfrak
R=\{AB^{-1}:A,B\in \hat W_0, B \mbox{ invertible}\}$, and hence all invertible
matrices in $\mathfrak R$ have positive determinant.\end{proof}

\begin{proof}[Proof of Theorem   \ref{P:ClassifyingmfR}] Assume first that every
non-zero $A\in \mathfrak R$ is invertible, namely that $\mathfrak R$ is a definable
division ring. It follows from \cite[Theorem 4.1]{PeSte} that $\mathfrak R$ is
definably isomorphic to either $R$ or $R(\sqrt{-1})$ or the ring of quaternions over
$R$. Because $\dim \mathfrak R\geq 2$, we are left with the last two possibilities.
The ring of quaternions, $H(R)$, is not definably isomorphic to a definable subring of $M_2(R)$. Indeed, if $A\subseteq M_2(R)$ is such a ring then, as we have seen in the proof of Proposition \ref{definability}, $A$ contains all scalar matrices. Because $H(R)$ has o-minimal dimension $4$ the same must be true of $A$, so $A$ is a $4$-dimensional $R$-vector subspace of $M_2(R)$. So $A=M_2(R)$, but the latter is not isomorphic to $H(R)$ since it is not a division ring.

%Indeed, an easy calculation shows that the only matrices in $M_2(R)$ satisfying $M^2=-1$ are of the form $\begin{pmatrix}
%0 & a \\
%-a^{-1} & 0
%\end{pmatrix}$ and the product of three such matrices is never $-1$, so there are no $M_1,M_2,M_3\in M_2(R)$ satisfying the equations \[
%M_1^2=M_2^2=M_3^2=M_1M_2M_3=-1.
%\]
So $\mathfrak R$ is necessarily isomorphic to $R(\sqrt{-1})$. Since $R(\sqrt{-1})\cong R\oplus iR$ and $\mathfrak R$ is a subring of $M_2(R)$, we immediately see that $\mathfrak R$ is generated, as a vector space over $R$ by the diagonal matrices and some matrix $M(i)$ such that $M(i)^2=-1$. It follows that the eigenvalues of $M(i)$ are $\pm i$, so $M(i)$ is diagonalizable and conjugate to $\begin{pmatrix}
    0 & 1 \\
    -1 & 0 \end{pmatrix}$, say via some matrix $M$.
It is now immediate that $\mathfrak R$ is of
the form (1) with respect to this matrix $M$.

We thus assume that there exists  at least one matrix $A$ that is not invertible,
of rank $1$. We want to show that there exists an invertible $M\in M_2(R)$ such that
$\mathfrak R$ has form as in (2).

%We fix
% some open $W\subseteq \mathfrak R$, a
% neighbourhood of $0$ such that  $-W=W$.

% By Lemma \ref{invertible}, we may replace $A$ with $AC$ for some invertible $C\in \mathfrak R$ close enough to
%$0$, and thus
% we may assume that $A\in W$.
We conjugate $\mathfrak R$ by some fixed matrix so that $A$, written in columns, has
the form $(w,0)$ for some $w\in R^2$. We now show that every matrix in $\mathfrak R$
is of the form $\begin{pmatrix}
        a & 0\\
        b & a
        \end{pmatrix}$ for some $a,b\in R$.
Consider  the set
 $$H=\{(u,0)\in \mfR: u\in R^2\}.$$ It is a
definable $R$-vector subspace of $\mathfrak R$ that is also closed under ring multiplication. As an $R$-vector space it has positive dimension over $R$ (since $H$ is a non-trivial subring of $\mathfrak R$) and $\dim_R H\leq 2$. \\

\noindent\textbf{Claim 1.} {\em  $\dim_R(H)=1$ }

\begin{proof}[Proof of Claim 1]
Write the matrices in $H$ in the form $B=\begin{pmatrix}
      a & 0\\
        b & 0
        \end{pmatrix}$, and note that for $C=\begin{pmatrix}
        c & d\\
        e & f
        \end{pmatrix}$,

\begin{equation}\label{det} |B+C|=|C|+(af-bd).\end{equation}

Assume towards a contradiction that $\dim H=2$, and then $H$ consists of {\em all}
matrices of the form $\begin{pmatrix}
      a & 0\\
        b & 0
        \end{pmatrix}$. We may now take $C=\begin{pmatrix}
        c & d\\
        e & f
        \end{pmatrix}\in \mathfrak R$ invertible, sufficiently close to $0$, and since $d,f$ cannot be both $0$, it easy
to see that by choosing $a,b$ appropriately, we may obtain a matrix $B+C\in
\mathfrak R$ whose determinant is negative, a contradiction.
\end{proof}

  Thus, $H$ is a $1$-dimensional $R$-vector space.\\

\noindent\textbf{Claim 2.} {\em The matrices in $H$ are not of the form
    $B=\begin{pmatrix}
    a & 0\\
    \alpha a & 0\\
    \end{pmatrix}$, for some fixed $\alpha \in \mathbb R$. }

\begin{proof}[Proof of Claim 2] Towards a contradiction, assume that there is such an $\alpha$. Take any invertible $C=\begin{pmatrix}
    c & d \\
    e & f
    \end{pmatrix}\in \mfR$. Then for every
 $B=\begin{pmatrix}
    a & 0 \\
    \alpha a & 0
    \end{pmatrix}\in H$, we have $|B+C|=|C|+(af-\alpha a d)$. By choosing $a$ appropriately, we
obtain $|B+C|<0$ (contradicting Proposition \ref{C:no-neg}), unless $f= \alpha d$. Hence $f=\alpha d$ and $C=\begin{pmatrix}
    c & d \\
    e & \alpha d
    \end{pmatrix}\in \mfR$. We  have
$$ \begin{pmatrix}
    c & d \\
    e & \alpha d
    \end{pmatrix}\cdot \begin{pmatrix}
    a & 0\\
    \alpha a & 0\\
    \end{pmatrix}=\begin{pmatrix}
    a(c+\alpha d)& 0 \\
    a(e+\alpha^2d) & 0
    \end{pmatrix}. $$
But the left hand side of the above equation is in $\mathfrak R$, and $H$ is the collection of all matrices in $\mathfrak R$ of the form $(u,0)$ for $u\in R^2$. Hence, $\begin{pmatrix}
a(c+\alpha d)& 0 \\
a(e+\alpha^2d) & 0
\end{pmatrix}\in H$.  By assumption, $\alpha(c+\alpha d)=e+\alpha^2 d$,
implying that $e=\alpha c$. However, this would make $C$ non-invertible, a
contradiction.
\end{proof}

We are thus left with the case that matrices in $H$ are of the form $\begin{pmatrix}
    0 & 0\\
    a & 0\\
    \end{pmatrix}$.  If we now take an arbitrary  $\begin{pmatrix}
    c & d \\
    e & f
    \end{pmatrix}\in \mathfrak R$ and multiply it on the right by a non-zero element of $H$,
we obtain another element of $H$, forcing $d$ to be $0$. Thus all matrices in
$\mathfrak R$  are lower triangular. Because $H$ contains all matrices of the form
$\begin{pmatrix}
    0 & 0\\
    a & 0\\
    \end{pmatrix},$ every matrix in $\mathfrak R$ can be written as the sum of a
diagonal matrix in $\mathfrak R$ and a matrix in $H$.\\

\noindent\textbf{Claim 3.} The diagonal matrices in $\mathfrak R$ are precisely the scalar matrices.
\begin{proof}[Proof of Claim 3]
    The set of diagonal matrices in $\mathfrak R$ is a definable additive subgroups of the ring of all diagonal matrices
    \[
    \Delta:=\left \{ \begin{pmatrix}
    a & 0 \\
    0 & b
    \end{pmatrix}: a,b\in R\right\}.
    \]
    We have already seen that all scalar matrices are in $\mathfrak R$, so if $\mathfrak R$ contains any matrix in $\Delta$ that is no--scalar, it contains all of $\Delta$, contradicting Proposition \ref{C:no-neg}.
%We already know (Lemma \ref{the ring}) that all scalar matrices are in $\mathfrak R$. We know from the form of $H$ that there are no diagonal matrices of the form $\begin{pmatrix}
%a & 0 \\
%b & 0
%\end{pmatrix}$ for $a\neq 0$. Now consider a matrix $\begin{pmatrix}
%0 & 0 \\
%b & a
%\end{pmatrix}\in \mathfrak R$.
%Since $\begin{pmatrix}
%0 & 0 \\
%b & 0
%\end{pmatrix}\in H$, we get that $\begin{pmatrix}
%0 & 0 \\
%0 & a
%\end{pmatrix}\in \mathfrak R$, and since the scalar matrices are in $\mathfrak R$, all matrices of this form are in $\mathfrak R$. So, if $\begin{pmatrix}
%a & 0 \\
%b & c
%\end{pmatrix}\in \mathfrak R$, then also $\begin{pmatrix}
%a & 0 \\
%b & c
%\end{pmatrix}-\begin{pmatrix}
%0 & 0 \\
%0 & c
%\end{pmatrix}=\begin{pmatrix}
%a & 0 \\
%b & 0
%\end{pmatrix}\in \mathfrak R$, implying that $a=0$, but since all scalar matrices are in $\mathfrak R$ this is impossible.
%
%So all diagonal matrices in $\mathfrak R$ are either in $H$ or  of the form
%$\begin{pmatrix}
%a & 0 \\
%0 & c
%\end{pmatrix}$, with $a,c\neq 0$. So all diagonal matrices in $\mathfrak R$ are invertible, and therefore form a field that is a subring of $\mathfrak R$. But since $\mathfrak R$ does not contain any matrices of the form $\begin{pmatrix}
%0 & a \\
%-a^{-1} & 0 \\
%\end{pmatrix}$, it cannot be isomorphic to $R(\sqrt{-1})$ (see the proof fo Proposition \ref{C:no-neg}), so it must be definably isomorphic to $R$. Since it contains all the scalar matrices as a subfield, it must equal the field of scalar matrices.
\end{proof}

It follows that the matrices in $\mathfrak R$ are of the form $\begin{pmatrix}
    a & 0 \\
    b & a
    \end{pmatrix},$  as required. This ends the proof of Theorem \ref{P:ClassifyingmfR}.\end{proof}

\subsection{From ring to field}
\label{sec:from ring to field}
\begin{definition}
    We say that $\mathfrak R$ is of \emph{analytic form} if it satisfies (1) of Theorem \ref{P:ClassifyingmfR}.
\end{definition}

Our goal in this section is to prove that Case (2) of Theorem \ref{P:ClassifyingmfR}
contradicts the strong minimality of $\CD$. Thus, our negation assumption is that
there exists a matrix $M\in GL(2,R)$ such that all matrices in $M^{-1}\mathfrak R M$
are of the form
 \begin{equation}\label{form} \begin{pmatrix}
a & 0\\
b & a
\end{pmatrix}\end{equation} for $a,b\in R$.

Let us first note that we may assume that all matrices in $\mathfrak R$ itself are
in form (\ref{form}). Indeed, if $f\in \mathfrak F$, we  have defined  $J_0 f$ with respect to some fixed atlas giving $G$ its definable differentiable manifold structure.  Let $g$ be the chart in that atlas mapping a neighborhood $U\subseteq R^2$ onto a neighborhood of $0$. Consider $h: M^{-1}U\to G$ given by $x\mapsto g(Mx)$. Since $h$ is a diffeomorphism, we can replace $g$ in our atlas with $h$. Denoting $\mathfrak R_h$ the ring $\{J_0(f): f\in \mathfrak F\}$ with respect to this new atlas we see that $\mathfrak R_h=M^{-1}\mathfrak R M$, as needed.

%We may also assume that $MU=U$. Now
%consider the definable bijection $h:G\to G$ which is the identity outside $U$ and
%$h(x)=Mx$ on $U$. The push-forward of the structure $\CD$ under $h$
% is an isomorphic, definable, strongly minimal structure $\CD'$, and it is easy to verify that for
%  any $\CD$-function $f\in \mathfrak F$, its image in $\CD'$
% is a function whose differential at $0$ is $M^{-1}J_0f M$.
%  Thus the ring of Jacobians at $0$ of all smooth $\CD'$-functions $f$ with $f(0)=0$ consists of matrices
% as in (\ref{form}). We now replace $\CD$ with $\CD'$.

In the case where $G=\la R^2,+\ra$, \cite[Corollary 2.18]{HaKo}
immediately eliminates the possibility that $\mathfrak R$ is the ring of upper triangular matrices. The goal of this subsection is to prove an
analogue of that result in the context of an arbitrary group $G$.

We first need the following version of the uniqueness of definable solutions to
definable ODEs. It can be  easily deduced from \cite[Theorem 2.3]{OtPePi}.

\begin{proposition}\label{ODE}
    Let $Gr(k,n)$ be the space of all $k$-dimensional linear subspaces of $R^n$. Let $U\sub R^{n}$ be an open set
     and assume that $L:U\to  Gr(k,n)$ is a definable $\CC^3$-function assigning to each $p\in U$ a $k$-linear space $L_p$.
      Assume that $C_1,C_2\sub U$ are definable $k$-dimensional smooth manifolds such that  for every $p\in C_1\cap C_2$,
       the tangent space of $C_i$  at $p$
       equals $L_p$.
       Then for every $p\in C_1\cap C_2$ there exists a neighbourhood $V\ni p$ such that $C_1\cap V=C_2\cap V$.
\end{proposition}

\begin{defn} {\em A definable vector field} on an open $U\sub G$, is given by a definable partial
function $F:U\to T(U)$ from $U$ to its tangent bundle $T(U)$, such for every $g\in
G$, $F(g)\in T_g(G)$.

 Every definable non-vanishing vector field $F$ on $U$ gives rise to a {\em a definable line
 field}, still denoted by $F$, where to each $g\in U$ we assign the $1$-dimensional subspace of $T_{g}(U)$ spanned
 $F(g)$.

We say that a line field $F$ is (left) $G$-invariant if if for every $g,h \in U$,
$$F(h)=d_{g}(\ell_{hg^{-1}})\cdot F(g).$$

 Given a line field $F$, we
say that a definable smooth $1$-dimensional set $C\sub U$ is {\em a trajectory of
$F$} if for every $g\in C$, the tangent space to $C$ at $g$ is $F(g)$.
\end{defn}
\begin{lemma}\label{traject}  Let $F$ be
a definable non-vanishing $G$-invariant line field.
    Assume that $C \sub G$ is a definably connected smooth $1$-dimensional trajectory
of  $F$. Then $C$ is a coset of a definable local subgroup of $G$.
\end{lemma}
\begin{proof}  Recall that we identify an open neighborhood $U$ of $0$ with an
open subset of $R^2$, and  $T(U)$ is identified with $U\times R^2$. The line field
can be viewed as a map $F:U\to Gr(1,2)$.

It will suffice to show that if $h\in C$, then $h\ominus C$ is a local subgroup. Hence we may assume that $0\in C$. Since $F$ is left-invariant, for any $g\in G$,
 $g\oplus C$ is also a
trajectory of $F$. By Proposition \ref{ODE}, $C$ and $g\oplus C$ coincide on some
neighborhood of $g$, provided that $g\in C$. It follows that every $x\in C$  and $g\in C$ sufficiently
small, we also have $x\oplus g\in C$. Thus $C$ is a local subgroup of $G$.
\end{proof}

We can now return to our main goal: proving that $\mathfrak R$ is of analytic form.
Recall that we assume that for a $\CD$-function $f$ and $b\in \dom(f)$ we can write
\begin{equation}\label{E:a-b}
\tilde J_b(f)=\begin{pmatrix}
\alpha_f(b) & 0\\
\beta_f(b) & \alpha_f(b)
\end{pmatrix}.
\end{equation} When $f$ is clear from the context we  omit the subscript $f$.

Let $v_0=\begin{pmatrix} 0 \\ 1
\end{pmatrix}\in T_0(G)$ and  consider the non-vanishing $G$-invariant vector field $F$ given by
\[\{d_0 \ell_b \cdot v_0:b\in G\}.\]
For $b\in G$, let $v_b=d_0\ell_b \cdot v_0\in T_b(G)$.

\begin{lemma}\label{bad1} For  every $\CD$-function $f:U\to G$ and $b\in \dom(f)$, we have
\[d_b f \cdot v_b=\alpha (b) v_{f(b)}\in R v_{f(b)},\] namely the line-field induced by $F$ is invariant under $df$.

If, in addition, $\alpha(b)=0$, then
    $$d_b f\cdot T_b(G)\sub Rv_{f(b)}.$$
\end{lemma}
\begin{proof}
By assumption on the form of matrices in $\mathfrak R$, we have $\tilde d_b f \cdot
v_0 =\alpha(b)\cdot  v_0$. Writing $\tilde d_b f$ explicitly (and composing on the
left with $d_0\ell_{f(b)}$), we obtain
\[(d_bf) (d_0 \ell_b) \cdot v_0=\alpha (b) (d_0\ell_{f(b)}\cdot v_0),\] which implies the first
clause.

For the second clause, notice the special form of $\tilde d_bf$ implies that when
$\alpha(b)=0$,  then for every $v\in T_0(G)$, we have $\tilde d_b f\cdot v\in R
v_0$. The result easily follows.
\end{proof}

\begin{lemma}\label{bad1.5} Assume that $f:U\to G$ is a $\CD$-function, and that $C\sub U$ is a definable smooth curve
which is a trajectory of $F$. Then so is $f(C)$.
\end{lemma}
\begin{proof} By the first clause of Lemma \ref{bad1}, the image of $C$ under $f$ is also a
trajectory of $F$.
\end{proof}

\begin{lemma} \label{bad2} Assume that $f$ is a $\CD$-function, and $C\sub \dom(f)$ is a definable smooth curve
 such that at every $b\in C$ we have $\alpha(b)=0$ (in the above notation). Then for
    every generic $b\in C$, the tangent space of  $f(C)$ at $f(b)$ is the $R$-span of $v_{g(b)}$. Namely,
 $f(C)$ is a trajectory of $F$, in a neighborhood of
    $f(b)$.\end{lemma}

\begin{proof}
 Consider the restriction of $f$ to $C$, and pick a generic $b$ in $C$. Since
$b$ is generic, the map $f|C\, :\,C\to f(C)$ is a submersion, namely
$T_{f(b)}(C)=d_b f\cdot T_b(C)$. By the second clause of Lemma \ref{bad1}, we
conclude that $T_{f(b)}(f(C))$ equals $Rv_{f(b)}$.
\end{proof}

\begin{lemma} \label{bad3} There exists a $\CD$-function $h$ and a  definable curve $C\sub
    G$ such that $h(C)$ is a trajectory of $F$.
\end{lemma}

\begin{proof} This is similar to the proof of the  claim in Proposition \ref{thm-dim}. Fix any $\CD$-function $f$,
 and $\alpha_f$ as in (\ref{E:a-b}) above.  \\

 \noindent\textbf{Claim.} There is $a_0\in G$ such that for infinitely many $b\in G$, we have $\alpha_f(b)=\alpha_f(b\oplus a_0)$.
 \begin{proof}[Proof of Claim]
For $r\in R$, let $C_r=\{b\in G:\alpha_f(b)=r\}.$ Pick $r$ generic in the
image of $\alpha_f$. By Proposition \ref{thm-dim} (and in the notation of that proposition), $\dim(J(U))=2$. Therefore,  genericity of $r$ implies that $C_r$ is $1$-dimensional. Consider $D_r=C_r\ominus
C_r$. If $D_r$ is still 1-dimensional, then as we have already seen several times, $C_r$ is
contained in a coset of a $\bigvee$-definable subgroup $H_r$ and then picking $a_0\in H_r$ small enough will work with any $b\in C_r$.

Otherwise, $D_r$ is $2$-dimensional. We may now pick $a_0\in D_r$ generic over $r$.
Since $r$ is still generic over $a_0$, there are infinitely many $r'$ such that
$a_0\in D_{r'}$. For each such $r'$, there exists $b\in C_{r'}$ with $b\oplus a_0\in
C_{r'}$. \end{proof}

Fix $a_0$ as above, and consider the $\CD$-function $h(x)=f(x\oplus a_0)\ominus
f(x)$. It is easily verified that for each $b\in G$ we have $$\tilde J_b h=\tilde
J_b(f(x\oplus a_0)\ominus f(x))=\tilde J_b f(x\oplus a_0)- \tilde J_b f(x).$$

It follows that $\alpha_h(b)=0$ for every $b\in G$ such that  $\alpha_f(b\oplus
a_0)=\alpha_f(b)$. Let $C$ be the collection of all those elements $b$.  By the claim, $C$ is a curve. By Lemma
\ref{bad2}, the curve $h(C)$ is a trajectory of $F$ near $b$.
\end{proof}

We can now conclude the following theorem.
\begin{theorem}\label{P:ClassifyingmfRbad}
    The ring  $\mathfrak R$ is of analytic form.
\end{theorem}
\begin{proof} We still work under the negation assumption that we are in Case 2 of
Theorem \ref{P:ClassifyingmfR}.    Using Lemma \ref{bad3} and Lemma \ref{traject} we
obtain a definable local
    subgroup $H$ which is a trajectory of the vector field $F$, and thus all of its
    cosets are also trajectories of $F$. Let $U$ be a neighbourhood of $0$ which can be
    covered by cosets of $H$, all trajectories of $F$.

Fix any $\CD$-function $f\in \mathfrak F$ which is not $G$-affine. By Lemma
\ref{bad1.5}, for every $a\in U$ such that $f(a)\in U$,  the image $f(H\oplus a)$ is
also a coset of $H$. Fix $a_0\in H$ close enough to $0$ and consider the
$\CD$-function $k(x)=f(x)\ominus f(x\oplus a_0)$. Since $f$ is not $G$-affine, the
function $k$ is not constant.

Notice that for every $x$ sufficiently close to $0$, the elements $x$ and $x\oplus
a_0$ belong to the coset $x\oplus H$, and therefore as we just noted, $f(x)$ and
$f(x\oplus a_0)$ belong to the same coset of $H$. It follows that $k(x)\in H$ and
therefore $k$ sends an open subset of $G$ into $H$, contradicting strong minimality
(the pre-image of some point will be infinite).
\end{proof}

Note that the above argument does not really use the definability of the trajectory
$C$ but merely its existence. Thus, if we worked over the reals, then we could have
used the usual existence theorem for solutions to differential equations in order to
derive a contradiction.

\section{Some intersection theory for $\CD$-curves}

\label{sec-intersection}

Our ultimate goal is to show, under suitable assumptions, that if two plane curves
$C, D\sub G^2$ are tangent at some point $p$, and $C$ belongs to a $\CD$-definable
family $\CF$ of plane curves, then by varying $C$  within $\CF$ one gains additional
intersection points with $D$, near the point $p$ (see Proposition
\ref{intersection}~(2)). This will allow us to detect tangency $\CD$-definably.

 The main tool towards this end is  the following theorem, whose proof will be
 carried out in this section via a sequence of lemmas.

\begin{theorem}\label{main-noninject} Assume that $f$ is in $\mathfrak F$.
If $J_0(f)=0$, then there is no neighborhood of $0$ on which $f$ is injective.
\end{theorem}

We now digress to report on an unsuccessful strategy, which nevertheless may be of
some interest.

\subsection{Digression: on almost complex structures}
\label{almost-complex}
 Let $K=R(\sqrt{-1})$. In analogy to the notion of an almost complex structure on a real manifold, we
call {\em a definable almost $K$-structure} on a definable $R$-manifold $M$, a
definable smooth linear $J:TM\to TM$ sending each $T_x(M)$ to $T_x(M)$, such that
$J^2=-1$.

Note that every definable $K$-manifold admits a natural almost $K$-structure,
induced by multiplication of each $T_x(M)$ by $i=\sqrt{-1}$. It is known that when
$K=\mathbb C$ any $2$-dimensional almost complex structure is isomorphic, as an
almost complex structure, to a complex manifold. The proof of this result seems to
be using integration and thus we do not expect it to hold for almost $K$-structures
in arbitrary o-minimal expansions of real close fields.

Returning now to our $2$-dimensional group $G$, we can endow $G$ with a definable
almost $K$-structure in the following way. Just as we did at the beginning of
Section \ref{sec:from ring to field}, we may first assume that every matrix in
$\mathfrak R$ has the form $\begin{pmatrix}
    a & -b \\
    b & a
    \end{pmatrix}$.
 Next, we
identify naturally $T_0(G)$ with $R^2\sim K$ and let $J:T_0(G)\to T_0(G)$ be defined
by $J(x,y)=(x,-y)$. Next, use the differential of $\ell_a$ to obtain $J:TG \to TG$
as required. Note that since $TG$ is a trivial tangent bundle, this step can be
carried out for any definable group of even dimension.  However in the case of $G$,
our choice of $J$ and the fact that for each $\CD$-function $f$, $\tilde J_af$ has
analytic form, implies that $f$ is so-called $J$-holomorphic, namely that each each
$a\in \dom f$ we have
$$J\circ d_af=d_{f(a)}\circ J.$$

Now, if our underlying real closed field $R$ were the field of real numbers, then $G$
would be isomorphic as an almost complex structure to a complex manifold $\hat G$,
and this isomorphism would send every $J$-holomorphic function from $G$ to $G$ to a
holomorphic function from $\hat G$ to $\hat G$. In particular, by our above
observation every $\CD$-function would be sent to a holomorphic function. This would
give an immediate proof of Theorem \ref{main-noninject}, due to the fact that the
result is true for holomorphic maps.

Unfortunately, we do not know how to prove for arbitrary $K$ that every
$2$-dimensional almost $K$-manifold is (definably) isomorphic to a $K$-manifold, and
hence we cannot use the theory of $K$-holomorphic maps in order to deduce Theorem
\ref{main-noninject}. We thus use a different strategy.

\subsection{A motivating example}
If $f$ were holomorphic, then the above theorem would follow from the argument
principle and the open mapping theorem. Since our functions are not necessarily holomorphic, we describe  a different, more topological proof of Theorem \ref{main-noninject} for an analytic function  $f$: we let $h(z)=f(z)/z$ (complex division)
 for $z\neq 0$ and $h(0)=0$. The assumption that $J_0f=0$ implies that $h$ is continuous at
 $0$ and hence holomorphic.  Thus $h$ is either locally constant or an open map in a
 neighborhood of $0$. Now, if $h$ were locally constant, then $f\equiv 0$ near $0$
 and thus clearly non-injective, so assume that $h$ is an open map.

We now consider the complex function
 $M(z,w)=z\cdot w$, and for $a,b\in \mathbb C$ near $0$, let
 $M_{a,b}(z)=M(z-a,h(z)-b)$. Notice that $M_{0,0}(z)=f(z)$.
 Let $\deg_0(f)$ be the local degree of $f$ at $0$ (see details below).
 Since the local degree is preserved under definable homotopy (see Fact \ref{degree} below), it follows from the general theory that  $\deg_0(M_{a,b})=\deg_0(f)$ for sufficiently
 small $a,b$.
Because each $M_{a,b}$ is
 holomorphic, the sign of $|J_zM_{a,b}|$ is positive at a generic $z$ in a small disc around $0$,
 and therefore
 $$\deg_0(M_{a,b})\geq |M^{-1}_{a,b}(w)|,$$ for all $w$ close to $0$.

If we take  $w=0$,  then we get
 \[|M^{-1}_{a,b}(0)|\geq 2\]
  (the points $a$ and $h^{-1}(b)$ being two such pre-images),
 implying $\deg_0(f)=\deg_0(M_{a,b})\geq 2$. This implies
 that  $f$ is not
 locally injective near $0$.\smallskip

Our objective is to imitate the above proof, using $\CD$-functions instead of
holomorphic ones. The main obstacle is the fact that we do not have multiplication
or division in $\CD$, so we want to produce a $\CD$-function which sufficiently
resembles the multiplication function $M$.

\subsection{Topological preliminaries}
Throughout this section we will be using implicitly the o-minimal version of
Jordan's plane curve theorem (see \cite{WoePhD}). We  recall some definitions
and results (see \cite[Section 2.2-2.3]{PeStExpansions}). Given a circle $C\sub
R^2$, a definable continuous $f:R^2\to R^2$  and $w\notin f(C)$, we let $W_C(f,w)$
denote the winding number of $f$ along $C$ around $w$. If $f^{-1}(w)$ is finite,
$p\in R^2$ and $f(p)=w$, then  $\deg_p(f)$ is defined to be $W_C(f,f(p))$ for all
sufficiently small $C$ around $p$. We need the following results.

\begin{fact} \label{degree}
    Let $C\sub R^2$ be a circle oriented counter clockwise.
\begin{enumerate}
\item  If $\{f_t:t\in T\}$ is a definable continuous family of functions with $w
\notin f_t(C)$ for any $t\in T$ and $T$ definably connected, then
$W_C(f_{t_1},w)=W_C(f_{t_2},w)$ for all $t_1,t_2\in T$.

\item Assume that $C$ is a circle around $p$, $f:C\to R^2$ definable and continuous,
and $w_1, w_2$ are in the same component of $R^2\setminus f(C)$. Then
$W_C(f,w_1)=W_C(f,w_2)$.

\item If $f$ is definable and $R$-differentiable at $p$ and  $J_p(f)$ is invertible,
then $\deg_p(f)$ is either $1$ or $-1$, depending on whether $|J_p(f)|$ is positive
or negative.

\item Assume that $f$ is a definable $\CM$-smooth, open map, finite-to-one in a
neighbourhood $U$ of $p$ and that $f(z)\neq f(p)$ for all $z\neq p$ in $U$. Assume
also that $J_z(f)$ is invertible of positive determinant for all generic  $z \in U$.

Let $C\sub U$ be a circle around $p$. Then for all $w\in f(\intr(C))$, if $w$ and
$f(p)$ are in the same component of $R^2\setminus f(C)$, then $W_C(f,f(p))\geq
|f^{-1}(w)\cap \intr(C)|$, and if $w$ is also generic, then
$W_C(f,f(p))=|f^{-1}(w)\cap \intr(C)|$.

%\item Assume that  $f:U\to R^2$ is a definable continuous map on an open $U\sub
%R^2$. If $f$ is locally injective at all but finitely many points of $U$ then $f$ is
%an open map.
\end{enumerate}
\end{fact}

\proof (1) follows from \cite[Lemma 2.13(4)]{PeStExpansions}. (2) is just
\cite[Lemma 2.15]{PeStExpansions}. The proof of (3) is the same as the classical
one, so we omit it.

% The result is easy for linear maps. Next, let $T(z)=(J_pf)\cdot z$. One first
%shows that for a sufficiently small $r$, the restrictions of  $f(z)-f(p)$ and $T(z)$
%to $C_r$, the circle around $p$ of radius $r$, are homotopic in $R^2\setminus
%\{0\}$. It follows that $\deg_p(f)=\deg_p(T)$.

(4)  It follows from (2) that $W_C(f,f(p))=W_C(f,w)$. We let $\{z_1,\ldots,
z_k\}=f^{-1}(w)\cap \intr(C)$.  By \cite[Lemma 2.25]{PeStExpansions},
$W_C(f,w)=\Sigma_{i=1}^k \deg_{z_i}(f),$ so it is sufficient to see that
$\deg_{z_i}(f)\geq 1$, for each $i$. We  fix a small circle, $C_i$, around $z_i$ such
that $W_{C_i}(f,w)=\deg_{z_i}(f)$, and then fix a generic $w_0\in f(\intr(C_i))$
sufficiently close to $w$, so in particular, the Jacobian of $f$ at each pre-image
of $w_0$ is invertible of positive determinant. By \cite[Lemma
2.25]{PeStExpansions}, $\deg_{z_i}(f)=\Sigma_j \deg_{p_j}(f)$, where the $p_j$ are the
pre-images of $w_0$ in $\intr(C)$. By (3), for each $p_j$, we have $\deg_{p_j}(f)=1$,
thus $\deg_{z_i}(f)=|f^{-1}(w_0)|\geq 1$.

The same argument shows that for generic $w_0$ near $p$, we have
$W_C(f,f(p))=|f^{-1}(w_0)|.$ \qed

 \subsection{Back to $\CD$-functions}

We still identify an open neighborhood of $G$ with an open subset of $R^2$ and
identify $0_G$ with $0=(0,0)$. For an open set $U\sub G$ and a  function $f:U\to G$,
sending $x_0$ to $y_0$, we say that \emph{$f$ is generically $k$-to-$1$ at $x_0$} if
for every open $V\ni x_0$ and $W\ni y_0$ there exists an   open $y_0\in W_0\sub W$
such that for any generic $y\in W_0$, $|f^{-1}(y)\cap V|=k$.

Below we use the notion of a $\CD$-function $M$ from an open $U\sub G^2$ into $G$.
By that we mean that there exists a $\CD$-definable set $S\sub G^2\times G$ of
Morley rank $2$ containing the graph of $M$.
\begin{lemma}\label{fact4} Let $U\sub G^2$ be a definable open neighborhood of
$(0,0)$.
    Assume that $M:U\to G$ is a continuous $\CD$-function such that $M(0,y)=M(x,0)=0$ for all $x,y$ close enough to $0$.
    Assume that $f,h\in \mathfrak F$ and we have:
    \begin{enumerate}
    \item  For every $a,b$ in some neighborhood of $0$, the function $g_{a,b}(x)=M(f(x)\om a, h(x)\om b)$ is not locally constant
    near $0$.
    \item $f$ and $h$ are, respectively,  generically $k$-to-one and  $m$-to-one near $0$,
    \end{enumerate}

    Then $g(x)=M(f(x),h(x))$ is, generically, at least $k+m$-to-one near $0$.
\end{lemma}
\begin{proof} By assumptions on $f,h$ and $M$, for every $a,b$ in some neighborhood of $0$,
the function $g_{a,b}$ is a non-locally constant $\CD$-function on some neighborhood of $0$, namely it is continuous and its graph is
contained in a rank one $\CD$-definable set. By Theorem \ref{C:continuous}, it is open as well.
Since it is definable in $\CD$ and not locally constant,  it is finite-to-one near $0$.
    Also, it follows from
     Corollary \ref{C:no-neg} and Proposition \ref{P:ClassifyingmfRbad} that $g_{a,b}$
     has positive determinant of the Jacobian at every point where the Jacobian matrix does not
     vanish, which by Lemma \ref{C:sgn} is a co-finite set.

    We now fix a simple closed curve $C$ around $0$ such that $0\notin g(C)=g_{0,0}(C)$ and $\deg_0(g)=W_C(g,0)$. By continuity of $M$ and $g$ we can find an open
    $U_1\ni 0$, and an open disc $U_2\ni 0$, such that for all $a,b\in U_1$, $g_{a,b}(0)\in U_2$ and  $g_{a,b}(C)\cap U_2=\emptyset$. It follows that $g_{a,b}(0)$ and $0$ are in the
    same component of $R^2\setminus
g_{a,b}(C)$.

Take $a,b\in U_1$ independent generics. By Fact \ref{degree},
$$\deg_0(g)=W_C(g,0)=W_C(g_{a,b},0)=W_C(g_{a,b},g_{a,b}(0))\geq |g_{a,b}^{-1}(0)| .$$

Because  $a,b$ are independent generics $f^{-1}(a)\cap
h^{-1}(b)=\emptyset$. Also, by our assumptions on $M$ and the definition of
$g_{a,b}$, we have $$f^{-1}(a)\cup h^{-1}(b)\sub g_{a,b}^{-1}(0).$$
 Hence,  $|g^{-1}_{a,b}(0)|\ge m+k$. It follows from
 Fact  \ref{degree} (4) that  $\deg_0(g)\geq m+k$ and that $g$ is generically, at least, $k+m$-to-one near $0$.
\end{proof}

\subsection{Producing the function $M$}\label{Producing M}

We now proceed to construct the desired $\CD$-function $M$ as in Lemma \ref{fact4}.
We start with a
  $\CD$-function $k(x)$ which is not $G$-affine and fix a  generic $a_0\in \dom k$.
  Define
$$M(x,y)=(k(a_0\op x\op y)\om k(a_0\op x))\om (k(a_0\op y)\om k(a_0)).$$
We write $M_a(y)=M(a,y)$.

By definition, we have
\begin{description}
    \item[(A)]  For $x,y$ near $0$,  $M(0, y)=M(x, 0)=0$.
\end{description}

Our next goal is to show that $M$ can be used, similarly to multiplication, to
``divide (an appropriate) function $f$ by $x$''. Namely, that we can implicitly
solve $M(x,y)=f(x)$ in some neighborhood of $x=0$. This is the purpose of the next
few results.

By Theorem \ref{P:ClassifyingmfRbad} and the discussion in Section \ref{sec:from
ring to field}, we may assume that for a smooth $f\in \mathfrak F$, the matrix
$J_0(f)$ has the form
$$\left(
\begin{array}{cc}
  c& -e \\
  e & c \\
\end{array}\right),$$ with $c,e\in R$.

 We consider the partial definable map $d:G\to R^2$, mapping $a$ to the first column of the Jacobian matrix  $J_0(M_a)$ (so if $J_0(M_a)=\begin{pmatrix}
 c & -e \\
 e & c \\
 \end{pmatrix}$, then $d(a)=(c,e)$). Note that $d(a)$ completely determines $J_0(M_a)$, and in particular $d(0)=0$ if and only if $J_0(M)=0$.  By Lemma \ref{basics} (and using the fact that
$\tilde J_0(f)=J_0(f)$), we get that $J_0(M_a)$ is equal to:
\[J_0 M_a =\tilde J_0([k(a_0\oplus a\oplus y)\om k(a_0\op a)]) -\tilde J_0 ([k(a_0\op
y)\om k(a_0)])=\tilde J_{a_0\op a}(k)-\tilde J_{a_0}(k). \tag{$\dagger$}\] By
Proposition \ref{thm-dim}, applied to $k(x)$, the image of every open $U\ni 0$ under
$x\mapsto \tilde J_x(k)$ is a 2-dimensional subset of $\mathfrak R$, hence by
o-minimality this map is locally injective near the generic $a_0$. Equivalently, the
map $x\mapsto \tilde J_{a_0\op x}(k)$ is locally injective near $0$. Since $\tilde
J_{a_0}(k)$ is constant, it follows that $d(x)$ is locally injective at $0$. In
particular, we have

\begin{description}
    \item[(B)] $d(0)=0$, and there is a neighborhood of $0$ where $d(a)\neq 0$ for all
 $a\neq 0$.
    \end{description}.

We are going to use several different norms in the next argument, so we set
$$||(x,y)||=\sqrt{x^2+y^2},$$ and for a linear map $T$ we denote the operator norm by
$$||T||_{\textup{op}}=\max\{||T(x)||/||x||:x\neq 0\}.$$ Observe that   $||d(a)||=||J_0(M_a)||_{\textup{op}}$. It is well-known (and easy to see) that if we
identify every linear map with a $2\times 2$ matrix, then $||T||_{\textup{op}}$ and $||T||$
are equivalent norms.

 We need an additional property of $M$. Given two
  functions $\alpha,\beta:U^*\to R^{\geq 0}$ on a punctured neighborhood $U^*\sub R^2$ of  $0$, we write $\alpha\sim
 \beta$ if $\lim\limits_{t \to 0} \alpha(t)/\beta(t)$ is a positive element of $R$. We
 will show:
\begin{description}
 \item[(C)] There are definable $R^{>0}$-valued functions $e(a)$ and $\delta(a)$,
 in some punctured neighborhood $U^*$ of $0$, with $e(a)\sim
    ||d(a)||$ and $\delta(a)\sim ||d(a)||^2$, such that for every $a\in U^*$, the function
    $M_a=M(a,-)$ is invertible on the disc $B_{e(a)}$ and its image contains
    the disc $B_{\delta(a)}$ (recall that for $a=0$ we have $M_a(x)\equiv 0$ near $0$).
\end{description}

  In order to prove (C), we use an
  effective version of the inverse function theorem, as appearing in \cite[\S7.2]{vdDries}. We give the
  details, with references to \cite{vdDries}.

\begin{proposition}\label{IFT}  There exists a constant $C>0$ such that setting $e(a)=||d(a)||/4C$ and $\delta(a)=e(a)^2/2$, we have that for all $a$ in a small  %definable
%functions $e(a)$ and $\delta(a)$ from a
punctured neighborhood of $0$
%into $R^{>0}$, such that for $a\neq 0$ in a
%small neighborhood of $0$,
the function $M_a(y)$ is injective on $B(0;e(a))$ and
its image contains a ball of radius $\delta(a)$ around $0$. %Furthermore, there is a
%    constant $C>0$ such that  $e(a)=||d(a)||/4C$ and $\delta(a)=e(a)^2/2$.
\end{proposition}

\begin{proof} We start with some
observations. If
$$A=J_0(M_a)=\left(
\begin{array}{cc}
  c & -e \\
  -e & c \\
\end{array}\right)$$
then $||A||_{\textup{op}}=\sqrt{c^2+e^2}=||d(a)||.$ And if $A$ is invertible, then
$||A^{-1}||_{\textup{op}}=1/||d(a)||$.
%Furthermore, the norm $||A||_{\textup{op}}$ is equivalent to its Euclidean norm. {\bf I got
%very confused by the various norms since I am not sure which norms vdD is using in
%the text we refer to. There are three norms around $| |$, $|| ||$, and $|| ||_{\textup{op}}$.
% Can we clarify this point?}

Consider  the partial map $D:G\times G\to R^4$, defined by $D(a,y)=J_y(M_a)\in
M_2(R)$. For each $a,y$, we view $D(a,y)$ both as a linear operator and a vector in
$R^4$.
%(in fact the image of the function $d(y,a)$ is contained in a 2-dimensional subspace of the $2\times2$ matrices).
Since $M$ is a $\CC^2$-function, $||J_{(a,y)}D||_{\textup{op}}$ is bounded by some constant
$C$, as $(a,y)$ varies in a neighborhood $B_1\times B_2$ of $(0,0)$, and we may assume that $C>1$. By \cite[Lemma 7.2.8]{vdDries} applied to
$D$, for every $(a_1,y_1), (a_2,y_2)\in B_1\times B_2$ we have
\[\tag{*}||J_{y_1}(M_{a_1})-J_{y_2}(M_{a_2})||<C||(a_1,y_1)-(a_2,y_2)||.\]

Note also that $D(0,0)=J_0 M_0=0$, so restricting further $B_1, B_2$ we may also
assume that $||D(a,y)||<1$ for all $(a,y)\in B_1\times B_2$.

We now need a version of \cite[Lemma 7.2.10]{vdDries}.

\begin{lemma}\label{invertible2}
    For every $a\in B_1$ such that $J_0 M_a$ is invertible, and for all $y_1,y_2\in B_2$,
    if $||y_1||,||y_2||\leq e(a)$,
     then
    \begin{enumerate}

        \item the matrices $J_{y_1}M_a, J_{y_2}M_a$ are invertible.

        \item
        $||M_a(y_1)-M_a(y_2)||\geq
        e(a)||y_1-y_2||.$ In particular, $M_a$ is injective on the disc $B_{e(a)}$.
    \end{enumerate}
%    Moreover, we can choose $e(a)=||d(a)||/C$ for some constant $C>0$.
\end{lemma}
\begin{proof} We fix  $a$ with $J_0(M_a)$ invertible and we write $J_0 M_a=\begin{pmatrix}
c & e\\
-e & c
\end{pmatrix}$.
 By (*), for every $y\in B_2$ and for every $E>0$, if $||y||<E/2C$, then
$$||J_{y}M_a-J_0M_a||\leq C||y||\leq C||(y,a)-(0,a)|| <\frac{E}{2}.$$
 In particular, since $J_0 M_a\neq 0$,  we may take $E<||d(a)||=| |J_0M_a||_{\textup{op}}$ and then  $J_y M_a$ must be non-zero. Because
$M_a$ is a $\CD$-function it follows that $J_yM_a$ is invertible.

Let $c'=1/||J_0(M_a)^{-1}||_{\textup{op}}$. As we pointed out earlier, in our case
$$||J_0(M_a^{-1})||_{\textup{op}}=||J_0(M_a)^{-1}||_{\textup{op}}=1/||d(a)||,$$ hence $c'=||d(a)||$.
Now, for all non-zero vectors $w$, we have $||J_0(M_a)^{-1}(w)||\leq \frac{1}{c'}
||w||$, so by substituting $w$ with $J_0(M_a)^{-1}(z)$, we get $c'||z||=||d(a)||
\cdot ||z|| \leq ||J_0M_a (z)||$.

 Hence, for any two $y_1,y_2\in
R^2$:
\[\tag{**}||J_0M_a \cdot (y_1-y_2)||\geq ||d(a)||\cdot ||y_1-y_2||.\]

By \cite[Lemma 7.2.9]{vdDries}, applied to the function $M_a$,  we also have for all
$y_1,y_2\in B_2$,
\[
||M_a(y_1)-M_a(y_2)-J_0 M_a (y_1-y_2)||\le ||y_1-y_2||\max_{t\in [y_1,y_2]} ||J_t
M_a - J_0 M_a||_{\textup{op}},
\] where $[y_1,y_2]$ is the line segment in $R^2$ connecting $y_1$ and $y_2$. Hence,
by the triangle inequality,
\[
||M_a(y_1)-M_a(y_2)||\ge ||J_0 M_a (y_1-y_2)||-||y_1-y_2||\max_{t\in [y_1,y_2]}
||J_t M_a - J_0 M_a||_{\textup{op}}.
\]

Putting this together with (*) and (**), we have: if $y_1,y_2\in B_2$ and
$||y_i||<E/2C$, for $i=1,2$, then
\[
||M_a(y_1)-M_a(y_2)||\ge  (||d(a)||-C||y_1-y_2||)||y_1-y_2||.
\]
If in addition $||y_1-y_2||<\frac{||d(a)||}{2C}$, then \[ \tag{***}
||M_a(y_1)-M_a(y_2)||\ge
(||d(a)||-\frac{||d(a)||}{2})||y_1-y_2||=\frac{||d(a)||}{2}||y_1-y_2||.
\]

We summarize what we have shown so far: for any $E<||d(a)||$, if
$||y_1||,||y_2||<E/2C$ and $||y_1-y_2||<||d(a)||/2C$, then $J_y(M_a)$ is invertible
and (***) holds.

We now fix the parameters as follows: set $E=||d(a)||/2$, $e(a)=||d(a)||/4C=E/2C$.
So, if $||y_1||, ||y_2||< E/2C$, then $||y_1-y_2||<E/C=||d(a)||/2C$, so we may apply
(***) and conclude that $J_{y_i}$ are invertible for $i=1,2$ and
\[
||M_a(y_1)-M_a(y_2)||\ge \frac{||d(a)||}{2}||y_1-y_2||\geq e(a)||y_1-y_2||.
\]
\end{proof}

By the proof of \cite[Theorem 2.11]{vdDries},
\[
    \{y: ||y-M_a(0)||<\frac{e^2(a)}{2}\}\subseteq \{M_a(z): ||z||<e(a)\}
\]
(apply the claim on the second line of p.113 with $\epsilon, c$ there both
substituted with $e(a)$ here, and our $M_a$ substituting $f$ there). Thus, the image
of the disc $B_{e(a)}$ under $M_a$ contains a disc of radius $\frac{e(a)^2}{2}$
around $M_a(0)=0$. We do not repeat the proof here.

%Setting  $\delta(a)=\frac{e(a)^2}{2}$
% completes the proof of Proposition \ref{IFT}.
\end{proof}

\subsection{Proving Theorem \ref{main-noninject}}
 We now fix a $\CD$-function
$M:G^2\to G$ satisfying conditions (A), (B) and (C) as above, with  $d(x)$ the first column of $J_0 M_x$.
We first need a simple observation.
\begin{fact}\label{square} Assume that $f:U\sub R^2\to R^2$ is a definable $C^2$-function sending $0$ to $0$.
If $J_0f=0$, then $\lim\limits_{x\to 0}||f\circ f(x)||/||x||^2=0$.
\end{fact}
\begin{proof}
    As already mentioned above, the operator norm and the Euclidean norm on $R^2$ are equivalent -- so we may work with either.
%We first compare the two norms on $A\in M_2(R)$, namely the operator norm $||A||_{\textup{op}}$ and
%the Euclidean norm $||A||$, where $A$ is viewed as an element of $R^4$. It is not
%hard to see that the two are equivalent, namely there are constants $c_1,c_2>0$ such
%that for every $A$, $c_1||A||\leq ||A||_{\textup{op}}\leq c_2||A||$.

We first claim that there exists some neighborhood $U$ of $0$, and a constant $C$ such that for all $x\in U$, $||f(x)||\leq C||x||^2.$

We use the following corollary of \cite[Lemma 7.2.8]{vdDries}: Let $g$ be a definable $C^1$-map from an open ball $B\sub R^m$ centered at $0$ into $R^n$ with $g(0)=0$. Then for all $x\in B$,
 \begin{equation} \label{bounds} ||g(x)||\leq \sup\limits_{a\in B(0;||x|)} ||J_a g|| \,||x||. \end{equation}

We now consider the map $\alpha(a)= J_a(f) $, as a map from an open ball $B$ around $0\in G$ (identified with a ball centered at $(0,0) \in R^2$ ) into $R^4$. Since $f$ is a $C^2$-map the map $\alpha$ is a $C^1$-map and hence, by (\ref{bounds}), there is some constant $C$ (a bound on the norm of $d_a(\alpha)$ as $a$ varies in $B$), such that for all $a\in B$, $$||J_af||=||\alpha(a)|| \leq C ||a||.$$
It follows that for all $x\in B$, $\sup\limits_{a\in B(0;||x|)} ||J_a f||\leq C ||x||$.

Next, we apply (\ref{bounds}) to the map $f$ itself and conclude, using what we have just shown, that for all $x\in B$, $$||f(x)||\leq \sup\limits_{a\in B(0;||x|)} ||J_a f|| \,  ||x||\leq C||x||^2.$$ This ends the proof of our first claim.

It now follows that $||f( f(x))||\leq C ||f(x)||^2\leq  C^2||x||^4$. Thus
$\lim\limits_{x\to 0}||f\circ f(x)||/||x||^2=0$.
\end{proof}

%Recall that for a $C^{k+1}$ vector field $F: \CM^n\to \mathcal \CM^m$ we have a Taylor expansion around $0$ of the form
%\begin{align*}
%F(x) = & F(0)+\sum_{i=1}^n F_{x_i}x_i+\frac{1}{2}\sum_{i,j=1}^{n}F_{x_ix_j}x_ix_j+\dots \\ & +
%\frac{1}{k!}\sum_{i_1,\dots, i_k=1}^nF_{x_{i_1}\dots x_{i_k}}x_{i_1}\cdot \dots \cdot x_{i_k}+R_k(x)
%\end{align*}
%where the remainder term $R_k(x)$ satisfies $|R_k(x)|\le \frac{M}{(k+1)!}|x|^{k+1}$ for
%$M\ge \max\limits_{|z|\le |x|} \left |\sum\limits_{i_1,\dots, i_{k+1}=1}^nF_{x_{i_1}\dots x_{i_{k+1}}}(z)\right |$.

%Under the assumptions that $f(0)=0$ and $J_0(f)=0$ we get that $f(x)=R_2(x)$. Substituting $R_2(x)$ for $x$ we get
%that $f\circ f(x)=R_2(R_2(x))$ and $|R_2(R_2(x))|\le O(|x|^4)$, with the desired conclusion. \qed \\

%\begin{remark}
%    If $f: \CM^2\to \CM^2$ is differentiable at $0$ and $f(0)=0$ then by definition
%    $J_0 f =\lim\limits_{x\to 0}\frac{f(x)}{|x|}$.
%\end{remark}

%Applying the above observation to $d$, we get:
We also need the following lemma.

\begin{lemma}\label{limit1} If $x(t):(a,\epsilon)\to  R^2$ is a definable curve tending to $0$ as $t\to 0$,
then $\lim\limits_{t\to 0} \frac{||d(x(t))||}{||x(t)||}\neq 0$.
\end{lemma}
\begin{proof} Recall that $d$ is  a map from $U$ into $R^2$ mapping $a$ to the first column of $J_0 M_a$, and  recall that $d(a)$ completely determines $J_0(M_a)$.
%Since   $J_0(M_a)=\left(%
%\begin{array}{cc}
%  c & -e \\
%  e & c \\
%\end{array}%
%\right)$, we have
%$d(a)=(c,e)$.
%Since $d(a)$ completely determines $J_0(M_a)$ we we may also view $d$ as a map from $G$ into $\mathfrak R$.

We claim that $J_0(d)$  is invertible. Indeed, we have seen in
($\dagger$) above (Section \ref{Producing M}), that $J_0(M_a)=\tilde J_{a_0\op
a}k-\tilde J_{a_0}k$.
     By Proposition \ref{thm-dim},  the function $a\mapsto \tilde J_a f$
     is a  diffeomorphism in a small neighbourhood
     of the generic point $a_0$ onto an open subset of $\mathfrak R$.
     Since $x\mapsto a_0\op x$ is  a diffeomorphism (between open subsets of $G$)
     in a neighbourhood of $0$, we get that $a\mapsto \tilde J_{a_0\op a} k$ is a diffeomorphism near $0$ between
     an open subset of $G$ and $\mathfrak R$.
     Since $\tilde J_{a_0}k$ is a constant matrix, it follows that $a\mapsto J_0(M_a)$ is a diffeomorphism near $0$.
     Thus, from the special form of $J_0(M_a)$ and the definition of $d$ we get that $J_0 d$ is invertible.

It follows from the definition of the differential that
$$\lim\limits_{t\to 0} \left(\frac{d(x(t))}{||x(t)||}-\frac{J_0(d)\cdot x(t)}{||x(t)||}\right)=0.$$ Since
$J_0(d)$ is invertible, the limit of $\frac{J_0(d)\cdot x(t)}{||x(t)||}$ is a
non-zero vector, and hence $\lim\limits_{t\to 0} \frac{||d(x(t))||}{||x(t)||}\neq
0$.
\end{proof}

\begin{corollary} \label{fact3}
    Let $e(a)$ and $\delta(a)$ be as in Proposition \ref{IFT}. Assume that $f:G\to G$ is a smooth non $G$-linear $\CD$-function such that
    $f(0)=0$ and $J_0f=0$. Let $g=f\circ f$. Then there is an open neighborhood $U\ni 0$,
    such that for all non-zero $a\in U$, we have

    (i) $||g(a)||<\delta(a)$.

    (ii) There exists a unique $y\in B(0;e(a))$ such that $M(a,y)=g(a).$

\end{corollary}
\begin{proof} Assume that (i) fails.
    Then there exists a definable function $x(t)$
    tending to $0$ in $G$, such that for all $t$,
    \[||g(x(t))||\geq \delta(x(t))=||d(x(t))||^2/32C^2.\]

    Because $J_0f=0$, Fact \ref{square} implies that $\lim\limits_{t\to 0}||g(x(t))||/||x(t)||^2=0$.
     Combined with the above inequality we get $\lim\limits_{x(t)\to 0} ||d(x(t))||^2/||x(t)||^2=0$,
     hence $\lim\limits_{t\to 0} ||d(x(t))||/||x(t)||=0$, contradicting Lemma \ref{limit1}. Thus
     there exists $U\ni 0$ such that for all $a\in U$ we have $||g(a)||<\delta (a)$.
     It now follows from our choice of $\delta(a)$ that there exists a unique $y\in
     B(0;e(a))$ such that $M(a,y)=g(a)$.
\end{proof}

\begin{corollary}\label{the-function} Let $f$ and $g$ be as above, and $e(a)$, $\delta(a)$
as in \emph{\textbf{(C)}} above.  Let
    $U=\{x:|g(x)|<\delta(a) \}$ and $U^*=U\setminus \{0\}$.
    For every $x\in U^*$, let $h(x)$ be the unique $y$ in $B_{e(a)}$
    such that  $M(x,y)=g(x).$ Then

    (i) $U$ contains an open disc around $0$.

    (ii) $h$ is differentiable on $U^*$ and $\lim\limits_{x\to 0}h(x)=0$ (so it extends continuously to
    $0$). Moreover, if $g$ is not constant, then neither is $h$.

    (iii) The continuous extension of $h$ to $U$ is a $\CD$-function.

\end{corollary}
\begin{proof}
    Clause (i) is just  Corollary \ref{fact3}.  To see that $h$ is
    differentiable everywhere we apply the Implicit Function Theorem to $M(x,y)-g(x)$.
     By Lemma \ref{invertible2}, $J_yM_x$ is invertible for every $x\in U^*$ and $|y|<e(a)$,
      so indeed $h(x)$, the solution to $M(x,y)-g(x)=0$,  is differentiable at $x$.

    To see that the limit of $h$ at $0$ is $0$, we compute the limit along an arbitrary curve $x(t)$ tending to $0$.
    By definition, $|h(x(t))|<||e(t)||\sim ||d(t)||$, so since $d(0)=0$,
    also $h(x(t))$ must tend to $0$. The second clause of (ii) follows since if $h$
    were constant with $h(0)=0$ necessarily $h$ would vanish on its domain, implying that $g$ was  identically $0$ (because $M(x,0)=0$ for all $x$). Because $f$ is not constant its image is infinite, and because it is a $\CD$-function it follows that also $g=f\circ f$ is non-constant, a contradiction.

For (iii), note that the graph of $h$ is contained in the
    plane curve  $B=\{(x,y):(x,y,g(x)) \in \hat M(x,y)\}$ where $\hat M$ is the
     $\CD$-definable set of Morley rank $2$ containing the graph of  $M$.
\end{proof}
We  note that locally near the point $(0,0)$ itself, the $\CD$-definable set $B$
need not be the graph of a function, but this does not come up in the argument.
\\

\begin{proof}[Proof of Theorem \ref{main-noninject}]  Assume that
 $f\in \mathfrak F$  and $J_0(f)=0$. We will show that $f$ is not injective near $0$.

  Consider $g(x)=f(f(x))$, and assume towards a contradiction that $f$ and thus also $g$ is injective.
near $0$. By Corollary \ref{the-function}, there exists a $\CD$-function $h$ in a
neighborhood $U$ of $0$, with $h(0)=0$ such that for all $x\in U$,
$$M(x,h(x))=g(x).$$

% By \ref{ }, $h$ is an open map at $0$. Assume that $h$ is $k$-to-one near $0$, with $k$
%possibly greater than $1$.

%By definition for each $b$ near $0$ we have $M(x-b,h(x-b))=g(x-b)$, hence the function
%$h(x-b)$ is $\CD$-definable around every point different than $-b$. Let us see that
%$h_{ab}=M(x-a,h(x-b))$ satisfies the assumptions in \ref{fact4}. The function
%$h_{ab}$ is locally $\CD$-definable around each $x\neq -b$ hence it can have only have
%finitely many non-injective points, and at each regular point $x\neq -b$ the
%Jacobian has positive determinant.

We now wish to apply Lemma \ref{fact4} to the functions $x\mapsto x$ and $x\mapsto h(x)$. For that we
just need to note that for $a$ and $b$ near $0$ the function $g_{a,b}(x)=M(x\om
a,h(x)\om b)$ is non-constant near $0$. Indeed, we can find a fixed definably
connected open $W\ni 0$, such that for $a,b$ close to $0$, $W\sub dom(g_{a,b})$.
Since each $g_{a,b}$ is a $\CD$-function, its graph is contained in a strongly
minimal set, and hence if it were constant near $0$, then it would have to be
constant on the whole of $W$. But then, by the continuity of $M$, the function
$g=g_{0,0}$ must also be constant on $W$, contradiction.

 By applying Lemma \ref{fact4}, we conclude that $g(x)$ is at least $1+k$-to-one near
$0$, where $k\geq 1$. This contradicts the assumption that $f$ and thus $g$ were
locally injective.  Contradiction.
\end{proof}

The following example shows that the proof of Theorem \ref{main-noninject} uses more
than just the basic geometric properties of the function $f$.
\begin{example}
A crucial point in our above argument was that $f(z)/z$, or in the language of our
proof, the implicitly defined function $h$, is an open map. This followed from the
fact that it was a $\CD$-function.

Consider the function $f(z)=|z|^2 z$ from $\mathbb C$ to $\mathbb C$. The function
is smooth everywhere,  $J_0f=0$, and yet it is injective everywhere. However, the
function $f(z)/|z|$ is clearly not an open map.

\end{example}

\subsection{Intersection theory in families}
Based on the topological properties we established thus far we can develop some
intersection theory resembling that of complex analytic curves.

\begin{defn}
Let $X,Y$ be two plane curves, and $p=(p_1,p_2)\in X\cap Y$. We say that {\em $X$ and
$Y$ are tangent at $p$} if there are $\CD$-functions $f,g$ which are $C^1$ in a
neighborhood of $p_1$, with $\Gamma_f\sub X$ and $\Gamma_g\sub Y$, such that
$$f(p_1)=p_2=g(p_1)\,\,\ \mbox{ and } \,\, J_{p_1} f=J_{p_1}g.$$
\end{defn}

The following proposition is the key technical tool for identifying tangency in the reduct $\CD$. The first part of the proposition uses mainly the topological properties of $\CD$-functions to show that if $X,E$ are $\CD$-plane curves intersecting generically enough, then the number of intersection points cannot drop under slight perturbations of the curves. The second part of the proposition uses the differential properties of $\CD$-functions (and in particular Theorem \ref{main-noninject}) to show that if $X$ and $E$ are tangent at a point, then the number of intersection points is expected to increase under slight perturbations.

\begin{proposition} \label{intersection} Let $\CF=\{E_a:a\in T\}$ be a
$\CD$-definable almost faithful family of plane curves, $\CD$-definable over
$\emptyset$, and let $X$ be a strongly minimal plane curve not almost a straight line.

Assume that $a$
 is generic in $T$ over $\emptyset$, $E_a$ strongly minimal,
 $X\cap E_a$ is finite and $p=(x_0,y_0)\in E_a\cap X$.

 \begin{enumerate}
\item If  $p$ is $\CD$-generic in $E_a$ over $a$, non-isolated on $E_a$,
    non-isolated in $X$ and
also $\CD$-generic in $X$ over $[X]$, then for every neighborhood $U\ni p$, there is
a neighborhood $V\ni a$ in $T$, such that for every $a'\in V$, $E_{a'}$ intersects
$X$ in  $U$.

\item (Here we do not make any genericity assumptions on $p$). Assume that for some
open $W\ni a$ whenever $a'\in W$ the set $E_{a'}$ represents a $\CD$-function
$f_{a'}$ in a neighborhood of $(x_0,y_0)$ and that the map $(a',x)\mapsto f_{a'}(x)$ is
continuous at $(a,x_0)$.
 Assume also that $X$
represents a function $g$ at $p$ and that $J_{x_0}(f_a)=J_{x_0}(g)$.
 Then for every neighborhood $U\ni p$ there is a neighborhood $V\ni a$ in $T$,
such that for every $a'\in V$, either $E_{a'}$ and $X$ are tangent at some point in
$U$ or $|E_{a'}\cap X\cap U|>1.$
\end{enumerate}
\end{proposition}
\proof (1) Fix an open $U=U_1\times U_2 \ni p$ definably connected. Since $p$ is non-isolated and
$\CD$-generic in $E_a$ over $a$, it follows from Corollary \ref{open-rel2}, applied to $E_a$, that  there  are three possibilities: (i) $E_a$ is locally at $p$ the graph of a
constant function in the first variable,  (ii) $E_a$ is locally at $p$ the graph of a constant function in the second variable, or (iii) $E_a$ is locally at $p$ the graph of a
homeomorphism.

In all cases, $E_a$ is locally at $p$ either the graph of a continuous function from $x_0$ to $y_0$ or vice versa.
Since our assumptions on $p$ are symmetric with respect to the coordinates, we may assume
that there is an open  $U=U_1\times U_2\ni p$ so that $E_a$ is locally the graph of a
continuous function $f_a:U_1\to U_2$.

 Since, in addition, $a$ is generic in $T$ over $\emptyset$ we may shrink
 $U$ and find an open definably connected $V_0\ni a$ in $T$ such that for every $a'\in V_0$, the set
  $E_{a'}\cap U_1\times U_2$
 is the graph of a $\CD$-function $f_{a'}:U_1\to U_2$ and furthermore, the map $(a',x)\mapsto f_{a'}(x)$
 is continuous on $V_0\times U_1$.

Since $p$ is not isolated in $X$, $\CD$-generic in $X$ over $[X]$, and the
projections of $X$ on both coordinates are finite-to-one, it follows from Corollary
\ref{open-rel2}, applied to $X$, that, after possibly shrinking $U$ further, the set $X\cap U$ is the
graph of an open continuous map $g:U_1\to U_2$.

%There are two cases to consider:
%\\

%\noindent{\bd Case 1} The function $f_a$ is constant near $x_0$.
%\\

%In this case, the genericity of $a$ in $T$ over $x_0$, implies, after possibly
%shrinking $V$, that for all $a'\in V$, the function $f_{a'}$ is constant near $x_0$,
%whose value near $x_0$ is $c(a')$.

% Because $g$ is an open map the image of $U_1$ under $g$ is an open set, so by
% taking $V$ sufficiently small we can ensure that $c(a')\in Im(g)$ for all $a\in V$.
% Thus implies
% that $E_a\cap X\cap U\neq\0$.
% \\

% \noindent{\bf Case 2} The function $f_a$ is not constant near $x_0$.

 Notice that for every $a'\in V_0$, and $(x,y)\in U$,
$$(x,y)\in E_{a'}\cap X\Leftrightarrow f_{a'}(x)\ominus g(x)=0.$$
Because $X\cap E_a$ is finite, the function $f_a\ominus g$ is not constant on its
domain, so by Theorem \ref{C:continuous}, $f_a\om g$ is open on $U_1$.
\\

\noindent{\bf Claim.} {\em There exists $V\ni a$ such that for every $a'\in V\setminus
\{a\}$, the function $f_{a'}\ominus g$ is an open map on $U_1$.}
\begin{proof}[Proof of Claim]
Indeed, assume towards a contradiction that for  $a'\in V_0$ arbitrarily close to
$a$ the map $f_{a'}\ominus g$ is not an open map. Thus, by Theorem
\ref{C:continuous}, it is constant on $U_1$. It follows from continuity that $f_a\om
g$ is constant on $U_1$, contradicting out assumption.
\end{proof}

Thus, we showed that there exists $V\ni a$ such that for all $a'\in V$, the function
$f_{a'}\ominus g$ is open and finite-to-one on $U_1$. In addition, the the map
$(a',x)\mapsto (f_{a'}\om g)(x)$ is continuous in a neighborhood $(a,x_0)$. Because
$0\in (f_a\ominus g)(U)$ it follows from Fact \ref{degree}(1), (4) that for some
open $V_0\ni a$ small enough and for all $a'\in V_0$, the set $(f_{a'}\ominus g)(U)$
contains $0$, namely $X\cap E_{a'}\cap U\neq \emptyset$. This ends the proof of (1).

(2) Let $g$ be a $\CD$-function with $g(x_0)=y_0$ such that $\Gamma_g\sub X$, and
$J_{x_0}f_a=J_{x_0}g$. Note that $(f_a\ominus g)(x_0)=0$ and $f_a\neq g$. So, for
$C\sub G$ a sufficiently small circle around $x_0$  the only zero
 of $f_a\ominus g$
in the closed ball $B$ determined by $C$ is $x_0$. By continuity of $(x,a')\mapsto
f_{a'}(x)$, we may find some neighborhood $V\sub W$ of $a$ such that for every
$a'\in V$, $0\notin (f_{a'}\ominus g)(C)$. It follows from Fact \ref{degree}(1)
that
$$W_C(f_{a'}\ominus g,0)=W_{C}(f_{a}\ominus g,0)$$ for every $a'\in V$.

By our assumptions, $J_{x_0}(f_a\ominus g)=0$ and therefore by Theorem
\ref{main-noninject}, $f_a\ominus g$ is not injective in any neighborhood of $x_0$,
 that is, for every generic $y$ near $0$, $|(f_a\ominus g)^{-1}(y)|>1$. It follows from Fact \ref{degree}(4) that $W_C(f_a\ominus g,0)>1$. Thus, for every $a'\in V$,
 $W_C(f_{a'}\ominus
g,0)>1$.

We can now conclude that  for every $a'$, either $0$ is a regular value of the
function $f_{a'}\ominus g$ on $\intr(C)$, in which case it has more than one
pre-image and then $E_{a'}$ and $X$ intersect more than once in $\intr(C)$, or $0$
is a singular value, in which case the curves $E_{a'}$ and $X$ are tangent at some
point in $\intr(C)$.\qed

\section{The main theorem}

\label{sec-last}

We are now ready to prove our main result. Our proof follows  that of
\cite[Theorem 7.3]{HaKo}.
% As some of the technicalities of that proof were dealt with in earlier stages of the present paper, the proof will be
%  somewhat simpler.
  We begin with a series of useful technical facts. Throughout this section we let $K:=\mfR$.

%We identify $\mathfrak R$ with the ring $K:=R(\sqrt{-1})$. For the purposes of the present
% discussion  we introduce, in analogy with the notion of an open function:
%\begin{definition}
%    A set $R\subseteq G^2$ is \emph{an open relation} if for all open $V\sub G$, the sets
%    $\{x: (\exists v\in V)R(x,v)\}$ and
%    $\{y: (\exists v\in V)R(v,y)\}$ are open in $G$.
%\end{definition}

%Recall that by Corollary \ref{C:open} any $\CD$-definable plane curve, viewed as a binary relation, is -- up to finitely
%many corrections -- open in the above sense. Obviously, if a strongly minimal $C\subseteq G^2$ is open then so is
%$C(x\op a, y\op b)$ for all $a,b\in G$.

\begin{lemma}\label{L:nice families}
    There exist $\CD$-definable families $\CC_0=\{E^0_a:a\in T_0\}$, $\CC_1=\{E^1_b:b\in T_1\}$,
    of plane curves all passing through $(0,0)$ such that:
    \begin{enumerate}
    \item For $i=0,\!1$,  $T_i$ is strongly minimal and  $\CC_i$ is  almost faithful.
        \item Every generic curve in $\CC_i$, $i=0,\! 1$, is closed,  strongly minimal and has no isolated points.
\item There are definable open neighborhoods $U\sub G$ of $0$, and there are
definable
        open sets  $T_0'\sub T_0$, $T_1'\sub T_1$ such that for every $i=0,\! 1$, and $a\in T_i'$,
        the curve $E^i_a$
        represents a function $f^i_t:U\to G$ in $\mathfrak F$,
\item For $i=0,\! 1$, the sets
$$W_i:=\{J_0\, f_a^i:a\in T_i'\}$$ are open subsets of $K$,  with $0\in \cl(W_0)$ and $1\in \cl(W_1)$.
         \item For each $i=0,\! 1$, the map $(a,x)\mapsto f^i_a(x)$ is continuous on  $T_i'\times U$.
    \end{enumerate}
\end{lemma}
\begin{proof} By Claim \ref{last claim}, there exists a $\CD$-function $f:U\to G$
which is not $G$-affine, such that $J_0f=0$.  Let $S\sub G^2$ be a strongly minimal
set representing $f$. By Theorem \ref{main-frontier}, we may assume that $S$ is
closed, and by allowing parameters we may assume that $S$ has no isolated points.
Let
$$\mathcal C_0=\{S\om p:p\in S\}.$$
Let $T_0:=S$ and for $a\in S$ let $E_a^0:=S\om a$.

For every $a=(x_0,f(x_0))\in S$, the curve $S_a$ represents the $\CD$-function
$f(x\op x_0)\om f(x_0)$. By Proposition
       \ref{thm-dim},  the set of elements of $K$ $$W=\{J_0(f(x\oplus x_0)\ominus f(x_0)): x_0\in U\}$$
       has dimension $2$, and by applying the same proposition to a smaller $U$, we see that
       $J_0f=0$ is in the closure of a $2$-dimensional component of $W$. By
       o-minimality, we can find an open $U'\sub U$ such that the
       set $W_0=\{J_0(f(x\oplus x_0)\ominus f(x_0)): x_0\in U'\}$ is an open subset of
       $K$ with the $0$ matrix in its closure. We let  $T_0':=\{(x_0,f(x_0)):x_0\in
       U'\}$. By its definition, the sets $U'$, $T_0'$ and  $W_0$  satisfy all clauses of the lemma.

       In order to obtain $\CC_1$, we replace $f$ with the function $h(x)=f(x)\op x$. It
       is a $\CD$-function which is not $G$-affine, with $J_0h=1\in K$. We repeat
       the above process and obtain the rest of  the lemma.\end{proof}

Our aim is to construct a field configuration in $\CD$ (see Definition \ref{D:gpconf}). We will pull a field configuration from $K$ into $\CD$ by using the properties of Jacobians of $\CD$-functions as studied in the previous sections. Lemma \ref{L:nice families} provides us with the families of curves we will be using to construct the field configuration. For simplicity of notation we will absorb into the language all the parameters needed to define all the objects appearing in Lemma \ref{L:nice families}.

Observe that although $W_0$ and $W_1$ from Lemma \ref{L:nice families}, are not
neighborhoods of $0$ and $1$, respectively, it is still the case that  for every
$B\in W_0$, if $A\in W_1$ and $C\in W_0$ are sufficiently close to $1$ and $0$,
respectively, then $AB+C$ is still in $W_0$ (since $W_0$ is open). Similarly, for
every $A\in W_1$, if $C\in W_1$ is sufficiently close to $1$, then $AC\in W_1$.

Let $e=(1,0)$ be the identity of
 $\mathbb G_m \ltimes \mathbb G_a$, and  choose $b\in W_0$ and
 $h,g$ in $W_1\times W_0\sub \mathbb G_m\ltimes \mathbb G_a$
 sufficiently  close to $e$ so that $gh\in W_1\times W_0$, and $h\cdot b$ and
 $hg\cdot b$ are in $W_0$.
  Note that we may choose $g,h,b$ to be independent generics in the sense of $\CM$
  (and thus also
   independent in the sense of $K$).

To simplify notation,  we denote the functions in $\CC_0$ by $f_t$ and the functions
in $\CC_1$ by $g_s$, and abusing notation, we will sometimes write   $f\in \CC_i$  for a $\CD$-function $f$ which
is represented by a curve in $\CC_i$.  In particular, let us denote, for $i=1, 2$,
$$\CC_i'=\{f_t^i: t\in T_i'\}.$$

We will construct a field configuration of jacobian matrices of $\CD$-functions
in $\CC_0'$ and $\CC_1'$, and show that the corresponding configuration of parameters of $\CD$-definable curves representing those $\CD$-functions  is a field configuration in $\CD$.

%Identifying (as sets) $\mathbb G_m \ltimes \mathbb G_a$ with $K^\times \times K$ we write $g=(a_1,b_1)$ and
%$h=(a_2,b_2)$, $a_i, b_i\in K$.
% In this notation the identity, $e$, of $\mathbb G_m \ltimes \mathbb G_a$ is $(1,0)$.
%  We may choose $g$ and $h$ in $W_1\times W_0$.

We get the following corollary to Lemma \ref{L:nice families}.

 \begin{corollary}  There are $a_1, a_2\in W_1\sub K$ and $b, b_1, b_2\in W_0\sub K$, such that $g=(a_1,b_1),h=(a_2,b_2)\in W_1\times W_0$ and the following hold:
    \begin{enumerate}
        \item There exist $g_1,g_2\in \CC_1'$ and $f_1,f_2,k_1\in \CC_0'$
        with $J_0 g_i=a_i$
        (for $i=1,2$) and $J_0 f_i=b_i$ (for $i=1,2$) and $J_0 k_1=b$.
        \item $hg\in W_1\times W_0$, and there are $f_3\in \CC_0'$ and $g_3\in \CC_1'$ with $(J_0 g_3, J_0 f_3)=hg$.
        \item There are $k_2,k_3\in \CC_0'$ such
        that $J_0 k_2= h\cdot b$ and $J_0 k_3 =hg \cdot b$.
    \end{enumerate}
 \end{corollary}
%\begin{proof}
%    Clause (1) is  Lemma \ref{L:nice families}(4).
%    Clause (2) follows from the fact that $\mathbb G_m \ltimes \mathbb G_a$
%    is a topological group and (3) is again just Lemma \ref{L:nice families}(4).
%     Finally, (5) is the the same as (2) and (3) using the continuity of the action of the group on $K$,
 %     rather than the continuity of group multiplication.
%\end{proof}
%  So our choice of $g,h$ and $b$ implies
% (using (3) of Lemma \ref{L:nice families})   We may choose  $g,h$ as close as we wish to the neutral element, $e$, of
%  $\mathbb G_m \ltimes \mathbb G_a$, such that there are $f_3\in \CC_0$ and $g_3\in \CC_1$ with $(J_0 f_3, J_0 g_3)=hg$ still in $W_1\times W_0$. Similarly, we can find $k_2,k_2\in
%\CC_0$ such
%    that $J_0 k_2= h\cdot b$ and $J_0 k_3 =hg \cdot b$ (throughout, these equalities should be understood via our identification
%     of $\mfR$ with $K$).

For a $\CD$-function $\Psi$, we denote by $[\Psi]$ the $\CD$-canonical parameter of
some fixed strongly minimal set  representing it. Our goal is to prove the following
proposition.

%For a $\CD$-function $\Psi$, let $[\Psi]$ denote the $\CD$-canonical parameter of the strongly  minimal set  representing it.

\begin{proposition}\label{L:fieldConf}
    Keeping the above notation,
    \[\tag{*}
    \mathcal Y:=\{([f_1], [g_1]), ([f_2],[g_2]), ([f_3], [g_3]),  [k_1], [k_2], [k_3] \}
    \]
    is a field configuration in $\CD$.
\end{proposition}
\proof
We have to verify that the following diagram satisfies (1) - (4) from Definition \ref{D:gpconf}:
%    We have to verify the following diagram (in the sense of $\CD$):
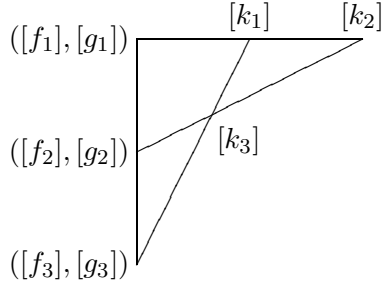
    \[
    \groupconflarge{([f_1], [g_1])}{([f_2],[g_2])}{([f_3], [g_3])}{[k_1]}{[k_2]} {[k_3]}
    \]
       The families $\CC_0$ and $\CC_1$ are almost faithful,
    so for a function $f_a\in \CC_0'$
    we have $\acl_{\CD}(a)=\acl_{\CD}([f_a])$. Since field configurations are stable with respect to $\CD$-inter-algebraicity (over $\0$), we may
    assume that $a=[f_a]$. The same is true for $\CC_1'$.

    By construction, $\mathcal Y$ satisfies the assumptions of Lemma \ref{L:Mcon2Dcon}. So we are reduced to proving (3) of Definition \ref{D:gpconf}.
%Note that for every $z\in \CC_i'$, $J_0 z$ is interdefinable (in $\CM$) with $[z]$.
%Therefore, if $(J_0 z_1, J_0 z_2)$ is an independent tuple (in $\CM$), then the
%tuple $([z_1], [z_2])$ is also independent (in $\CM$), and hence also
%$\CD$-independent. For example,  since the above tuple $(a_1, b_1)$ is independent,
%we obtain that $\mr( [f_1], [g_1])=2$. In other words, the above diagram is vertex by vertex inter-algebraic in $\CM$ with the field configuration $\CF$. Hence all the independence requirements (clauses (1), (2) and (4) of Definition \ref{D:gpconf}) still hold in our diagram and a fortiori also in $\CD$.
%
%Similarly,    we can see that  clauses (1),
%(2) and (4) of Definition \ref{D:gpconf} are true.
     %({\bf WHY $4$? -- because those are inter-definable
     % in the full structure with generic elements of $AGL_1(K)$ and for those this is standard.
     % So (4) is satisfied even in the full structure.})
%      of Definition \ref{D:gpconf} are true.
%       (since they were true in the sense of $\CM$).
%\noindent{\bf I still feel quizzy about (4), but I left it as it was orginally}.
%So in order
%      to verify the diagram, we have to show that the dependencies implied by (3) of that definition remain valid in $\CD$.
       That is, we have to show that all lines in the above diagram represent $\CD$-dependencies. For example, we have to show that $\{[k_2], [k_3],  ([f_2],[g_2])\}$ is $\CD$-dependent, and similarly
       $\{[k_1], [k_3], ([f_3],[g_3])\}$, etc. Since all the arguments are similar, we only prove in detail the latter case.

    It will suffice to show the following statement.
%    We will show that $[k_3]\in \acl_{\CD}([k_1],([f_3],[g_3]))$. Note that as
%     $$(J_0 f_3,J_0 g_3)\cdot J_0 k_1=J_0 f_3 J_0 k_1+J_0 g_3$$ ($\cdot$ on the left is the action of $\mathbb G_m \ltimes \mathbb G_a$ on
%      $\mathbb G_a$ and addition and multiplication on the right are to be understood in the sense of $K$) this element of $K$
%      is exactly  $J_0 (f_3\circ k_1 \op g_3)$. By the choice of $k_3$ we get that $J_0 k_3 = J_0 (f_3\circ k_1 \op g_3)$.
%       So the
%        functions $k_3 $ and $f_3\circ k_1 \op g_3$ are tangent at $0$. Thus, in order to conclude the proof of the proposition
%        we need to show that the tangency
%         of $k_3$ and  $f_3\circ k_1 \op g_3$ implies that $k_3\in \acl_{\CD}([f_3,],[k_1], [g_3]).$
%The proof is given in the next
%          lemma and the geometric idea behind it goes back to Eugenia Rabinovich's work, \cite{Ra}.

\begin{lemma}\label{P:tangency}
$k_3\in \acl_{\CD}([f_3],[k_1],[g_3])$.
\end{lemma}

%   Fix $\CC_0=\{E_b:b\in G\}$ as above and assume it is $\emptyset$-definable. Let $X$ be a $\CD$-definable curve through $(0,0)$
%    which is an open relation and a closed subset of $G^2$, such that $\dim([X]/\emptyset)=2 \mr([X]/\emptyset)$.
%Assume that $J_0 X\in U$  is generic in $K$ over $\emptyset$.
%       If $f_a\in \CC_0$ and $J_0 (f_a)=J_0 X$ then $a\in \acl_{\CD}([X])$.
\proof The geometric idea behind it goes back to Eugenia Rabinovich's work
\cite{Ra}. Write $f_a=k_3$, with $a\in T_0$. By our assumptions, $J_0 f_a$ is
generic in $K$ over $\emptyset$.  To simplify the notation, we denote the curves in
$\CC_0$ by $E_{a'}$, $a'\in T_0$, and the curves in $\CC_1$ by $C_g$, $g\in T_1$.
% (so we identify the function representing $C_{g_3}$ with $g_3$).

Let $X$ be a strongly minimal subset of $S:=(E_{f_3}\circ E_{k_1})\boxplus C_{g_3},$
representing the function $(f_3\circ k_1)\oplus g_3$ (see Lemma \ref{the ring} for
the notation). We want to show that $a\in \acl_{\CD}([S])$. Assume towards a
contradiction that this is not the case.
\\

\noindent{\bf Claim 1.} {\em The projections of $X$ on both coordinates is infinite
and all isolated points of $X$ are in $\acl_\CD([S])$ .}

\begin{proof}[Proof of Claim 1] By our choice of $\CC_0$, the curves $E_{f_3}$ and $E_{k_1}$
 are strongly minimal without isolated points.
It follows that each of these curves has a finite intersection with every straight line, and thus $E_{f_3}\circ E_{k_1}$ has no isolated
points. Indeed, if $(a,b)\in E_{f_3}\circ E_{k_1}$ there is some $c$ such that $(a,c)\in E_{k_1}$ and $(c,b)\in E_{f_3}$. Since both curves are not straight lines and have no isolated points, they are open over $c$ at $(a,c)$ and $(c,b)$,  respectively (Corollary \ref{open-rel2}). So for every open $U\ni (a,c)$  there is $c'\in \pi_2(U\cap E_{k_1})$ distinct from $c$. So there is some $a'$ such that $(a',c')\in U\cap E_{k_1}$. A similar argument will provide us with some $(c',b')\in E_{f_3} $ so $(a',b')\in E_{f_3}\circ E_{k_1}$ with $(a',b')$ arbitrarily close to $(a,b)$.

Note also that the curve $C_{g_3}$ has no isolated points. An argument similar to the one in the previous paragraph shows that the $\boxplus$-sum $S$ of $E_{f_3}\circ E_{k_1}$ and $C_{g_3}$ has no isolated points either. Let $I(X)$ be the set of isolated points of $X$. Let $X':=X\setminus I(X)$. Then, as $S$ has no isolated points, $I(X)\subseteq \cl(S\setminus I(X))$. But $\cl(S\setminus I(X))=\cl(S\setminus X)\cup \cl(X')$, and since $I(X)\cap \cl(X')=\0$ we get that $I(X)\subseteq  \fr(S\sm X)$. Since $[X]\in
\acl_\CD(S)$, we get from Theorem \ref{main-frontier} that $\fr(S\setminus X)\subseteq
\acl_\CD([S])$.

Since $X$ is a strongly minimal set representing a function $(f_3\circ k_1)\oplus
g_3$, its projection on the first coordinate is finite-to-one. Since the function has
a non-zero jacobian at $0$, it is non-constant and hence has an infinite projection
in the second coordinate as well.\end{proof}

%To see that the projection of $X$ on the first coordinate is finite-to-one, we  note that this follows from the
%fact that each of the curves $E_{f_3}$, $E_{k_1}$ and $C_{g_3}$ has this property.
%\\

It follows from Claim 1 that the (finite) set of isolated points of $X$ is contained
in a $\CD$-algebraic set $\CD$-definable over $[S]$. Thus, by removing this
$\CD$-definable set there is no loss of generality in assuming that $X$ contains no
isolated points.

Note that the assumption that $a\notin \acl_\CD([S])$ implies that  $a\notin
\acl_{\CD}([X])$. We will ultimately show that this leads to a contradiction.
Because $E_a$ is strongly minimal, it follows that  $E_a\cap X$ is finite.

    %By the assumption that $\dim\{J_0 f: f\in T_0'\}=2$, we see that
    % not all curves in $\CC_0$ are tangent to
    % $X$ at $0$.

    Since $T_0$ is strongly minimal, there exists some natural number $n$ such that $|X\cap E_b|=n$ for all
       $b\in T_0$ which is $\CD$-generic over $[X]$.
       Thus, the set $$F=\{b\in T_0:|X\cap E_b|<n\}$$ is finite, defined in $\CD$  over $[X]$.
      We will show that $a\in F$, thus reaching a contradiction.

 By our choice of $a$, $\dim (J_0 f_a/\emptyset)=2=\dim G$ and since $J_0 f_a \in \dcl(a)$,
    we also have $\dim(a/\emptyset)=2$. Thus we also have $a\in \dcl(J_0 f_a)$.\\

%{\bf I am still puzzled as to why we need the next claim. It complicates things.}

\noindent{\bf Claim 2.} \label{anotherclaim} \emph{Let $\{x_1,\dots, x_k\}:=X\cap
E_a$. Then for
 every $i=1,\ldots, k$, either  $\mr(x_i/a)=1$, or $x_i\in \acl_{\CD}(\emptyset)$.}
\begin{proof}[Proof of Claim 2] We consider the family
$$\CF'=\{(E_{a_1}\circ E_{a_2})\boxplus C_b:a_1,a_1\in T_0\,\, ,\,\, b\in T_1\},$$ and for simplicity
write the members of $\CF'$ as $\{X_t:t\in T\}$.  By our choice of $X$, there is
$t_0\in T$ generic such that $X$ is a strongly minimal subset of $X_{t_0}$, so
definable over $\acl_{\CD}(t_0)$. We may now replace $\CF'$ by another family of the
same dimension, defined over $\0$, such that the generic member of $\CF'$ is
strongly minimal and $X$ belongs to the family. We
call this new family $\CF$.

Thus $X=X_{t_0}$, with $\CF=\{X_t:t\in T\}$ a $\CD$-definable almost faithful family
of plane curves, and $t_0$ generic in $T$ over $\0$. Our underlying negation assumptions implies that $\mr(a/t_0)=1$.

Assume now that $\mr(x_i/a)\neq 1$. Since $x_i\in E_a$ it follows that  $x_i\in \acl_{\CD}(a)$.
Because $\mr(a/t_0)=1$ it follows that $t_0$ is $\CD$-generic in $T$ over $a$ and hence also over $x_i$.
But then $x_i$ is in $X_t$ for every $t$ which is $\CD$-generic in $T$. This necessarily implies that $x_i\in \acl_{\CD}(\0)$ becausre there can be only finitely many points in
$G\times G$ belonging to every $\CD$-generic curve $X_t$.
 This ends the proof of Claim 2.\end{proof}

We now return to the proof of Lemma \ref{P:tangency}. By Claim 2 we may  assume that
for $i=1,\ldots, r$,
 we have $\mr(x_i/a)=1$ and for $i=r+1,\ldots, k$, we have
$x_i\in \acl_{\CD}(\emptyset)$. Without loss of generality, $x_{k}=0$.

In order to show that $a\in F$,     we have to show that $k<n$. Towards that end, we
will show that there are infinitely many  $a'\in T_0$ such
    that $n=|X\cap E_{a'}|\ge k+1$.

    Let $U_1\dots U_r, U_k$ be pairwise disjoint open neighborhoods of $x_1,\ldots,x_r, x_{k}$, respectively.
Since $x_{r+1}, \ldots , x_k$ are in $\acl_{\CD}(\emptyset)$, each of these
points belongs to all but finitely many $E_{a'}$.

Because $X$ and $E_a$ have no isolated points, we can apply Proposition \ref{intersection}. We
first apply Proposition \ref{intersection} (2) to $0=x_k$, and obtain $V\ni a$ such
that for every $a'\in V$, $|E_{a'}\cap X\cap U_k|\geq 2$, counted with multiplicity.
Because $J_0(E_a)$ is generic in $K$, it is attained at most finitely many times and
hence by choosing $V$ sufficiently small and $a'\in V$, $a'\neq a$, the
 curves $E_{a'}$ and $X$ are not tangent at $0$, so there exists $p\in E_{a'}\cap X\cap U_k$ which is different than $0$.
It follows that for all but finitely many $a'\in V$, $|E_{a'}\cap X\cap U_k|\geq 2$.

We now  apply Proposition \ref{intersection}(1) to $x_1,\ldots, x_r$, and  obtain a
sub-neighborhood $V'$ of $a$ such that for every $a'\in V'$, and $i=1,\ldots, r$,
$E_{a'}\cap X\cap U_i\neq \emptyset$.

Summarizing, we see  that for every $a'\neq a$ close to $a$, we have  $|E_{a'}\cap
X|\geq k+1$, and therefore $a$ is in the finite set $F$ which we defined above. This
ends the proof of Lemma
 \ref{P:tangency}, and with it the proof of Proposition \ref{L:fieldConf}. \qed \\

We can now prove our main result.
\begin{theorem}
    Let $\CD=\la G;\oplus,\cdots\ra $ be a strongly minimal expansion of a group $G$,
    interpretable in an o-minimal expansion $\CM$ of a field $R$, with $\dim_{\CM}(G)=2$. If $\CD$ is not
     locally modular, then there exists in $\CD$ an interpretable algebraically closed field $K\simeq R((\sqrt{-1})$,
     and there exists a $K$-algebraic group $H$, such that $G$ and $H$ are definably isomorphic in $\CD$ and
     every $\CD$-definable subset of $H^n$ is
     $K$-constructible.

     Moreover, the structure $\CD$ and the field $K$ are bi-interpretable.
\end{theorem}
\begin{proof}
    By  Proposition   \ref{L:fieldConf},
      the configuration $\mathcal Y$ of $(*)$ is a field configuration in $\CD$.

     By Fact \ref{Hrushovski}, an algebraically
     closed field $K$ is interpretable in $\CD$. By strong minimality, there exists a $\CD$-definable function
     $f:G\to K$ with finite fibres (this is standard using the symmetric functions on $K$).
     By \cite[Lemma 4.6]{PePiSt}
     (and using strong minimality of $G$), there exists a  finite subgroup $F\le G$ such that $G/F$ is
     internal to $K$ in the structure $\CD$  (it is in fact, the proof of the
     lemma which provides us with the finite subgroup $F$).
By \cite[Theorem 3.1]{PeStExpansions}, every $\CD$-definable subset of $K^n$ is
$K$-constructible, and therefore $G/F$ is $\CD$-definably isomorphic  to a
$K$-constructible group. By Weil-Hrushovski, \cite[Theorem 1]{Bous-WeilHr}, it is therefore definably isomorphic to
a $K$-algebraic group $H$ (of algebraic dimension $1$). It is known  that $H$, as an algebraic curve with
all its induced $K$-algebraic structure, is bi-interpretable with $K$ (this follows, for example, from the main result of \cite{HrZil}). For the sake of completeness let us sketch this argument.

If $C$ is an algebraic curve in $K$, then clearly, $C$ is interpretable in $K$. Since, up to finitely many points, $C$ is affine, it is inter-algebraic in $K$ with $K$. Using this inter-algebraicity, we can pull-back any field configuration from $K$ to $C$ allowing us to interpret a field $K'$ in $C$ (with its $K$-induced structure). By \cite[Theorem 4.15]{PoiGroups}, $K'$ is $K$-definably isomorphic to $K$, so in particular $K$ is interpretable in $C$. Finally, the isomorphism from $K$ to $K'$ takes $C$ to a $C$-interpretable curve $C'$. The induced map from $C$ to $C'$ is $K$-definable hence it is definable in $C$. This shows that $C$ and $K$ are bi-interpretable.

So, $H$ and hence also $G/F$, with all its induced $\CD$-structure is bi-interpretable with $K$.
By Lemma \ref{G and G/F}, the structure $\CD$ is  also bi-interpretable
with $K$.
\end{proof}

\subsection{Concluding remarks}\label{sec-concl}

Note that as a result of the main theorem, the almost $K$-structure on $G$ which we
introduced in Section \ref{almost-complex} turns out to be definably isomorphic to
the $K$-structure of the algebraic group $H$. Thus, in this very special setting, we
are able to mimic the classical result about the integrability of $2$-dimensional
almost complex curves.

Also, note that the general o-minimal version of Zilber's conjecture remains open
for general strongly minimal structures whose universe has dimension $2$. As we
noted earlier, the more general conjecture, allowing underlying sets of arbitrary
dimension is open even for reducts of the complex field.

\bigskip
\footnotesize
\noindent\textit{Acknowledgments.}
The first author wishes to thank Rahim Moosa for
many enlightening discussions on the subject, as well as the model theory groups at
Waterloo and McMaster for running a joint working seminar on the relevant literature
during the academic year 2011-2012. The authors wish to thank Sergei Starchenko for
his important feedback, and Michael Wan for pointing out an earlier mistake in Section \ref{sec-frontier}. Many thanks also to the Oberwolfach Mathematical Institute for
bringing the authors together during the Workshop in Model Theory in 2016, and to
the Institute Henri Poincare in Paris, for its hospitality during the trimester
program ``Model theory, Combinatorics and valued fields''  in 2018. Finally, we thank the referee for a very careful reading of the manuscript and for providing us with numerous comments that have contributed significantly to the presentation of this paper.

The first author was supported by a Research Grant from the German Research Foundation (DFG) and a Zukunftskolleg Research Fellowship. The second author was
partially supported by an Israel Science  Foundation grants No. 1156/10 and 181/16.

\end{document}